\documentclass[10pt,twoside]{article}

\usepackage[latin1]{inputenc}
\usepackage{amsmath, amssymb, amsthm, amsfonts, amsgen, stmaryrd}
\usepackage{indentfirst}
\usepackage{fancyhdr}
\usepackage{graphicx}

\numberwithin{equation}{section}

\def\rightharpoonupfill@{\arrowfill@\relbar\relbar\rightharpoonup}
\newcommand{\xrightharpoonup}[2][]{\ext@arrow 0359\rightharpoonupfill@{#1}{#2}}

\newcommand{\ds}{\displaystyle}

\newcommand{\Nb}{{\mathbb{N}}}
\newcommand{\Qb}{{\mathbb{Q}}}
\newcommand{\Rb}{{\mathbb{R}}}
\newcommand{\Zb}{{\mathbb{Z}}}

\newcommand{\lb}{{\llbracket}}
\newcommand{\rb}{{\rrbracket}}

\def\leq{\leqslant}
\def\geq{\geqslant}

\let\e=\varepsilon

\newcommand{\W}{{\mathcal{W}}}
\newcommand{\We}{{\mathcal{W}_{\e}}}

\let\O=\Omega
\let\o=\omega
\let\G=\Gamma
\let\a=\alpha
\let\go=\rightarrow

\newtheorem{thm}{Theorem}[section]

\newtheorem{rmk}[thm]{Remark}
\newtheorem{lemma}[thm]{Lemma}
\newtheorem{proposition}[thm]{Proposition}
\newtheorem{coro}[thm]{Corollary}
\makeatletter

\def\debaixodolim#1#2{\mathrel{}\mathop{\lim}\limits^{#1}_{#2}}

\def\debaixodoliminf#1#2{\mathrel{}\mathop{\liminf}\limits^{#1}_{#2}}
\def\debaixodolimsup#1#2{\mathrel{}\mathop{\limsup}\limits^{#1}_{#2}}
\def\debaixodosup#1#2{\mathrel{}\mathop{\sup}\limits^{#1}_{#2}}
\def\debaixodauniao#1#2{\mathrel{}\mathop{\bigcup}\limits^{#1}_{#2}}

\begin{document}

\title{\bf Multiscale nonconvex relaxation and \\ application to thin films}
\author{ Jean-Fran\c cois Babadjian and Margarida Ba\'{\i}a \footnote{Corresponding
author}}
\date{}

\maketitle

\chead[{\sc Jean-Fran\c cois Babadjian and Margarida Ba\'{\i}a}]
{\sc Multiscale nonconvex relaxation and application to thin films}
\thispagestyle{empty}

\begin{abstract}
\noindent $\G$-convergence techniques are used to give a
characterization of the behavior of a family of heterogeneous
multiple scale integral functionals. Periodicity, standard growth
conditions and nonconvexity are assumed whereas a stronger uniform
continuity with respect to the macroscopic variable, normally
required in the existing literature, is avoided.  An application
to dimension reduction problems in reiterated homogenization of
thin films is presented.

\vspace{0.3cm}


\noindent {\bf MSC 2000 (${\cal A}{\cal M}{\cal S}$):} 35E99,
35M10, 49J45,  74B20,  74G65, 74Q05,  74K35.
\end{abstract}


\section{Introduction}

\noindent In this work we study the $\e$-limit behavior of an
elastic body whose microstructure is periodic of  period $\e$ and
$\e^2$, and whose volume may also depend on  this small parameter
$\e$, by a $\G$-limit argument. We refer to the books of Dal Maso
\cite{DM}, Braides \cite{B} and Braides and Defranceschi \cite{BD}
for a comprehensive treatment on this kind of variational
convergence.  We seek to approximate in a $\G$-convergence sense the
microscopic behavior of such materials by a macroscopic, or average,
description. The asymptotic analysis of media with multiple scale of
homogenization is referred to  as {\it Reiterated Homogenization}.

Let ${\O}_\e$ denote the reference configuration of this elastic
body that we assume to be a bounded and open subset of $\Rb^{N}$
($N\geq 1$). In the sequel we identify  $\Rb^{d \times N}$ (resp.
$\Qb^{d \times N}$) with  the space of real (resp. rational)-valued
$d \times N$ matrices and  $Q$ will stand for the unit cube
$(0,1)^N$ of $\Rb^N$.

To take into account the periodic heterogeneity of this material,
we suppose that its stored energy density,  given  by a function
$f: \O_\e \times \Rb^N \times \Rb^N \times \Rb^{d \times N} \to
\Rb$, is $Q$-periodic with respect to  its second and third
variables, and we treat the nonconvex case under  standard growth
and coercivity conditions of order $p$, with $1<p<\infty$.   Under
a deformation $u: \O_\e \to \Rb^d$ the elastic energy of this body
turns out to be given  by the functional

\begin{equation}\label{1027}
\int_{\O_\e}f\left(x,\frac{x}{\e},\frac{x}{\e^2};\nabla
u(x)\right)dx.
\end{equation}
\noindent Its dependence upon the small parameter $\e$ allow us to
consider materials whose microscopic heterogeneity scales like $\e$
and $\e^2$. The generalization  of this study to any number of
scales $k\geq 2$  follows  by an iterated argument similar to the
one used in Braides and Defranceschi (see Remark 22.8 in \cite{BD}).

\begin{figure}[h]
\begin{center}
\scalebox{.20}{\includegraphics{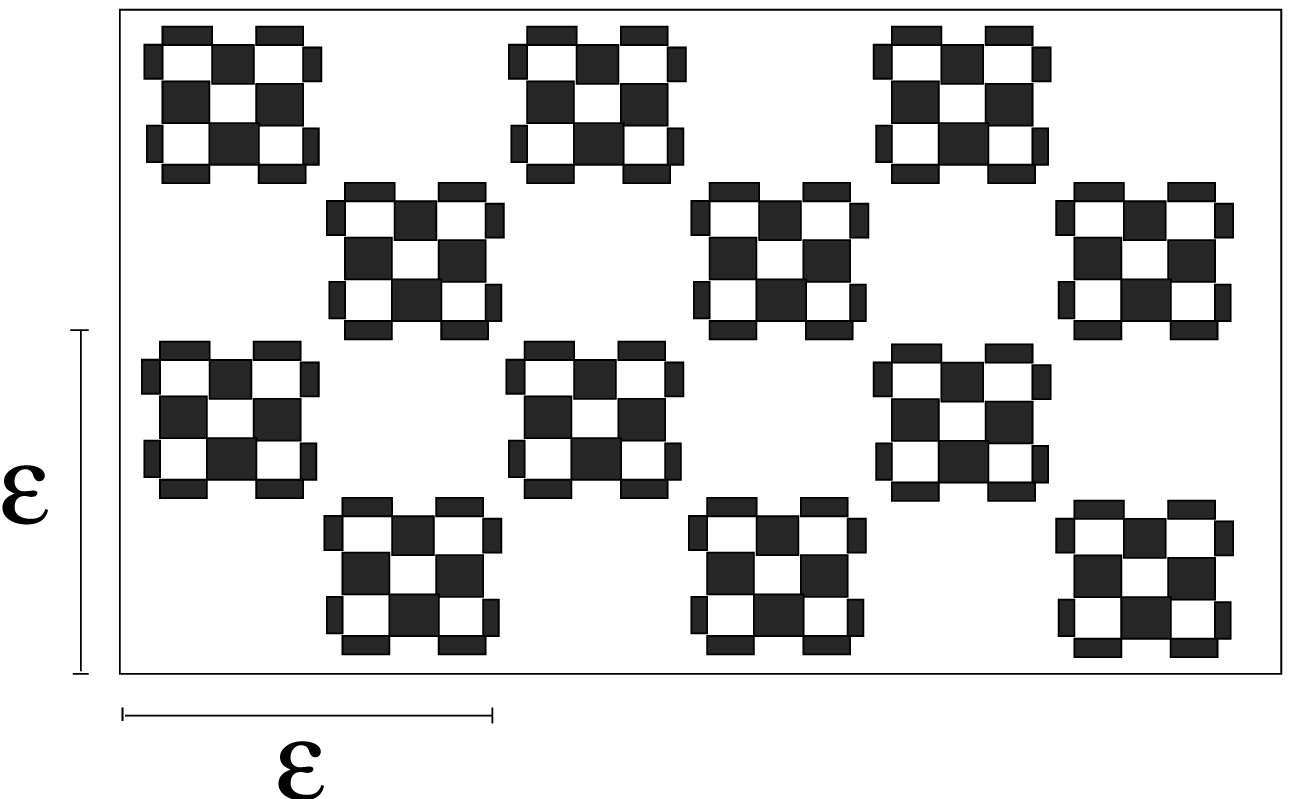}} \hspace{1cm}
\scalebox{.75}{\includegraphics{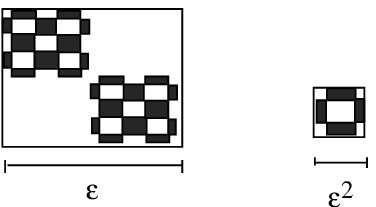}}
\end{center}
\caption{Example of a domain with periodic microstructure of order
$\e$ and $\e^{2}$}\label{phc}
\end{figure}

We will address  two independent problems: In Section \ref{part2} we
will analyze the behavior of bodies with periodic microstructure
whose volume does not depend on $\e$, yielding to a pure reiterated
homogenization problem (Figure \ref{phc}). The originality of this
part is that we do not require any uniform continuity hypotheses on
the first and second variables of $f$, as it is customary in the
existing literature (see references below).

In Section \ref{part3} we consider three dimensional cylindrical
bodies with similar periodic properties as before, whose thickness
scales like $\e$, leading to a homogenization and dimension
reduction problem. The main contribution here is that our arguments
allow us to homogenize this  material  in the reducing  direction.\\

More precisely, in Section \ref{part2} we  describe the case where
$\O_{\e}=\O$  and our family of energies  is of the form
\begin{equation}\label{intfunc}
\int_{\O}f\left(x,\frac{x}{\e},\frac{x}{\e^2};\nabla u(x)
\right)\, dx.
\end{equation}
\noindent This kind of asymptotic problems  can be seen as a
generalization of the {\it Iterated Homogenization Theorem} for
linear integrands, proved by Bensoussan, Lions and Papanicolau
\cite{BLP},  in which  the homogenized operator is derived by  a
formal two-scale asymptotic expansion method. This result has been
recovered in several ways via other types  of convergence such as
$H$-convergence, $G$-convergence, and multiscale convergence (see
Lions, Lukkassen, Persson and  Wall \cite{LLPW1,LLPW2} and
references therein). In the framework of $\G$-convergence, Braides
and Lukkassen (see Theorem 1.1 in \cite{BL}  and also \cite{L})
investigated the nonlinear setting for integral functionals of the
type
$$\int_{\Omega}f\left(\frac{x}{\e},\frac{x}{\e^2};\nabla u(x) \right)\, dx,$$
where the integrand $f : \Rb^N \times \Rb^N \times \Rb^{d \times N}
\to \Rb$ is assumed to satisfy usual periodicity and growth
conditions and
\begin{itemize}
\item[-] $f(y,\cdot\,;\xi)$ is measurable for all $(y,\xi) \in
\Rb^N \times \Rb^{d \times N}$;
\item[-] $f(y,z;\cdot)$ is continuous and convex for all $(y,z)
\in \Rb^N \times \Rb^N$;
\item[-] there exist a locally integrable function $b$   and a
continuous positive real function $\omega$, with $\omega(0)=0$, such
that
\begin{equation}\label{unif-cont} |f(y,z;\xi)-f(y',z;\xi)| \leq
\omega(|y-y'|)\big[b(z)+f(y,z;\xi)\big]\end{equation} \noindent for
all $y$, $y'$, $z \in \Rb^N$, $\xi \in \Rb^{d \times N}$.
\end{itemize}
\noindent This result has been  extended to the case of nonconvex
integrands depending explicitly on the macroscopic variable $x$, as
in (\ref{intfunc}), under the above strong uniform continuity
condition (\ref{unif-cont}) (see Theorem 22.1 and Remark 22.8 of
Braides and Defranceschi \cite{BD}). Using techniques of multiscale
convergence and restricting the argument to the convex and
homogeneous case (no dependence on the variable $x$), Fonseca and
Zappale were able to weaken the continuity condition
(\ref{unif-cont}). Namely, they only required $f$ to be continuous
(see Theorem 1.9 in \cite{FZ}). Recently, in an independent work,
Barchiesi \cite{Bar} studied the $\G$-limit of functionals of the
type (\ref{intfunc}) in the convex case under  weak regularity
assumptions on $f$ with respect to the oscillating variables.

All the above results share the same property: The homogenized
functional, as it is referred in the literature for the
$\Gamma$-limit of (\ref{intfunc}), is obtained by iterating twice
the homogenization formula derived in the study of the $\G$-limit
of functionals of the type
\begin{equation}\label{1-micro-s}
\int_{\Omega} f\left(x,\frac{x}{\e};\nabla
 u(x)\right)\, dx.\end{equation}

\noindent This can be seen from formula (14.12) in Braides and
Defranceschi \cite{BD}, or in Ba\'{\i}a and Fonseca \cite{B&Fo}
where $f: \Omega \times\Rb^{N}\times\Rb^{d\times N} \rightarrow \Rb$
is assumed to satisfy standard $p$-growth and $p$-coercivity
conditions, $Q$-periodicity with respect to the oscillating
variable, and
\begin{itemize}
\item[-] $f(x,\cdot\,;\cdot)$ is continuous for a.e. $ x\in
\Omega$; \item[-] $f(\cdot,y;\xi)$ is measurable for all $y\in
\Rb^{N}$ and all $\xi \in \Rb^{d\times N}$
\end{itemize}
\noindent (see Theorem 1.1 in \cite{B&Fo}; it will be of use for the
proofs of Theorems \ref{jf-m1} and \ref{jf-m3} below).

As these previous  results seem to show, it is not clear what is the
natural regularity on $f$  for  the integral (\ref{intfunc}) to be
well defined. Such problems have been discussed by Allaire (Section
5 in \cite{A}). In particular,  the measurability of the function $x
\mapsto f(x,x/\e,x/\e^2;\xi)$ is ensured  whenever $f$ is continuous
in its second and third variable. Following the lines of Ba\'{\i}a
and Fonseca \cite{B&Fo} we will assume that
\begin{itemize}
\item[$(H_1)$] $f(x,\,\cdot\,,\,\cdot\,;\,\cdot\,)$ is continuous
for a.e. $x \in \Omega$;
\item[$(H_2)$] $f(\,\cdot\,,y,z;\xi)$ is measurable for all $(y,z,\xi) \in \Rb^N \times \Rb^N \times \Rb^{d\times N}$;
\item[$(H_3)$]  $f(x,\,\cdot\,,z;\xi)$ is
$Q$-periodic for all $(z,\xi)\in \Rb^{N}\times \Rb^{d\times N}$ and
a.e. $x \in \Omega$, and $f(x,y,\,\cdot\,;\xi)$ is $Q$-periodic for
all $(y,\xi)\in \Rb^{N }\times\Rb^{d\times N}$ and a.e. $x \in
\Omega$;
\item[$(H_4)$] there exists $\beta>0$ such that
$$ \frac{1}{\beta}|\xi|^p-\beta\leq f(x,y,z;\xi)\leq
\beta(1+|\xi|^p) \quad \text{for  all  }(y,z,\xi) \in \Rb^{N} \times
\Rb^N \times \Rb^{d\times N} \text{ and a.e. }x \in \Omega.$$
\end{itemize}
\noindent From the applications point of view it would be
interesting to consider  functions that are continuous with respect
to the first variable and only measurable with respect to some of
the oscillating variables, as it is relevant, for instance, in the
case of mixtures. Nevertheless, the arguments we use here do not
allow us to treat this case.

In what follows we write $\Gamma(L^p(\Omega))$-limit whenever we
refer to the $\G$-convergence with respect to the usual metric in
$L^{p}(\O;\Rb^{d})$. The above considerations lead us to the main
result of Section \ref{part2} that in  particular recovers Theorem
1.9 in Fonseca and Zappale \cite{FZ}.

\begin{thm}\label{jf-m1}  Let $f: \Omega \times \Rb^{N
}\times\Rb^{N}\times\Rb^{d\times N} \rightarrow \Rb$ be a function
satisfying $(H_1)$-$(H_4)$. For each $\e>0$ define ${\cal F}_{\e}:
L^{p}(\Omega;\Rb^{d})\rightarrow \overline \Rb$ by
\begin{equation}\label{Iepsilon}
{\cal F}_{\e}(u):=\left\{\begin{array}{ll} \ds \int_{\Omega}
f\left(x,\frac{x}{\e},\frac{x}{\e^{2}};\nabla u(x)\right)\, dx &
\text{if } u\in W^{1,p}(\Omega;\Rb^{d}),\\ +\infty &
\text{otherwise}.
\end{array}\right.
\end{equation}
\noindent Then the  $\Gamma (L^p(\Omega))$-limit of the family
$\{\mathcal F_\e\}_{\e>0}$ is given by the functional
$$\mathcal F_{\rm hom}(u): =
 \left\{\begin{array}{ll} \ds \int_{\Omega} \overline f_{\rm hom}(x;\nabla u(x))\,
 dx & \text{if } u\in W^{1,p}(\Omega;\Rb^{d}), \\
+ \infty &\text{otherwise},
\end{array}\right.$$
\noindent where $\overline f_{\rm hom}$ is defined for all $\xi\in
\Rb^{d\times N}$ and a.e. $x \in \Omega$ by
\begin{equation}\label{fhombarra}
\overline f_{\rm hom}(x;\xi):=\lim_{T \to  + \infty} \inf_{\phi}
\left\{ \frac{1}{T^{N}}\int_{(0,T)^{N}}  f_{\rm hom}(x,y;\xi+\nabla
\phi(y))\, dy: \; \phi \in W_{0}^{1,p}((0,T)^{N}
;\Rb^{d})\right\},\end{equation} \noindent and
\begin{equation}\label{fhom}
f_{\rm hom}(x,y;\xi):=\lim_{T \to  + \infty} \inf_{\phi} \left\{
\frac{1}{T^{N}}\int_{(0,T)^{N}}  f(x,y,z;\xi+\nabla \phi(z))\, dz:\;
\phi \in W_{0}^{1,p}((0,T)^{N} ;\Rb^{d})\right\}\end{equation}
\noindent  for a.e. $x\in \O$ and  all $(y,\xi)\in \Rb^N \times
\Rb^{d\times N}$.
\end{thm}
\noindent As mentioned before one homogenizes first with respect to
$z$, considering $y$ as a parameter, and then one homogenizes with
respect to $y$.

We remark that most of the proofs presented in this section follow
the lines of the ones in Braides and Defranceschi \cite{BD} (Theorem
22.1 and Remark 22.8), and that our main contribution is to use
arguments that  allow us to weaken the strong uniform continuity
hypothesis (\ref{unif-cont}). Let us briefly describe how we
proceed: The idea consists in proving the result for integrands
which do not depend explicitly on $x$ (see Theorem \ref{jf-m2} in
Section \ref{2homo}), and then to treat the general case  by
freezing this macroscopic variable (Section \ref{2hetero}). To do
this, we start by claiming that under hypotheses
$(H_{1})$-$(H_{4})$, $f$ is uniformly continuous up to a small
error. Indeed, since $f$ is a Carath\'eodory integrand,
Scorza-Dragoni's Theorem (see Ekeland and Temam \cite{Ek&Te})
implies that the restriction of $f$ to $K \times \Rb^N \times \Rb^N
\times \Rb^{d \times N}$ is continuous, for some compact set $K
\subset \O$ whose complementary has arbitrarily small Lebesgue
measure. Then the periodicity of $f$ with respect to its second and
third variable leads  $f$ to  be uniformly continuous on $K \times
\Rb^N \times \Rb^N \times \overline B$, for some closed ball
$\overline B$ of $\Rb^{d \times N}$ of sufficiently large radius.
Finally, to ensure that the energy remains arbitrarily small on the
complementary of $K$ and on the set of $x$'s such that the gradient
of the deformation does not belong to $\overline B$, we use the
Decomposition Lemma (see Fonseca, M\"uller and Pedregal \cite{FMP}
or  Fonseca and Leoni \cite{FL}) which allows us to select
minimizing sequences with $p$-equi-integrable gradient. Thus, in
view of the $p$-growth character of the integrand, the energy over
sets of arbitrarily small Lebesgue measure tends to zero.\\

\begin{figure}
\begin{center}
\scalebox{.60}{\includegraphics{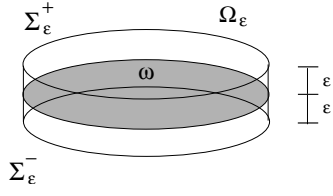}}
\end{center}
\caption{Cylindrical thin domain of thickness $\e$}\label{td}
\end{figure}

In Section \ref{part3} we consider the case where $\O_{\e}$ is a
cylindrical thin domain of the form $\O_{\e}:=\o\times (-\e, \e)$
(Figure. \ref{td}), whose heterogeneity may depend periodically
upon its thickness. We assume that its basis, $\o$, is a bounded
open subset of\, $\Rb^2$ and we seek to characterize  the behavior
of the elastic energy (\ref{1027}) when $\e$ tends to zero.

\begin{figure}[h]
\begin{center}
\scalebox{.50}{\includegraphics{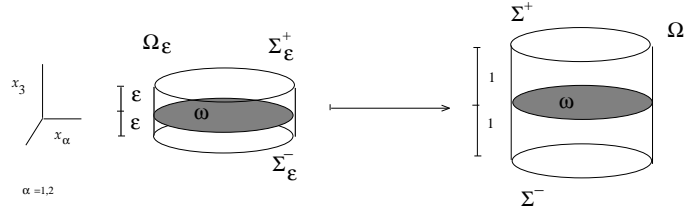}}
\end{center}
\caption{Rescaled domain of unit thickness}\label{rd}
\end{figure}

Two  simultaneous features occur in this case: a reiterated
homogenization and a  dimension reduction process. As usual, in
order to study this problem as $\e\go 0$ we rescale the
$\varepsilon$-thin body into a reference domain of unit thickness
(see e.g.\! Anzellotti, Baldo and Percivale \cite{ABP}, Le Dret
and Raoult \cite{LDR}, and Braides, Fonseca and Francfort
\cite{BFF}), so that the resulting energy will be defined on a
fixed body, while the dependence on $\varepsilon$ turns out to be
explicit in the transverse derivative. For this, we consider the
change of variables
$$ \Omega_{\varepsilon}\rightarrow \Omega:=\omega\times I, \quad (x_1, x_2, x_3)
\mapsto  \left(x_1, x_2, \frac{1}{\varepsilon}x_3\right),$$
\noindent  and define $v(x_\alpha,x_3 / \e)=u(x_\alpha,x_3)$ on the
rescaled cylinder $\O$, where $I:=(-1,1)$ and $x_\alpha:=(x_1,x_2)$
is the in-plane variable (Figure. \ref{rd}).

In what follows we denote by $\nabla_i={\partial}/{\partial x_i}$
for $i \in \{1,2,3\}$ and $\nabla_\alpha= (\nabla_1,\nabla_2)$.  For
all $\overline{\xi}=(z_1|z_2)\in \Rb^{3\times 2}$ and $z\in \Rb^3$,
$(\overline{\xi}|z)$ is the matrix whose first two columns are $z_1$
and $z_2$ and whose last one is $z$. After replacing $v$ by $u$ in
(\ref{1027}), changing variables and dividing by $\e$, our goal is
to study  the sequence of rescaled energies
\begin{equation}\label{1005}
\int_{\O} W\left(x,\frac{x}{\e},\frac{x_\a}{\e^2};\nabla_\a
v(x)\Big| \frac{1}{\e} \nabla_3 v(x) \right)dx,
\end{equation}
\noindent where  to simplify notations we set
\begin{equation}\label{1845}
f(x,y,z;\xi)=W(x_\a,y_3,y_\a,z_3,z_\a;\xi)
\end{equation}
\noindent  for $W: \O \times \Rb^3 \times \Rb^2 \times \Rb^{3 \times
3} \to \Rb$  (or equivalently
$W(x',y',z'_\a;\xi)=f(x'_\a,x_3,y'_\a,x'_3,z'_\a,y'_3;\xi)$).

Similar results have been already studied independently by Shu
\cite{S} (Theorem 3 (i)) and by Braides, Fonseca and Francfort
\cite{BFF} (Theorem 4.2) for energies of the type
$$\int_\O W\left(x_3,\frac{x_\a}{\e};\nabla_\a v \Big|\frac{1}{\e} \nabla_3 v \right) dx$$
\noindent where  $W : I \times \Rb^2 \times \Rb^{3 \times 3} \to
\Rb$ is a Carath\'eodory integrand.  Later Babadjian and Ba\'{\i}a
\cite{BB1} (Theorem 1.2) treated the case where the integrand $W$
depends also on  the macroscopic in-plane variable $x_\a$ under
measurability  hypothesis with respect to $x=(x_\a,x_3)$, and a
continuity  requirement with respect to the oscillating variable.
Integral functionals of the form
$$\int_\O W\left(\frac{x_\a}{\e^2};\nabla_\a v \Big|\frac{1}{\e} \nabla_3 v \right) dx,$$
\noindent have been studied in Shu \cite{S} (Theorem 5) under
different length scales for the film thickness and the material
microstructure. Homogenization in the transverse direction $x_3$
remained an open question and comes as a consequence of the main
result of this section.

In the sequel we denote by $\mathcal L^{N}$ the $N$-dimensional
Lebesgue measure in $\Rb^{N}$, by $Q':=(0,1)^2$ the unit cube in
$\Rb^{2}$ and  we will identify $W^{1,p}(\o;\Rb^3)$ (resp.
$L^p(\o;\Rb^3)$) with those functions $u \in W^{1,p}(\O;\Rb^3)$
(resp. $L^p(\O;\Rb^3)$) such that $\nabla_3 u(x)=0$ a.e. in $\O$.

Following the lines of Babadjian and Ba\'{\i}a \cite{BB1} and
Babadjian and Francfort \cite{Ba&Fr} we assume that
\begin{itemize} \item[$(A_1)$] $W(x,\,\cdot\,,\,\cdot\,;\,\cdot\,)$ is
continuous for a.e. $x \in \Omega$; \item[$(A_2)$]  $W(\,\cdot\,
,y,z_\a;\xi)$ is measurable for all $(y,z_\a,\xi) \in \Rb^3 \times
\Rb^2 \times \Rb^{3 \times 3}$;
 \item[$(A_3)$]
$\left\{\begin{array}{l} y_\a \mapsto W(x,y_\a,y_3,z_\a;\xi) \,\,
{\rm is}\,\,  Q'\text{-periodic\,\, for\,\, all} \,\,
(z_\a,y_3,\xi)\in \Rb^{3} \times \Rb^{3\times 3} \,\,
\text{ and  \,\, a.e. } x \in \Omega,\\\\
(z_\a,y_3) \mapsto W(x,y_\a,y_3,z_\a ;\xi) \,\, {\rm is}\,\, Q
\text{-periodic\,\, for\,\, all} \,\, (y_\a,\xi)\in
\Rb^{2}\times\Rb^{3\times 3} \,\,\text{ and\,\, a.e. } x \in \Omega;
\end{array}\right.$
\item[$(A_4)$] there exists $\beta>0$ such that
$$ \frac{1}{\beta}|\xi|^p-\beta\leq W(x,y,z_\a;\xi)\leq
\beta(1+|\xi|^p) \quad \text{for  all }(y,z_\a,\xi) \in \Rb^{3}
\times \Rb^2 \times \Rb^{3\times 3} \text{ and a.e. }x \in \Omega.$$
\end{itemize}
We prove the following theorem.
\begin{thm}\label{jf-m3}
Let $W: \Omega \times \Rb^{3}\times\Rb^{2}\times\Rb^{3\times 3}
\rightarrow \Rb$ be a function satisfying $(A_1)$-$(A_4)$. For each
$\e>0$, consider the functional $\We :
L^{p}(\Omega;\Rb^{3})\rightarrow \overline \Rb$ defined by
\begin{equation}\label{Wepsilon}
\We(u):=\left\{\begin{array}{ll} \ds \int_{\Omega}
W\left(x,\frac{x}{\e},\frac{x_\a}{\e^{2}};\nabla_\a
u(x)\Big|\frac{1}{\e} \nabla_3 u(x) \right)\, dx & \text{if } u\in
W^{1,p}(\Omega;\Rb^{3}),\\
 +\infty & \text{otherwise}.
\end{array}\right.
\end{equation}
\noindent Then the  $\Gamma (L^p(\Omega))$-limit of the family
$\{\We\}_{\e>0}$ is given by  the functional
$$\mathcal \W_{\rm hom}(u): =
\left\{\begin{array}{ll} \ds 2 \int_{\omega} \overline W_{\rm
hom}(x_\a;\nabla_\a u(x_\a))\, dx_\a & \text{if } u\in
W^{1,p}(\omega;\Rb^{3}), \\ + \infty &\text{otherwise},
\end{array}\right.$$
\noindent where $\overline W_{\rm hom}$ is defined, for all
$\overline \xi\in \Rb^{3\times 2}$ and a.e. $x_\a \in \omega$, by
\begin{eqnarray}\label{Whombarra}
\overline W_{\rm hom}(x_\a;\overline \xi) & := & \lim_{T \to  +
\infty} \inf_{\phi} \Bigg\{ \frac{1}{2 T^{2}}\int_{(0,T)^{2} \times
I} W_{\rm hom}(x_\a,y_3,y_\a;\overline \xi+\nabla_\a \phi(y)|
\nabla_3 \phi(y)
)\, dy: \nonumber\\
&&\hspace{2.0cm}\phi \in W^{1,p}((0,T)^{2} \times I ;\Rb^{3}), \quad
\phi=0 \text{ on }\partial (0,T)^2 \times I\Bigg\}\end{eqnarray}
\noindent and
\begin{eqnarray}\label{Whom}
W_{\rm hom}(x,y_\a;\xi) & := & \lim_{T \to  + \infty} \inf_{\phi}
\Bigg\{ \frac{1}{T^{3}}\int_{(0,T)^{3}}
W(x,y_\a,z_3,z_\a;\xi+\nabla \phi(z))\,
dz:\nonumber\\
& &  \hspace{4.0cm} \phi \in W_{0}^{1,p}((0,T)^{3};\Rb^{3})\Bigg\},
\end{eqnarray}
\noindent  for a.e. $x\in \O$ and  all $(y_\a,\xi)\in \Rb^2 \times
\Rb^{3\times 3}$.
\end{thm}
\vspace{0.2cm}
\begin{rmk}\label{infwhom}
{\rm It can be proved that the limits as $T \to +\infty$ in
(\ref{Whombarra}) and (\ref{Whom}) can
 be replaced by an infimum taken for every $T>0$ as in
Braides and Defranceschi \cite{BD} or Ba\'{\i}a and Fonseca
\cite{B&Fo}. }
\end{rmk}

Let us formally justify the periodicity assumptions $(A_3)$: Since
the volume of  $\O_\e$ is of order $\e$ and $\e^2 \ll \e$, in a
first step, we can think of  $\e$ as being a fixed parameter and let
$\e^{2}$ tend to zero.  Then at this point dimension reduction is
not occurring and (\ref{1027}) can be seen as a single one-scale
homogenization problem as in (\ref{1-micro-s}), in which it is
natural to assume $f(x_\a,y,\cdot\,;\xi)$, or equivalently,
$(z_\a,y_3) \mapsto W(x,y_\a,y_3,z_{\alpha};\xi)$ (see (\ref{1845}))
to be $Q$-periodic.  The homogenization formula for this case  give
us an homogenized stored energy density $W_{\rm
hom}(x,y_{\alpha};\xi)$ that, in a second step, is used as the
integrand of a similar problem than the one treated in Babadjian and
Ba\'{\i}a \cite{BB1}. In particular the required $Q'$- periodicity
of $W_{\rm hom}(x,\cdot\,;\xi)$ can be obtained from the
$Q'$-periodicity of $y_\a \mapsto W(x,y_\a,y_3,z_\a;\xi)$.

We finally observe  that the proof of  Theorem \ref{jf-m3} is very
closed to its $N$-dimensional analogue Theorem \ref{jf-m1},  the
main difference being the  use of the   scaled gradients
decomposition lemma derived by Bocea and Fonseca in place of the
usual Decomposition Lemma (compare Theorem 1.1 of \cite{Bo&Fo} with
Fonseca, M\"uller and Pedregal \cite{FMP} or Fonseca and Leoni
\cite{FL}). As it has been noted in \cite{Ba&Fr,BB1}, Theorem
\ref{jf-m3} and Theorem \ref{jf-m1} cannot be treated similarly. In
the former case we need to extend our Carath\'eodory integrands by a
continuous function by means of Tietze's  Extension Theorem (see
e.g. DiBenedetto \cite{DiB}). This argument was already used in
Babadjian and Ba\'{\i}a \cite{BB1} and Babadjian and Francfort
\cite{Ba&Fr} where the authors used a weaker extension result (see
Theorem 1, Section 1.2 in Evans-Gariepy \cite{EG} and Lemma 4.1 in
\cite{BB1}).

In the sequel, given $\lambda>0$ we denote by
$\overline{B}(0,\lambda)$  the closed ball of radius $\lambda$ in
$\Rb^{d \times N}$, that is the set $\{ \xi \in \Rb^{d\times N}:
|\xi|\leq \lambda\}$, and  the letter  $C$ stands for a generic
constant.  Throughout the text$\debaixodolim {} {n,m}:=
\debaixodolim {} {n} \hspace{-0.1cm} \debaixodolim {} {m}$ while
$\debaixodolim {} {m,n}:= \debaixodolim {} {m} \hspace{-0.1cm}
\debaixodolim {} {n}$ with obvious generalizations.


\section{Reiterated homogenization}\label{part2}

\noindent The main objective of this section is to prove  Theorem
\ref{jf-m1}.   In Subsection \ref{continuity} we state the main
properties of  $f_{\rm hom}$ and $\overline{f}_{\rm hom}$ that are
basic for our analysis.  In Subsection \ref{2homo} we present some
auxiliary results  for  the proof of the homogeneous counterpart
of Theorem \ref{jf-m1}, Theorem \ref{jf-m2}, in which we assume
that $f$ does not depend explicitly on $x$. The proof of Theorem
\ref{jf-m1} in its fully generality is presented in Subsection
\ref{2hetero}.   Finally, in Subsection \ref{cx-subsection}  we
remark an alternative proof for  convex integrands.

\begin{rmk}
{\rm In the sequel, and without loss of generality, we  assume
that $f$ is non negative. Indeed, it suffices to replace $f$ by
$f+\beta$ which is non negative in view of $(H_4)$. }
\end{rmk}

 As for  notations   $Q(a,\delta):=a+\delta
(-1/2,1/2)^N$ (cube of center $a \in \Rb^N$ and edge length
$\delta$) and  $Q:=(0,1)^{N}$ stands for  the unit cube in
$\Rb^N$.


\subsection{Properties of $\overline{f}_{\rm hom}$}\label{continuity}

\noindent Repeating the argument used in Ba\'{\i}a and Fonseca
\cite{B&Fo} (Theorem 1.1), we can see that the function $f_{\rm
hom}$ given in (\ref{fhom}) is well defined and it is (equivalent
to) a Carath\'eodory function:
\begin{eqnarray}
&&f_{\rm hom}(\cdot,\cdot\,;\xi) \text{ is }\mathcal L^N
\otimes \mathcal L^N\text{-measurable for all  }\xi \in \Rb^{d \times N}, \label{fhom1}\\
&&f_{\rm hom}(x,y;\cdot) \text{ is continuous for }\mathcal L^N
\otimes \mathcal L^N\text{-a.e. }(x,y) \in \Omega \times \Rb^N.
\label{fhombis2}
\end{eqnarray}
By condition $(H_3)$ it follows that
\begin{equation}\label{fhom3}
f_{\rm hom}(x,\cdot\,;\xi) \text{ is }Q\text{-periodic for a.e. }x
\in \O \text{ and all }\xi \in \Rb^{d \times N}.
\end{equation}
\noindent Moreover, $f_{\rm hom}$ is  quasiconvex in the $\xi$
variable and satisfies the same $p$-growth and $p$-coercivity
condition $(H_4)$ as $f$:
\begin{equation}\label{fhom4}
\frac{1}{\beta}|\xi|^p-\beta\leq f_{\rm hom}(x,y;\xi) \leq
\beta(1+|\xi|^p)\quad \text{for a.e. }x \in \O\text{ and all
}(y,\xi) \in \Rb^{N} \times \Rb^{d\times N},
\end{equation}
\noindent where $\beta$ is the constant in $(H_{4})$.

As a consequence of (\ref{fhom1}), (\ref{fhombis2}) and
(\ref{fhom4}), the function $\overline f_{\rm hom}$ given in
(\ref{fhombarra}) is also well defined, and is (equivalent to) a
Carath\'eodory function, which implies that the definition of
$\mathcal F_{\rm hom}$ makes sense on $W^{1,p}(\Omega;\Rb^d)$.

Finally, $\overline f_{\rm hom}$ is also quasiconvex in the $\xi$
variable and satisfies the same $p$-growth and $p$-coercivity
condition $(H_4)$ as $f$ and $f_{\rm hom}$:
\begin{equation}\label{fhom5}
\frac{1}{\beta}|\xi|^p-\beta\leq \overline f_{\rm hom}(x;\xi)
\leq \beta(1+|\xi|^p)\quad \text{for a.e. }x \in \O\text{ and all
}\xi \in \Rb^{d\times N},
\end{equation}
\noindent where, as before,  $\beta$ is the constant in $(H_{4})$.


\subsection{Independence of the macroscopic variable}\label{2homo}

\noindent  We assume that $f$  does not depend explicitly on $x$,
namely $f:\Rb^{N}\times \Rb^{N}\times \Rb^{d\times N}\rightarrow
\Rb^{+}$. In addition, according to $(H_{1})$-$(H_{2})$ and unless
we specify the contrary, we assume $f$ to be continuous and to
satisfy hypotheses $(H_{3})$ and $(H_{4})$.

For each $\e>0$ consider the functional ${\cal F}_{\e}:
L^{p}(\O;\Rb^{d})  \rightarrow [0,+\infty]$ defined by

\begin{equation}\label{Iepsilon2}
{\cal F}_{\e}(u):=\left\{\begin{array}{ll} \ds \int_{\O}
f\left(\frac{x}{\e},\frac{x}{\e^{2}};\nabla u(x)\right)\, dx &
\text{if } u  \in W^{1,p}(\O;\Rb^{d}),\\
+\infty & \text{otherwise}.
\end{array}\right.
\end{equation}
\noindent Our objective is to prove the following result.

\begin{thm}\label{jf-m2}
Under assumptions  $(H_3)$ and $(H_4)$ the  $\Gamma (L^p(\O))$-limit
of the family $\{\mathcal F_\e \}_{\e>0}$  is given by
$$\mathcal F_{\rm hom}(u) =
 \left\{\begin{array}{ll} \ds \int_{\O} \overline f_{\rm hom}(\nabla u(x))\,
 dx & \text{if } u  \in W^{1,p}(\O;\Rb^{d}),\\
+ \infty &\text{otherwise},
\end{array}\right.$$
\noindent where $\overline f_{\rm hom}$ is defined by
\begin{equation}\label{fhombarra2}
\overline f_{\rm hom}(\xi):=\lim_{T \to  + \infty} \inf_{\phi}
\left\{ \frac{1}{T^{N}}\int_{(0,T)^{N}}  f_{\rm hom}(y;\xi+\nabla
\phi(y))\, dy: \; \phi \in W_{0}^{1,p}((0,T)^{N}
;\Rb^{d})\right\}\end{equation} \noindent for all $\xi\in
\Rb^{d\times N}$, and where
\begin{equation}\label{fhombarra21}f_{\rm hom}(y;\xi):=\lim_{T \to  + \infty}
\inf_{\phi} \left\{ \frac{1}{T^{N}}\int_{(0,T)^{N}}
f(y,z;\xi+\nabla \phi(z))\, dz: \; \phi \in
W_{0}^{1,p}((0,T)^{N} ;\Rb^{d})\right\}\end{equation} for all
$(y,\xi)\in \Rb^N \times \Rb^{d\times N}$.
\end{thm}

This result can be seen as a generalization of Theorem 1.9 in
Fonseca and Zappale \cite{FZ} (for $s=1$), in which, as it is usual
for the convex case,  it is enough to consider variations that are
periodic in the cell $Q$. Their multiscale argument (see Subsection
\ref{cx-subsection} below) does not apply here since, as it is
expected in the non convex  case, the variations should be
considered to be periodic over an infinite ensemble of cells, as it
is seen from (\ref{fhombarra2}) and (\ref{fhombarra21}).\\

We start the proof of Theorem \ref{jf-m2} by localizing the
functionals given in (\ref{Iepsilon2}) in order to highlight their
dependence on the class of bounded, open subsets of $\Rb^N$,
denoted by $\mathcal A_0$. As it will be clear from the  proofs of
Lemmas \ref{lem<} and \ref{lem>} below, it would not be sufficient
to localize, as usual, on any open subset of $\O$.  Indeed,
formulas (\ref{fhombarra2}) and (\ref{fhombarra21}) suggest to
work in cubes of the type $(0,T)^{N}$, with $T$ arbitrarily large,
not necessarily contained in $\O$.\\

For each $\e>0$ consider ${\cal F}_{\e}: L^{p}(\Rb^N;\Rb^{d})
\times {\cal A}_0 \rightarrow [0,+\infty]$ defined by
\begin{equation}\label{Fe2}
{\cal F}_{\e}(u;A):=\left\{\begin{array}{ll} \ds \int_{A}
f\left(\frac{x}{\e},\frac{x}{\e^{2}};\nabla u(x)\right)\, dx &
\text{if } u \in W^{1,p}(A;\Rb^{d}),\\ +\infty & \text{otherwise}.
\end{array}\right.
\end{equation}

We will prove (Subsections \ref{e-gl} and \ref{c-gl} below) that the
family of functionals $\{\mathcal
F_\varepsilon(\cdot\,;A)\}_{\e>0}$, with  $A\in {\cal A}_{0}$,
$\G$-converges with respect to the strong
$L^{p}(A;\Rb^{d})$-topology to the functional ${\mathcal F}_{\rm
hom} (\cdot\,;A)$, where $\mathcal F_{\rm hom} : L^p(\Rb^{N};\mathbb
R^d) \times \mathcal A_0 \rightarrow [0,+\infty]$ is given by
$$\mathcal F_{\rm hom}(u;A) =
\left\{\begin{array}{ll} \ds \int_{A} \overline f_{\rm hom}(\nabla
u(x))\, dx & \text{if } u  \in W^{1,p}(A;\Rb^{d}),\\ + \infty
&\text{otherwise}.
\end{array}\right.$$
\noindent  As a consequence, taking $A=\O$ yields Theorem
\ref{jf-m2}.


\subsubsection{Existence and integral representation of the $\G$-limit}\label{e-gl}

\noindent Given $\{\e_j\} \searrow 0^+$  and $A \in \mathcal A_0$,
consider the $\Gamma$-lower limit  of $\{  \mathcal F_{\e_j}
(\cdot\,;A) \}_{j \in \Nb}$ for the $L^{p}(A;\Rb^{d})$-topology
defined  for $u \in L^{p}(\Rb^N;\Rb^d)$ by
$$\mathcal F_{\{\e_j\}}(u;A) := \inf_{\{u_j\}}\left\{\liminf_{j \to
+\infty} \mathcal F_{\e_j}(u_j;A) : \quad u_j \to u \text{ in
}L^p(A;\Rb^d)\right\}.$$ \noindent In view of the $p$-coercivity
condition $(H_{4})$  it follows that  $\mathcal F_{\{\e_j\}}(u;A)$
is infinite whenever $u \in L^p(\Rb^N;\Rb^d) \setminus
W^{1,p}(A;\Rb^d)$  for each $A \in \mathcal A_0$, so it suffices
to study the case where $u \in W^{1,p}(A;\Rb^d)$.

  The following result states the existence of $\G$-convergent
subsequences and can be proved by classical arguments on the lines
of those of  Braides, Fonseca and Francfort \cite{BFF}.

\begin{lemma}\label{glA0}
For any sequence $\{\e_j\} \searrow 0^+$, there exists a
subsequence $\{\e_n\} \equiv \{\e_{j_n}\}$ such that $\mathcal
F_{\{\e_n\}}(u;A)$ is the $\G(L^p(A))$-limit of $\{\mathcal
F_{\e_n}(u;A)\}_{n \in \Nb}$ for all $A \in \mathcal A_0$ and all $u
\in W^{1,p}(A;\Rb^d)$.
\end{lemma}


Our goal is to study the behavior of $\mathcal
F_{\{\e_n\}}(u;\cdot)$ in $\mathcal A(A)$ (family of open subsets of
$A$) for each $u \in W^{1,p}(A;\Rb^d)$. Following the proof of Lemma
2.10 in Ba\'{\i}a and Fonseca \cite{B&Fo}, it is possible to show
that $\mathcal F_{\{\e_n\}}(u;\cdot)$ is a measure on $\mathcal
A(A)$ for all $A \in \mathcal A_0$. Namely,  the following result
holds.

\begin{lemma}
For each $A \in \mathcal A_0$ and all $u \in W^{1,p}(A;\Rb^d)$,
the restriction of $\mathcal F_{\{\e_n\}}(u;\cdot)$ to $\mathcal
A(A)$ is a Radon measure, absolutely continuous with respect to
the $N$-dimensional Lebesgue measure.
\end{lemma}

For the moment, we are not in position to apply Buttazzo-Dal Maso
Integral Representation Theorem (see Theorem 1.1 in \cite{DMB})
because, a priori, the integrand would depend on the open set $A
\in \mathcal A_0$. The following result  prevents this dependence
from holding since it leads to an homogeneous integrand as it will
be seen in Lemma \ref{intrep} below.

\begin{lemma}\label{independency}
For all $\xi \in \Rb^{d\times N}$, $y_0$ and $z_0\in \Rb^N$, and
$\delta>0$
$$\mathcal F_{\{\e_{n}\}}(\xi\,\cdot\,;Q(y_0,\delta))=\mathcal F_{\{\e_{n}\}}(\xi\,\cdot\,;
Q(z_0,\delta)).$$
\end{lemma}

\noindent {\it Proof.}  Clearly, it suffices to establish the
inequality
$$\mathcal F_{\{\e_{n}\}}(\xi\,\cdot\,;Q(y_0,\delta))\geq
\mathcal F_{\{\e_{n}\}}(\xi\,\cdot\,;Q(z_0,\delta)).$$ \noindent Let
$\{w_{n}\}\subset W^{1,p}_0(Q(y_0,\delta);\Rb^{d})$, with
$w_{n}\rightarrow 0$ in $L^{p}(Q(y_0,\delta);\Rb^{d})$, be such that
$$\mathcal F_{\{\e_{n}\}}(\xi\cdot\,;Q(y_0,\delta))=\lim_{n \to
+\infty}\int_{Q(y_0,\delta)}f\Big(\frac{x}{\e_{n}},
\frac{x}{\e_{n}^{2}};\xi+\nabla w_{n}(x)\Big)\,dx$$ \noindent  (see
e.g. Proposition 11.7 in Braides and Defranceschi \cite{BD}). By
hypothesis $(H_4)$ and the Poincar\'{e} Inequality, we can suppose
that  the sequence $\{w_n\}$ is uniformly bounded in
$W^{1,p}(Q(y_0,\delta);\Rb^d)$. Thus by the Decomposition Lemma (see
Fonseca, M\"uller and Pedregal \cite{FMP} or Fonseca and Leoni
\cite{FL}), there exists a subsequence of $\{w_n\}$ (still denoted
by $\{w_n\}$) and a sequence $\{u_n\} \subset
W^{1,\infty}_0(Q(y_0,\delta);\Rb^{d})$ such that $u_n\rightharpoonup
0$ in $W^{1,p}(Q(y_0,\delta);\Rb^d)$,
\begin{equation}\label{equi3}
\{ |\nabla u_n|^p  \} \text{ is equi-integrable}
\end{equation}
\noindent and
\begin{equation}\label{equi4}
\mathcal L^N( \{ y\in Q(y_0,\delta):\quad u_n(y) \neq w_n (y) \})
\to 0.
\end{equation}
Then, in view of (\ref{equi3}), (\ref{equi4}) and the
$p$-growth condition $(H_4)$,
\begin{eqnarray}\label{gls}
\mathcal F_{\{\e_{n}\}}(\xi\cdot\,;Q(y_0,\delta))& \geq
&\limsup_{n \to +\infty}\int_{Q(y_0,\delta) \cap
\{u_n=w_n\}}f\Big(\frac{x}{\e_{n}},
\frac{x}{\e_{n}^{2}};\xi+\nabla u_{n}(x)\Big)\,dx\nonumber\\
& \geq & \limsup_{n \to
+\infty}\int_{Q(y_0,\delta)}f\Big(\frac{x}{\e_{n}},
\frac{x}{\e_{n}^{2}};\xi+\nabla u_{n}(x)\Big)\,dx
\end{eqnarray}
\noindent For all $n\in \Nb$ we  write
$$\frac{y_0-z_0}{\e_{n}}=m_{\e_{n}}+s_{\e_{n}}$$
\noindent  with $m_{\e_{n}}\in \Zb^{N}$ and $s_{\e_{n}}\in
[0,1)^{N}$,
\begin{equation}\label{quo-1}\frac{m_{\e_n}}{\e_n}=\theta_{\e_n}+l_{\e_n}\end{equation}
\noindent  with $\theta_{\e_{n}}\in \Zb^{N}$ and $l_{\e_{n}}\in
[0,1)^N$, and we define
\begin{equation}\label{quo-2}
x_{\e_{n}}:=m_{\e_{n}}\e_{n}-\e_{n}^{2}l_{\e_{n}}.\end{equation}
\noindent  Note that $x_{\e_n}=y_0-z_0-\e_n s_{\e_{n}}
-\e_{n}^{2}l_{\e_{n}} \to y_0-z_0$ as $n \to  +\infty$.  For all
$n\in \Nb$, extend $u_n$ by zero to the whole $\Rb^N$ and set
$v_n(x)=u_n(x+x_{\e_n})$ for $x \in Q(z_0,\delta)$. Then $\{v_n\}
\subset W^{1,p}(Q(z_0,\delta);\Rb^d)$,  $v_n \to 0$ in
$L^p(Q(z_0,\delta);\Rb^d)$, 
and  the sequence $\{|\nabla v_n|^p\}$ is equi-integrable by the
translation invariance of the Lebesgue measure. In view of
(\ref{gls}), (\ref{quo-1}), (\ref{quo-2}) and $(H_3)$,
\begin{eqnarray*}
\mathcal F_{\{\e_{n}\}}(\xi\cdot\,;Q(y_0,\delta))&\geq&\limsup_{n
\to
+\infty}\int_{Q(y_0-x_{\e_{n}},\delta)}f\Big(\frac{x+x_{\e_{n}}}{\e_{n}},
\frac{x+x_{\e_{n}}}{\e_{n}^{2}};\xi+\nabla
u_{n}(x+x_{\e_{n}})\Big)\,dx\\
 & = & \limsup_{n \to +\infty}\int_{Q(y_0-x_{\e_{n}},\delta)}f\Big(\frac{x}{\e_{n}}-\e_n l_{\e_n},
\frac{x}{\e_{n}^{2}};\xi+\nabla
v_{n}(x)\Big)\,dx\\
 & \geq & \limsup_{n \to +\infty}\int_{Q(z_0,\delta)}f\Big(\frac{x}{\e_{n}}-\e_n l_{\e_n},
\frac{x}{\e_{n}^{2}};\xi+\nabla
v_{n}(x)\Big)\,dx \\
& & \hspace{1.0cm}-\limsup_{n \to +\infty}\int_{Q(z_0,\delta)
\setminus Q(y_0-x_{\e_n},\delta)}f\Big(\frac{x}{\e_{n}}-\e_n
l_{\e_n}, \frac{x}{\e_{n}^{2}};\xi+\nabla v_{n}(x)\Big)\,dx.
\end{eqnarray*}
\noindent Since $v_n \equiv 0$ outside $Q(y_0-x_{\e_n},\delta)$, the
$p$-growth condition $(H_4)$ and the fact that $\mathcal
L^N(Q(z_0,\delta) \setminus Q(y_0-x_{\e_n},\delta)) \to 0$ yield
\begin{eqnarray*}
&& \limsup_{n \to +\infty}\int_{Q(z_0,\delta) \setminus
Q(y_0-x_{\e_n},\delta)}f\Big(\frac{x}{\e_{n}}-\e_n l_{\e_n},
\frac{x}{\e_{n}^{2}};\xi+\nabla
v_{n}(x)\Big)\,dx \\
&&\hspace{5.0cm}\leq \limsup_{n \to +\infty}\beta
(1+|\xi|^p)\mathcal L^N (Q(z_0,\delta) \setminus
Q(y_0-x_{\e_n},\delta)) =0,
\end{eqnarray*}
\noindent and therefore
\begin{equation}\label{2143}
\mathcal F_{\{\e_{n}\}}(\xi\cdot\,;Q(y_0,\delta)) \geq \limsup_{n
\to +\infty}\int_{Q(z_0,\delta)}f\Big(\frac{x}{\e_{n}}-\e_n
l_{\e_n}, \frac{x}{\e_{n}^{2}};\xi+\nabla v_{n}(x)\Big)\,dx.
\end{equation}
\noindent To eliminate   the term $\e_n \, l_{\e_n}$ in
(\ref{2143}),  and thus to recover $\mathcal
F_{\{\e_{n}\}}(\xi\cdot\,;Q(z_0,\delta))$,  we would like to apply a
uniform continuity argument. Since $f$ is continuous on $\Rb^N
\times \Rb^N \times \Rb^{d \times N}$ and separately $Q$-periodic
with respect to its two first variables, by hypothesis $(H_3)$, then
$f$ is uniformly continuous on $\Rb^N \times \Rb^N \times
\overline{B}(0,\lambda)$ for any $\lambda>0$. We define
$$R^\lambda_n:= \{ x \in Q(z_0,\delta) : |\xi  + \nabla
v_n(x)| \leq \lambda\},$$ \noindent and we note that by  Chebyshev's
inequality
\begin{equation}\label{che-1}\mathcal L^N(Q(z_0,\delta) \setminus R^\lambda_n) \leq
C/\lambda^p,\end{equation} \noindent for some constant $C>0$
independent of $\lambda$ or $n$. Thus, in view  of (\ref{2143}) and
the fact that $f$ is nonnegative,
$$\mathcal F_{\{\e_{n}\}}(\xi\cdot\,;Q(y_0,\delta)) \geq \limsup_{\lambda, n}
\int_{R^\lambda_n}f\Big(\frac{x}{\e_{n}}-\e_n l_{\e_n},
\frac{x}{\e_{n}^{2}};\xi+\nabla v_{n}(x)\Big)\,dx.$$ \noindent
Denoting by $\omega_\lambda : \Rb^+ \to \Rb^+$ the modulus of
continuity of $f$ on $\Rb^N \times \Rb^N \times
\overline{B}(0,\lambda)$, we get that  for any $x \in R^\lambda_n$
$$\left| f\Big(\frac{x}{\e_{n}},
\frac{x}{\e_{n}^{2}};\xi+\nabla
v_{n}(x)\Big)-f\Big(\frac{x}{\e_{n}}-\e_n l_{\e_n},
\frac{x}{\e_{n}^{2}};\xi+\nabla v_{n}(x)\Big)\right| \leq
\omega_\lambda(\e_n l_{\e_n}).$$ Then, the continuity of
$\omega_\lambda$ and the fact that $\omega_\lambda(0)=0$ yield
\begin{eqnarray*}
 \mathcal F_{\{\e_{n}\}}(\xi\cdot\,;Q(y_0,\delta))&\geq&
\limsup_{\lambda, n}\left\{
\int_{R^\lambda_n}f\Big(\frac{x}{\e_{n}},
\frac{x}{\e_{n}^{2}};\xi+\nabla
v_{n}(x)\Big)\,dx - \delta^N \omega_\lambda(\e_n l_{\e_n})\right\}\\
&  = &  \limsup_{\lambda, n}
\int_{R^\lambda_n}f\Big(\frac{x}{\e_{n}},
\frac{x}{\e_{n}^{2}};\xi+\nabla v_{n}(x)\Big)\,dx.
\end{eqnarray*}
\noindent The equi-integrability of $\{|\nabla v_n|^p\}$, the
$p$-growth condition $(H_4)$ and (\ref{che-1}), imply that
\begin{eqnarray*}
&&\limsup_{\lambda, n} \int_{Q(z_0,\delta) \setminus
R^\lambda_n}f\Big(\frac{x}{\e_{n}},
\frac{x}{\e_{n}^{2}};\xi+\nabla
v_{n}(x)\Big)\,dx \\
&&\hspace{2.0cm}\leq \beta \limsup_{\lambda \to +\infty}\sup_{n \in
\Nb}\int_{Q(z_0,\delta) \setminus R^\lambda_n}(1+|\nabla v_n(x)|^p
)dx = 0,
\end{eqnarray*}
\noindent and since $v_n \to 0$ in $L^p(Q(z_0,\delta);\Rb^d)$,
\begin{eqnarray*}
\mathcal F_{\{\e_{n}\}}(\xi\cdot\,;Q(y_0,\delta)) & \geq &
\limsup_{n \to +\infty}
\int_{Q(z_0,\delta)}f\Big(\frac{x}{\e_{n}},
\frac{x}{\e_{n}^{2}};\xi+\nabla
v_{n}(x)\Big)\,dx\\
 & \geq & \mathcal F_{\{\e_{n}\}}(\xi\cdot\,;Q(z_0,\delta)).
\end{eqnarray*}
$\hfill\blacksquare$\\

As a consequence of this lemma, we derive the following result.

\begin{lemma}\label{intrep}
There exists a continuous function $f_{\{\e_n\}}: \Rb^{d \times N} \rightarrow \Rb^+$
such that for all $A \in \mathcal A_0$ and all $u \in W^{1,p}(A;\Rb^d)$,
$$\mathcal F_{\{\e_{n}\}} (u;A)=\int_{A} f_{\{\e_n\}}(\nabla u(x))\, dx.$$
\end{lemma}

\noindent {\it Proof. } Let $A\in \mathcal A_0$. By Buttazzo-Dal
Maso Integral Representation Theorem (Theorem 1.1 in \cite{DMB}),
there exists a Carath\'eodory function $f^A_{\{\e_n\}}: A \times
\Rb^{d \times N} \rightarrow \Rb^+$ satisfying
$$\mathcal F_{\{\e_n\}} (u;U)=\int_{U} f^A_{\{\e_n\}}(x;\nabla u(x))\, dx$$
\noindent for all $U \in \mathcal A(A)$ and all $u \in
W^{1,p}(U;\Rb^d)$. Furthermore, for a.e. $x \in A$ and all $\xi \in
\Rb^{d \times N}$
$$f^A_{\{\e_n\}}(x;\xi)=\lim_{\delta  \to  0}\frac{ \mathcal F_{\{\e_{n}\}}
(\xi \cdot\,;Q(x,\delta))}{\delta^N}.$$ \noindent Define
$f_{\{\e_n\}}: \Rb^{d \times N} \rightarrow \Rb^+$ by
$$f_{\{\e_n\}}(\xi)=\lim_{\delta  \to  0}\frac{ \mathcal F_{\{\e_{n}\}}
(\xi \cdot\,;Q(0,\delta))}{\delta^N}.$$ \noindent As a consequence
of Lemma \ref{independency},
$f^A_{\{\e_n\}}(x;\xi)=f_{\{\e_n\}}(\xi)$ for a.e. $x \in A$ and for
all $\xi \in \Rb^{d \times N}$. It turns out that
$$\mathcal F_{\{\e_{n}\}} (u;A)=\int_{A} f_{\{\e_n\}}(\nabla u(x))\, dx$$
\noindent holds for all $u \in W^{1,p}(A;\Rb^d)$.
$\hfill\blacksquare$


\subsubsection{Characterization of the $\G$-limit}\label{c-gl}

\noindent  Our next objective is to show that $\mathcal
F_{\{\e_n\}}(u;A)=\mathcal F_{\rm hom}(u;A)$ for any $A \in
\mathcal A_0$ and all $u \in W^{1,p}(A;\Rb^d)$. In view of Lemma
\ref{intrep}, we only need to prove that ${\overline f}_{\rm
hom}(\xi)=
 f_{\{\e_{n}\}}(\xi)$ for all $\xi\in \Rb^{d\times N}$, and thus it suffices to work with affine functions
instead  of general Sobolev functions. In order to estimate $
f_{\{\e_{n}\}}$ from below in terms of ${\overline f}_{\rm hom}$,
we will need the following result, close in spirit to Proposition
22.4 in Braides and Defranceschi \cite{BD}.

\begin{proposition}\label{prop} Let $f:\Rb^{N}\times
\Rb^{N}\times \Rb^{d\times N}\rightarrow \Rb^{+}$ be a (not
necessarily continuous) function such that $f(x,\,\cdot\,;\,
\cdot)$ is continuous, $f(\,\cdot, y; \xi)$ is measurable, and
$(H_3)$ and $(H_4)$ hold. Let $A$ be an open, bounded, connected
and Lipschitz subset of\, $\Rb^{N}$. Given $M$ and $\eta$ two
positive numbers, and $\varphi: [0,+ \infty) \rightarrow [0,+
\infty]$  a continuous and increasing function satisfying $\varphi
(t)/t \to + \infty$ as $t \to  +\infty$, there exists
$\e_{0}\equiv \e_{0}(\varphi,M,\eta)>0$ such that for every  $0
<\e < \e_{0}$, every $a\in \Rb^{N}$ and every $u\in
W^{1,p}(a+A;\Rb^{d})$ with

\begin{equation}\label{prop-of-u}
\int_{a+A}\varphi (|\nabla u|^{p})\, dx\leq M,
\end{equation}
\noindent there exists $v \in W^{1,p}_{0}(a+A;\Rb^{d})$ with
$\|v\|_{L^{p}(a+A;\Rb^d)}\leq \eta$ satisfying
$$\int_{a+A} f\Big(x,\frac{x}{\e};\nabla u \Big)\,dx \geq
\int_{a+A} f_{\rm hom}(x;\nabla u + \nabla v)\, dx -\eta.$$
\end{proposition}

\noindent\textit{Proof.} The proof is divided into two steps.
First, we prove this proposition under the additional hypothesis
that $a$ belong to  a compact set of $\Rb^N$. Then, we conclude
the result in its full generality replacing $a$ by its decimal
part $a- \lb a \rb $ and using the periodicity of the integrands
$f$ and $f_{\rm hom}$.

{\it Step 1.}  For $a\in \mathcal [-1,1]^{N}$, the claim of
Proposition \ref{prop} holds. Indeed, if not  we may find
$\varphi$, $M$ and $\eta$ as above, and sequences $\{\e_{n}\} \to
0^{+}$, $\{a_{n}\}\subset [-1,1]^N$  and $\{u_{n}\}\subset
W^{1,p}(a_{n}+A;\Rb^{d})$ with
\begin{equation}\label{bound}
\int_{a_{n}+A}\varphi (|\nabla u_{n}|^{p})\, dx\leq M
\end{equation}
\noindent such that, for every $n\in \Nb$
\begin{eqnarray}\label{ifnot}
&&\int_{a_{n}+A} f\left(x,\frac{x}{\e_{n}};\nabla u_{n} \right)\,dx \nonumber\\
&&\hspace{0.5cm}< \inf_{v\in W^{1,p}_{0}(a_{n}+A;\Rb^{d})}
\left\{\int_{a_{n}+A} f_{\rm hom}(x;\nabla u_{n} + \nabla v)\, dx
: \,  \|v\|_{L^{p}(a_{n}+A;\Rb^d)}\leq \eta \right\} -\eta.
\end{eqnarray}
\noindent From (\ref{bound}) and the Poincar\'{e}-Wirtinger
Inequality, up to a translation argument, we can suppose that the
sequence $\{ \|u_n\|_{W^{1,p}(a_n+A;\Rb^d)} \}$ is uniformly
bounded. From this fact and since the set $a_n+A$ is an extension
domain, there is no loss of generality in assuming that $\{u_n\}$ is
bounded in $W^{1,p}(\Rb^N;\Rb^d)$ and that, due to (\ref{bound}),
\begin{equation}\label{bound1}\debaixodosup {}{n \in \Nb} \int_{\Rb^{N}}\varphi(|\nabla
u_{n}|^{p})\,dx\leq M_{1}\end{equation} for some constant
$M_{1}>0$ depending only on $M$ (see the proof of the  Extension
Theorem for Sobolev functions, Theorem 1, Section 4.4 in Evans and
Gariepy \cite{EG}). Passing to a subsequence, we can also assume
that $u_n \rightharpoonup u$ in $W^{1,p}(\Rb^N;\Rb^d)$. Let $B$ be
a ball of sufficiently large radius so that $a_n +A \subset B$ for
all $n \in \Nb$. De La Vall\'ee Poussin criterion (see e.g.
Proposition 1.27 in Ambrosio, Fusco and Pallara \cite{AFP}) and
(\ref{bound1}) guarantee that  the sequence $\{|\nabla u_{n}|^p\}$
is equi-integrable on $B$. This implies that there exists
$\delta=\delta (\eta)$ such that
\begin{equation}\label{*}
\sup_{n \in \Nb}\beta \int_{E} (1+|\nabla u|^{p}+|\nabla
u_{n}|^{p})\, dx \leq \frac{\eta}{2}
\end{equation}
\noindent  whenever $E$ is a measurable subset of $B$ satisfying
$\mathcal L^N(E) \leq \delta$ and where $\beta$ is the constant
given in $(H_4)$. As $\{a_{n}\}\subset [-1,1]^{N}$ we may suppose,
without loss of generality, that  $a_{n} \to  a \in \mathcal
[-1,1]^{N}$, and that for fixed  $0<\rho< 1$, with $\rho^{N} \ll
\delta$, the following hold for $n$ large enough:
\begin{equation}\label{**}
\left\{
\begin{array}{l}
a+(1-\rho)A\subset a_{n}+A \subset a+(1+\rho)A,\\
\\
\mathcal L^N(S_{n})\leq \delta, \,\, \text{where }S_{n}:=
[a_{n}+A]\setminus [a+(1-\rho)A] \subset B,\\\\
{\rm and}\quad \|u_{n}-u\|_{L^{p}(a+(1+\rho)A;\Rb^{d})}\leq \eta.
\end{array}
\right.
\end{equation}
\noindent Take now a sequence of cut-off functions $\varphi_{n} \in
\mathcal C^\infty_c(\Rb^N;[0,1])$ such that
$$\varphi_{n}=\left\{\begin{array}{rl}
1 & \text{ on } a+(1-\rho)A,
\\&\\
0 & \text { outside } a_{n}+A,
\end{array}\right.$$
\noindent and $\|\nabla \varphi_{n}\|_{L^{\infty}(\Rb^{N})}\leq
C/\rho$ for some constant  $C>0$. Let
$w_{n}=\varphi_{n}u+(1-\varphi_{n})u_{n}$. Then $w_{n}-u_{n}\in
W^{1,p}_0(a_{n}+A;\Rb^{d})$ and
$$\int_{a_{n}+A}|w_{n}-u_{n}|^{p}\,dx\leq \int_{a_{n}+A}\varphi_{n}|u-u_{n}|^{p}\,dx\leq
\int_{a+(1+\rho)A}|u-u_{n}|^{p}\,dx\leq \eta^p.$$ \noindent Then,
taking $v:=w_n-u_n$ as test function in (\ref{ifnot}), it follows
from (\ref{fhom4}), (\ref{*}), and (\ref{**}) that
\begin{eqnarray}\label{2307}
\int_{a_{n}+A} f\left(x,\frac{x}{\e_{n}};\nabla u_{n}
\right)\,dx &<& \int_{a_{n}+A} f_{\rm hom}(x;\nabla w_{n})\,dx-\eta\nonumber\\
& \leq & \int_{a+(1-\rho)A} f_{\rm hom}(x;\nabla u)\,dx \nonumber\\
&& \hspace{1.0cm} + \beta \int_{S_{n}}\left(1+|\nabla
u|^{p}+|\nabla
u_{n}|^{p}+\frac{C}{\rho}|u-u_{n}|^{p}\right)\,dx - \eta\nonumber\\
& \leq & \int_{a+(1-\rho)A} f_{\rm hom}(x;\nabla u)\,dx
-\frac{\eta}{2}+ \frac{\beta\,
C}{\rho}\int_{S_{n}}|u-u_{n}|^{p}\,dx.
\end{eqnarray}
\noindent Since $u_n \rightarrow u$ in $L^{p}(\Rb^N;\Rb^d)$, by
(\ref{**}) and (\ref{2307}) we have
\begin{eqnarray}\label{2308}
\debaixodolimsup {}{n \to  +\infty}\int_{a_{n}+A}
f\left(x,\frac{x}{\e_{n}};\nabla u_{n} \right)\,dx  \leq
\int_{a+(1-\rho)A} f_{\rm hom}(x;\nabla u)\,dx
 -\frac{\eta}{2},
\end{eqnarray}
\noindent and as  $u_n \rightharpoonup u$ in
$W^{1,p}(a+(1-\rho)A;\Rb^d)$, by Theorem 1.1 in Ba\'{\i}a and
Fonseca \cite{B&Fo} and (\ref{2308}), we get
\begin{eqnarray*}
\int_{a+(1-\rho)A} f_{\rm hom}(x;\nabla u)\,dx & \leq &
\debaixodoliminf {}{n \to  +\infty}\int_{a+(1-\rho)A}
f\left(x,\frac{x}{\e_{n}};\nabla u_{n}
\right)\,dx\\
&\leq & \debaixodoliminf {}{n \to  +\infty}\int_{a_{n}+A}
f\left(x,\frac{x}{\e_{n}};\nabla u_{n}
\right)\,dx\nonumber\\
& \leq & \int_{a+(1-\rho)A} f_{\rm hom}(x;\nabla u)\,dx
-\frac{\eta}{2}
\end{eqnarray*}
\noindent which is a contradiction.

{\it Step 2. (General case)} Let $a \in \Rb^N$. Then $a-\llbracket a
\rrbracket \in [-1,1]^N$. Given $u \in W^{1,p}(a+A;\Rb^d)$, set
$\tilde u(x):=u(x+\llbracket a \rrbracket)$ and thus $\tilde u \in
W^{1,p}(a- \llbracket a \rrbracket +A;\Rb^d)$. Applying Step 1 with
$\eta/3$, we get the existence of $0<\e_0'\equiv
\e_0'(M,\varphi,\eta)$ such that, for all $0<\e<\e'_0$, there exist
$\tilde v \in W^{1,p}_0(a- \llbracket a \rrbracket +A;\Rb^d)$
satisfying $\|\tilde v\|_{L^p(a-\llbracket a \rrbracket+A;\Rb^d)}
\leq \eta/3$ and
$$\int_{a- \llbracket a \rrbracket+A}f\left(x,\frac{x}{\e};\nabla \tilde u(x)\right) \,
dx \geq \int_{a-\llbracket a \rrbracket +A}f_{\rm hom}(x;\nabla
\tilde u(x)+\nabla \tilde v(x))
 \, dx -\frac{\eta}{3}.$$
\noindent Setting $v(x):=\tilde v(x-\llbracket a \rrbracket)$, then
$v \in W^{1,p}_0(a+A;\Rb^d)$ and  $\| v\|_{L^p(a+A;\Rb^d)}
\leq\frac{\eta}{3}\leq \eta$. Therefore, by a change of variables
\begin{equation}\label{2315}\int_{a+A}f\left(x,\frac{x-\llbracket a
\rrbracket}{\e};\nabla u(x)\right) \, dx \geq \int_{a+A}f_{\rm
hom}(x;\nabla u(x)+\nabla v(x))
 \, dx -\frac{\eta}{3},\end{equation}
\noindent where we have used condition $(H_3)$ and (\ref{fhom3}).
Writing
$$\frac{\llbracket a \rrbracket}{\e}=: m_\e+r_\e, \quad \text{ with }m_\e \in \Zb^N\,\,
\text{ and }\,\,|r_\e| <\sqrt N \e$$ \noindent by  $(H_3)$ the
inequality (\ref{2315})  reduces to
\begin{equation}\label{d/3_1}
\int_{a+A}f\left(x,\frac{x}{\e}-r_\e;\nabla u(x)\right) \, dx \geq
\int_{a+A}f_{\rm hom}(x;\nabla u(x)+\nabla  v(x)) \, dx
-\frac{\eta}{3}.
\end{equation}
Choose $\lambda>0$ large enough (depending on $\eta$) so that
\begin{equation}\label{2322}\beta \int_{\{|\nabla u|>\lambda\} \cap [a+A]}(1+|\nabla
u|^p)\, dx \leq \frac{\eta}{6}.\end{equation} \noindent Fixed $\rho
>0$, by  Scorza-Dragoni's  Theorem there exists a compact set
$K_{\rho}\subset a+A$ with $\mathcal L^{N} ([a+A]\setminus
K_{\rho})\leq \rho$ such that $f: K_{\rho}\times \Rb^{N}\times
\Rb^{d\times N} \to \Rb$ is continuous. Take $\rho\equiv \rho(\eta)$
small enough as for
\begin{equation}\label{1136}
\beta \int_{[a+A]\setminus K_{\rho}}(1+|\nabla
u|^p)\, dx \leq \frac{\eta}{6}.
\end{equation}
\noindent Then, from (\ref{2322}), (\ref{1136}) and the $p$-growth
condition $(H_4)$
\begin{equation}\label{d/3_3}
\int_{a+A}f\left(x,\frac{x}{\e}-r_\e;\nabla u(x)\right) \, dx \leq
\int_{\{|\nabla u|\leq \lambda\}\cap
K_{\rho}}f\left(x,\frac{x}{\e}-r_\e;\nabla u(x)\right) \,
dx+\frac{\eta}{3}.
\end{equation}
\noindent But since $f$ is $Q$-periodic in its second variable, then
$f$ is uniformly continuous on $K_{\rho} \times \Rb^N \times
\overline B(0,\lambda)$.  Thus, as $r_{\e} \to  0$, for any $\eta>0$
there exists $\e''_0\equiv \e''_0(\eta)>0$ such that for all
$\e<\e''_0$ and all $x \in \{|\nabla u|\leq \lambda\}\cap K_{\rho}$,
$$\left|f\left(x,\frac{x}{\e}-r_\e;\nabla
u(x)\right)-f\left(x,\frac{x}{\e};\nabla u(x)\right)\right| <
\frac{\eta}{3\mathcal L^{N}(A)}.$$ \noindent Hence
\begin{equation}\label{d/3_4}
\int_{\{|\nabla u|\leq \lambda\}\cap
K_{\rho}}f\left(x,\frac{x}{\e}-r_\e;\nabla u(x)\right) \, dx\leq
\int_{\{|\nabla u|\leq \lambda\}\cap K_{\rho}}
f(x,\frac{x}{\e},\nabla u)\,dx + \frac{\eta}{3},
\end{equation}
\noindent  and consequently, by (\ref{d/3_1}), (\ref{d/3_3}) and
(\ref{d/3_4}), for all $\e<\e_0:=\min\{\e'_0,\e''_0\}$ we have
$$\int_{a+A}f\left(x,\frac{x}{\e};\nabla u(x)\right) \, dx \geq \int_{a+A}f_{\rm hom}(x;\nabla u(x)+\nabla
 v(x)) \, dx -\eta.$$
$\hfill\blacksquare$\\

We are now in position to prove that $f_{\{\e_n\}}=\overline f_{\rm hom}$.

\begin{lemma}\label{lem<}
For all $\xi\in \Rb^{d\times N}$, ${\overline f}_{\rm hom}(\xi)
\leq  f_{\{\e_{n}\}}(\xi)$.
\end{lemma}

\noindent\textit{Proof. } By Lemma \ref{intrep}, given $\xi \in
\Rb^{d \times N}$  there exists a sequence $\{w_n\} \subset
W^{1,p}_0(Q;\Rb^d)$ such that $w_n \to 0$ in $L^p(Q;\Rb^d)$ and
\begin{equation}\label{f-e}
f_{\{\e_{n}\}}(\xi)=\lim_{n\to +\infty}\int_{Q}
f\left(\frac{x}{\e_{n}}, \frac{x}{{\e_{n}}^{2}};  \xi+ \nabla
w_{n}(x)\right)\, dx
\end{equation}
\noindent (see e.g. Proposition 11.7 in Braides and Defranceschi
\cite{BD}). Following the same argument as in Lemma
\ref{independency}, by the Decomposition Lemma, there is no loss of
generality in assuming that $\{ |\nabla w_n|^p\}$ is
equi-integrable. Thus, from De La Vall\'ee Poussin criterion (see
e.g. Proposition 1.27 in Ambrosio, Fusco and Pallara \cite{AFP})
there exists an increasing continuous function $\varphi:
[0,+\infty)\rightarrow [0,+\infty]$ satisfying $\varphi(t) / t \to
+\infty$ as $t \to  +\infty$ and such that
$$\debaixodosup {}{n \in \Nb}\int_{Q} \varphi(|\nabla
w_{n}|^{p})\,dx\leq 1.$$ \noindent Changing variables (\ref{f-e})
yields
$$f_{\{\e_n\}}(\xi)= \lim_{n \to +\infty}
\frac{1}{T_n^N}\int_{(0,T_n)^N} f\left( x,\frac{x}{\e_n};\xi+ \nabla
z_n(x) \right)\, dx$$
and
\begin{equation}\label{val-pou-bis}\sup_{n \in \Nb} \frac{1}{T_{n}^{N}}\int_{(0,T_{n})^{N}} \varphi
(|\nabla z_{n}|^{p})\,dx \leq 1,\end{equation} \noindent where we
set $T_n:=1/\e_n$ and $z_n(x):=T_n w_n (x/T_n)$ with $z_n \in
W^{1,p}_0((0,T_n)^N;\Rb^d)$.  For any $n \in \Nb$ define $I_n:=
\left\{1,..., \lb T_{n} \rb ^{N}\right\}$, and for $i\in I_{n}$
take  $a_{i}^{n}\in \Zb^{N}$ such that
\begin{equation}\label{lim0} \bigcup_{i\in I_{n}} (a_{i}^n+Q)\subseteq
(0,T_n)^N.\end{equation} \noindent Thus
\begin{equation}\label{equi8}{f}_{\{\e_n\}}(\xi) \geq \ds \debaixodolimsup {}{n \to
+\infty}\frac{1}{T_n^{N}}\sum_{i\in I_{n}}\int_{a_{i}^n+Q}f\left( x,
\frac{x}{\e_n};\xi+ \nabla {z}_n(x) \right)\,dx.\end{equation}
\noindent Let $M> 2$ and $\eta>0$.  For  $n\in \Nb$ define
$$I_{n}^{M}:=\left\{i\in I_{n}: \int_{a_{i}^n+Q}\varphi(|\nabla z_{n}|^{p})\,dx\leq M \right\}.$$
\noindent We note that  for any $M > 2$, there exists $n(M) \in \Nb$
such that for all $n \geq n(M)$ sufficiently large so that
$T_{n}>M$, $I_n^M \ne \emptyset$. In fact, if not we may find $M >
2$ and a subsequence $n_k \in \Nb$ satisfying

$$\int_{a_i^{n_k}+Q} \varphi(|\nabla z_{n_k}|^p)dx >M,$$
for all $i \in I_{n_k}$. Summation in $i$ and (\ref{lim0}) would
yield to
$$\int_{(0,T_{n_k})^N}\varphi(|\nabla z_{n_k}|^p)dx>M \lb T_{n_k} \rb^N$$
 \noindent which is in contradiction with (\ref{val-pou-bis}).
We also note that in view of (\ref{val-pou-bis})

$${\rm Card}(I_n \setminus I_n^M) M\leq \sum_{i\in I_n \setminus I_n^M}\int_{
 a_i^n +Q } \varphi(|\nabla
z_{n}|^{p})\,dx\leq  \int_{(0,T_{n})^{N}}\varphi(|\nabla
z_{n}|^{p})\,dx \leq  T_{n}^{N},$$

\noindent and so
\begin{equation}\label{card}
{\rm Card}(I_n \setminus I_n^M) \leq \frac{ T_n^N}{M}.
\end{equation}
\noindent By Lemma \ref{prop} there exists $\e_0 \equiv
\e_0(M,\eta)$ such that, for any $n$ large enough satisfying $0 \leq
\e_n< \e_{0}$ and for any $i \in I_{n}^{M}$, we can find
$v_i^{n,M,\eta} \in W_{0}^{1,p}(a_{i}^n+Q;\Rb^{d})$ with
$\|{v}_i^{n,M,\eta}\|_{L^{p}(a_{i}^n+Q;\Rb^{d})}\leq \eta$ and
$$\int_{a_{i}^n+Q}f\left(
x, \frac{x}{\e_n};\xi+ \nabla {z}_n(x) \right)\, dx \geq
\int_{a_{i}^n+Q}f_{\rm hom}\left( x; \xi +\nabla {z}_n+ \nabla
{v}_i^{n,M,\eta} \right)\, dx-\eta.$$ \noindent Consequently, for
$n$ large enough
$$\sum_{i\in I_{n}}\int_{a_{i}^n+Q}f\left( x, \frac{x}{\e_n};\xi+
\nabla {z}_n(x) \right)\, dx \geq  \sum_{i\in I_{n}^{M}}
\int_{a_{i}^n +Q} f_{\rm hom}( x;\xi+\nabla {z}_n+ \nabla
{v}_i^{n,M,\eta})\, dx - \eta \; {\rm card}(I_{n}^{M}).$$
\noindent As ${\rm Card}(I_{n}^{M})\leq \lb T_{n} \rb ^{N}$,
dividing by $T_n^N$ and passing to the limit when $n \to +\infty$
we obtain  from (\ref{equi8})

\begin{equation}\label{equi5}
f_{\{\e_n\}}(\xi)  \geq  \debaixodolimsup {}{M, \eta, n}
\frac{1}{T_n^N} \sum_{i \in I_n^M}\int_{a_{i}^n+Q}
 f_{\rm hom}( x;\xi+\nabla  {\phi}^{n,M,\eta})\, dx
\end{equation}
\noindent where $\phi^{n,M,\eta} \in W^{1,p}_0((0;T_n)^N;\Rb^d)$ is
defined by
$$\phi^{n,M,\eta}(x):=\left\{\begin{array}{ll}
z_{n}(x)+v_i^{n,M,\eta}(x) & \text{if } x \in a_{i}^n+Q \text{ and } i \in I_{n}^{M},\\
&\\
z_n(x) & \text{otherwise}.
\end{array}\right.$$
\noindent Now, in view of the definition of $\phi^{n,M,\eta}$, the
$p$-growth condition (\ref{fhom4}) and (\ref{card}),
\begin{eqnarray}\label{equi10}
\frac{1}{T_n^N}\sum_{i\in  I_n \setminus I_n^M} \int_{a_i^n +Q}
f_{\rm hom}( x; \xi+\nabla \phi^{n,M,\eta})\, dx & =
&\frac{1}{T_n^N} \sum_{i \in I_n \setminus I_n^M}
\int_{a_i^n +Q}  f_{\rm hom}( x; \xi+\nabla z_n)\, dx\nonumber\\
&\leq & \frac{\beta }{T_n^N}\sum_{i \in I_n \setminus I_n^M}
\int_{a_i^n +Q}( 1+|\nabla  z_n|^p)\, dx\nonumber\\
&\leq & \frac{\beta}{M} + \frac{\beta}{T_n^N} \sum_{i \in I_n
\setminus I_n^M}
\int_{a_i^n +Q}|\nabla  z_n|^p\, dx\nonumber\\
& \leq & \frac{\beta}{M}+ \beta \int_{\bigcup\limits_{i \in I_n
\setminus I_n^M} \frac{1}{T_n}(a_i^n +Q)} |\nabla  w_n|^p\, dx.
\end{eqnarray}
\noindent By (\ref{card}) we have that
$${\cal L}^{N} \left( \debaixodauniao {}{i \in I_n
\setminus I_{n}^{M}} \frac{1}{T_{n}}(a_{i}^{n} +Q) \right) \leq
\frac{1}{M}.$$ \noindent Consequently, in view of the
equi-integrability of $\{|\nabla w_n|^p\}$ and  (\ref{equi10}), we
get
\begin{equation}\label{equi11}
\debaixodolimsup {}{M, \eta, n} \frac{1}{T_{n}^{N}} \sum_{i\in I_n
\setminus I_{n}^{M}} \int_{a_{i}^{n} +Q} f_{\rm hom}( x; \xi+\nabla
{\phi}^{n,M,\eta})\, dx=0.
\end{equation}
\noindent Therefore, (\ref{equi5}) and (\ref{equi11}) imply
\begin{equation}\label{finaleq}
f_{\{\e_n\}}(\xi) \geq \debaixodolimsup {}{M, \eta, n}
\frac{1}{T_n^N} \sum_{i \in I_n}\int_{a_{i}^n+Q} f_{\rm hom}(
x;\xi+\nabla  {\phi}^{n,M,\eta})\, dx.
\end{equation}

\noindent Similarly, since
$${\cal L}^{N}\left(Q\setminus \left[\debaixodauniao{}{i\in I_{n}}\frac{1}{T_{n}}(a_{i}^{n}+Q)\right]\right) \to  0$$

\noindent  as $n \to +\infty$, we get that
\begin{equation}\label{finaleqbis}\debaixodolimsup
{}{M, \eta, n}
\frac{1}{T_{n}^{N}}\int_{\big(0,T_{n}^{N}\big)\setminus
\left[\debaixodauniao{}{i\in I_{n}} (a_{i}^{n}+Q) \right]} f_{\rm
hom}( x;\xi+\nabla  {\phi}^{n,M,\eta})\, dx=0.
\end{equation}

\noindent Hence by (\ref{finaleq}), (\ref{finaleqbis}) and
(\ref{fhombarra2}) we get that
$$f_{\{\e_{n}\}}(\xi)\geq \debaixodolimsup {}{M, \eta, n} \frac{1}{T_n^N}
\int_{(0,T_{n})^{N} } f_{\rm hom}( x;\xi+\nabla
{\phi}^{n,M,\eta})\, dx \geq {\overline f}_{\rm hom}(\xi).$$
$\hfill\blacksquare$\\

Let us now prove the converse inequality.

\begin{lemma}\label{lem>}
For all $\xi\in \Rb^{d\times N}$, ${\overline f}_{\rm hom}(\xi)
\geq  f_{\{\e_{n}\}}(\xi)$.
\end{lemma}
\noindent\textit{Proof. } In view of (\ref{fhombarra2}), for
$\delta>0$ fixed take  $T\equiv T_{\delta} \in \Nb$, with
$T_{\delta} \to  +\infty$ as $\delta  \to  0$, and $\phi\equiv
{\phi_{\delta}}\ \in W_{0}^{1,p}((0,T)^{N};\Rb^{d})$ be such that
\begin{equation}\label{delta}
{\overline f}_{\rm hom}(\xi)+\delta \geq \frac{1}{{T}^{N}}
\int_{(0,T)^{N}}f_{\rm hom}(x; \xi+\nabla {\phi}(x))\, dx.
\end{equation}
\noindent By Theorem 1.1 in Ba\'{\i}a and Fonseca \cite{B&Fo} and,
for instance, Proposition 11.7 in Braides and Defranceschi
\cite{BD}, there exists a sequence $\{\phi_{n}\}\subset
W_{0}^{1,p}((0,T)^{N};\Rb^{d})$ with $\phi_{n}\rightarrow \phi$ in
$L^{p}((0,T)^{N};\Rb^{d})$ such that
\begin{equation}\label{b-n}
\int_{(0,T)^{N}}f_{\rm hom}(x; \xi+\nabla {\phi}(x))\, dx=\lim_{n
\to +\infty}\int_{(0,T)^{N}}f\left(x,\frac{x}{\e_{n}};\xi+\nabla
\phi_{n}(x)\right) dx.
\end{equation}
\noindent Further, in view of the Decomposition Lemma (see
Fonseca, M\"uller and Pedregal \cite{FMP} or Fonseca and Leoni
\cite{FL}), we can assume -- upon extracting a subsequence still
denoted by $\{\phi_n\}$ -- $\{|\nabla \phi_n|^p\}$ to be
equi-integrable. Fix $n\in \Nb$ such that  $\e_{n} \ll 1$. For all
$i \in \Zb^N$ let $a_{i}^n \in \e_n \Zb^N \cap  (i(T+1)+
[0,\e_n)^N)$ (uniquely defined) (see Example in Figure 4).

\begin{figure}[h]
\begin{center}
\scalebox{.80}{\includegraphics{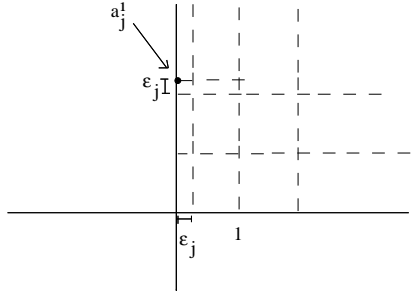}}
\end{center}
\caption{Example for $T=1$, $i=(0,1)$ and $N=2$}
\end{figure}

\noindent In particular, the cubes $a_i^n + (0,T)^N$ are not
overlapping because if $i,j \in \Zb^N$ with $i \neq j$, then $|i-j|
\geq 1$ and thus $|a_i^n - a_j^n|> T$. Set
$$\tilde \phi_n(x):=\left\{
\begin{array}{cl}
\phi_n(x-a_{i}^n) & \text{if }  x \in a_i^n + (0,T)^N \text{ and } i\in \Zb^{N},\\
&\\
0 & \text{otherwise},
\end{array}
\right.$$ then $\tilde \phi_n \in W^{1,p}(\Rb^N;\Rb^d)$. Let
$I_n:=\{ i \in \Zb^N : (0,T/\e_n)^N \cap (a_i^n + (0,T)^N) \ne
\emptyset \}$. Note that
\begin{equation}\label{card-a}
{\rm Card}(I_n) \leq \left(\left \llbracket \frac{1}{\e_n}
\right \rrbracket +1 \right)^N.
\end{equation}
\noindent If $\psi_n(x):=\e_n \tilde \phi_n(x/\e_n)$ then $\psi_n
\to 0$ in $L^p((0,T)^N;\Rb^d)$,  as $n  \to  +\infty$, because

\begin{eqnarray*}
\int_{(0,T)^N} |\psi_n(x)|^p dx & = & \e_n^p \int_{(0,T)^N} \left|\tilde \phi_n\left(\frac{x}{\e_n}\right)\right|^p dx\\
& = & \e_n^{p+N}\int_{(0,T/\e_n)^N} |\tilde \phi_n(x)|^p dx\\
& \leq & \e_n^{p+N} \sum_{i \in I_n} \int_{a_i^n + (0,T)^N}|\phi_n(x-a_i^n)|^p dx\\
& = & \e_n^{p+N} {\rm Card}(I_n) \int_{(0,T)^N} |\phi_n(x)|^p dx\\
& \leq & \e_n^{p+N} \left(\left \llbracket \frac{1}{\e_n}  \right
\rrbracket +1 \right)^N\int_{(0,T)^N} |\phi_n(x)|^p dx \to 0,
\end{eqnarray*}

\noindent where we have used the fact that $\tilde \phi_n \equiv 0$
on $(0,T/\e_n)^N \setminus \bigcup_{i \in I_n} (a_i^n + (0,T)^N)$
and that
$$\debaixodosup {} {n \in \Nb} \|\phi_{n}\|_{L^{p}((0,T)^N;\Rb^d)}<
+\infty$$ by the Poincar\'{e} Inequality and (\ref{b-n}).
Consequently,
\begin{eqnarray*}\mathcal F_{\{\e_n\}}(\xi \cdot\,;(0,T)^N) & \leq & \liminf_{n
\to
+\infty} \int_{(0,T)^N} f \left( \frac{x}{\e_n},\frac{x}{\e_n^2};\xi +\nabla \psi_n(x)\right)\, dx\\
& = & \liminf_{n \to +\infty} \int_{(0,T)^N} f \left(
\frac{x}{\e_n},\frac{x} {\e_n^2};\xi +\nabla \tilde \phi_n\left(\frac{x}{\e_n}\right)\right)\, dx\\
& = & \liminf_{n \to +\infty} \e_n^N \int_{(0,T/\e_n)^N} f \left(
x,\frac{x}{\e_n}; \xi +\nabla \tilde \phi_n(x)\right)\, dx.
\end{eqnarray*}
\noindent Note that since
\begin{eqnarray*}
\e_n^N \mathcal L^N \left( \left( 0 , \frac{T}{\e_n} \right)^N
\setminus \bigcup_{i \in I_n} (a_i^n + (0,T)^N) \right) & \leq &
T^N \left( 1 - \e_n^N \left\llbracket \frac{T}{\e_n(T+1)}
\right\rrbracket^N \right),
\end{eqnarray*}
from the $p$-growth condition $(H_4)$ it follows that
\begin{eqnarray*}
\mathcal F_{\{\e_n\}}(\xi \cdot\,;(0,T)^N) & \leq & 
\liminf_{n \to +\infty} \e_n^N \sum_{i \in I_n}\int_{a_i^n +
(0,T)^N} f \left( x,\frac{x}{\e_n};\xi +\nabla \tilde
\phi_n(x)\right)\, dx \\
&& \hspace{2.0cm} + \beta (1+|\xi|^p) T^{N} \left(1-
\left(\frac{T}{T+1} \right)^{N} \right).
\end{eqnarray*}
\noindent By a change of variables, for all $i\in I_{n}$
\begin{eqnarray*}
&&\int_{a_i^n + (0,T)^N} f \left( x,\frac{x}{\e_n};\xi +\nabla
\tilde
\phi_n(x)\right)\, dx\\
&&\hspace{3.0cm} = \int_{(0,T)^N} f \left(x+a_i^n ,\frac{x+a_i^n}{\e_n};\xi +\nabla \phi_n(x)\right)\, dx\\
&&\hspace{3.0cm} = \int_{(0,T)^N} f \left(x+a_i^n-i(T+1) ,\frac{x}{\e_n};\xi
+\nabla \phi_n(x)\right)\, dx,
\end{eqnarray*}
\noindent where we have used $(H_3)$, the fact that $T \in \Nb$ and
$a_i^n/\e_n \in \Zb^N$. By a similar uniform continuity argument
than the one used in Lemma \ref{independency} and since
(\ref{card-a}) holds, $|a_i^n-i(T+1)| \leq \e_n$, and $\{|\nabla
\phi_n|^p\}$ is equi-integrable,  it follows that
\begin{eqnarray}\label{glim2}
\mathcal F_{\{\e_n\}}(\xi \cdot\,;(0,T)^N) \leq \liminf_{n \to
+\infty} \int_{(0,T)^N} f \left( x,\frac{x}{\e_n};\xi +\nabla \tilde
\phi_n(x)\right)\, dx + c T^{N} \left(1- \left(\frac{T}{T+1}
\right)^{N} \right).
\end{eqnarray}

\noindent Consequently by (\ref{delta}), (\ref{b-n}, (\ref{glim2})
and Lemma \ref{intrep}
$$ f_{\{\e_n\}}(\xi) \leq \overline f_{\rm hom} (\xi) +\delta +  c \left(1-
\left(\frac{T}{T+1}\right)^{N}\right).$$ The result follows by
letting $\delta$ tend to zero.
$\hfill\blacksquare$\\

\noindent {\it Proof of Theorem \ref{jf-m2}.} From Lemma \ref{lem<}
and Lemma \ref{lem>}, we conclude that ${\overline f}_{\rm hom}(\xi)
= f_{\{\e_{n}\}}(\xi)$ for all $\xi\in \Rb^{d\times N}$. As a
consequence, $\mathcal F_{\{\e_n\}}(u;A)= \mathcal F_{\rm hom}(u;A)$
for all $A \in \mathcal A_0$ and all $u \in W^{1,p}(A;\Rb^d)$. Since
the $\G$-limit does not depend upon the extracted subsequence,
Proposition 8.3 in Dal Maso \cite{DM} implies that the whole
sequence $\mathcal F_\e(\cdot\,;A)$ $\G(L^p(A))$-converges to
$\mathcal F_{\rm hom}(\cdot\,;A)$. 
$\hfill\blacksquare$


\subsection{The general case}\label{2hetero}

\noindent Our aim here is to prove Theorem \ref{jf-m1}.

\subsubsection{Existence and integral representation of the $\G$-limit}

\noindent   The  idea in this case is to freeze the macroscopic
variable and to use Theorem \ref{jf-m2} through a blow up argument.
This leads us  to work on small cubes centered at convenient
Lebesgue points of  $\O$ which, contrary to Section \ref{part2},
allow us  to localize our functionals on $\mathcal A(\O)$, the
family of open subsets of $\O$. We define $\mathcal F_\varepsilon:
L^{p}(\O;\Rb^{d})\times {\mathcal A}(\Omega)\rightarrow [0,+\infty]$
by
$${\cal F}_{\e}(u;A):=\left\{\begin{array}{ll} \ds \int_{A}
f\left(x,\frac{x}{\e},\frac{x}{\e^{2}};\nabla
 u(x)\right)\, dx & \text{if } u\in W^{1,p}(A;\Rb^{d}),\\
  +\infty & \text{otherwise},
\end{array}\right.$$
\noindent   and we introduce the functional $\mathcal F_{\rm hom}:
L^p(\O;\mathbb R^d) \times \mathcal A(\O) \rightarrow [0,+\infty]$
$$\mathcal F_{\rm hom}(u;A) :=
\left\{\begin{array}{ll} \ds \int_{A} \overline f_{\rm hom}(x;\nabla
u(x))\,dx & \text{if } u  \in W^{1,p}(A;\Rb^{d}),\\ + \infty
&\text{otherwise}.
\end{array}\right.$$
\noindent Given $\{\e_j\} \searrow 0^+$  and $A \in \mathcal A
(\O)$, consider  the $\Gamma$-lower limit  of $\{  \mathcal F_{\e_j}
(\cdot\,;A) \}_{j \in \Nb}$ for the $L^{p}(A;\Rb^{d})$-topology
defined, for $u \in W^{1,p}(\O;\Rb^d)$, by
$$\mathcal F_{\{\e_j\}}(u;A) := \inf_{\{u_j\}}\left\{\liminf_{j \to
+\infty} \mathcal F_{\e_j}(u_j;A) : \quad u_j \to u \text{ in
}L^p(A;\Rb^d)\right\}.$$ Due to the $p$-coercivity condition in
$(H_{4})$, to prove Theorem \ref{jf-m1} it suffices to show that for
all $u\in W^{1,p}(\Omega;\Rb^{d})$
$$\Gamma ( L^{p}(\O))\text{-}\lim_{\e \rightarrow 0}  {\cal F}_{\e}
(u) = \int_{\Omega} {\overline f}_{\rm hom}(x,\nabla u(x))\, dx.$$
\noindent As a consequence  of Theorem 8.5 in Dal Maso \cite{DM},
there exists a subsequence $\{\e_n\} \equiv \{\e_{j_n}\}$ such that
for any $A \in \mathcal A(\O)$, ${\cal F}_{\{\e_n\}}(\cdot\,;A)$ is
the $\Gamma(L^p(A))$-limit of $\mathcal F_{\e_n}(\cdot\,;A)$ and,
for all $u \in W^{1,p}(\O;\Rb^d)$, the set function ${\cal
F}_{\{\e_n\}}(u;\cdot)$ is the restriction of a Radon measure to
$\mathcal A(\O)$. Furthermore, from  Buttazzo-Dal Maso Integral
Representation Theorem (see Theorem 1.1 in \cite{DMB}) it follows
that
\begin{lemma}\label{intrep2}
There exists a Carath\'{e}odory function  $ f_{\{\e_{n}\}}: \O
\times \Rb^{d\times N}\rightarrow \Rb$, quasiconvex in its second
variable, satisfying the same coercivity and growth conditions than
$f$, such that
$$\mathcal F_{\{\e_{n}\}}(u;A)=\int_{A}
f_{\{\e_{n}\}}(x;\nabla u(x))\,dx$$ \noindent for every
$A\in\mathcal A(\O)$ and $u\in W^{1,p}(\O;\Rb^{d})$. Moreover, for
all $\xi \in \Rb^{d \times N}$ and a.e. $x \in \O$,
$$ f_{\{\e_n\}}(x;\xi)=\lim_{\delta  \to  0}\frac{\mathcal F_{\{\e_{n}\}}(\xi\cdot\,;Q(x,\delta))}{\delta^N}.$$
\end{lemma}


\subsubsection{Characterization of the $\G$-limit}

\noindent Like in Section 4, we only need to prove that $
f_{\{\e_n\}}(x;\xi)=\overline f_{\rm hom}(x;\xi)$ for a.e.  $x$ and
all $\xi$. For this purpose let $L$ be the set of  Lebesgue points
$x_{0}$ for all functions $f_{\{\e_n\}}(\cdot\,;\xi)$ and $\overline
f_{\rm hom}(\cdot\,;\xi)$, for all $\xi \in \Qb^{d \times N}$. We
have  $\mathcal L^N(\O \setminus L)=0$ and  we will first show in
Lemma \ref{ineq<} and \ref{ineq>} below that the equality
$f_{\{\e_n\}}(x;\xi) = \overline f_{\rm hom}(x;\xi)$ holds for all
$x \in L$ and all $\xi \in \Qb^{d \times N}$. By definition of the
set $L$ it is enough to show that
\begin{equation}\label{009}
\int_{Q(x_0,\delta)} f_{\{\e_n\}}(x;\xi)\, dx=
\int_{Q(x_0,\delta)}\overline f_{\rm hom}(x;\xi)\, dx,
\end{equation}
for every $x_0\in L$ and each $\delta>0$ small enough so that
$Q(x_0,\delta) \in \mathcal A(\Omega)$.
\begin{lemma}\label{ineq<}
For all $\xi \in \Qb^{d \times N}$ and all $x_{0} \in L$,
$$\int_{Q(x_0,\delta)} f_{\{\e_n\}}(x;\xi)\,dx \geq \int_{Q(x_0,\delta)}\overline f_{\rm hom}(x;\xi)\, dx.$$
\end{lemma}
\noindent\textit{Proof. } Let $\xi\in \Qb^{d\times N}$ and $x_0
\in L$.
\noindent Let $\{u_n\} \subset
W^{1,p}(Q(x_0,\delta);\Rb^d)$ be a recovering sequence for
$\mathcal F_{\{\e_n\}}(\xi \cdot\,; Q(x_0;\delta))$, that is, a
sequence $\{u_{n}\}$ such that $u_n \to 0$ in
$L^p(Q(x_0,\delta);\Rb^d)$ and
$$\int_{Q(x_0,\delta)} f_{\{\e_n\}}(x;\xi)\, dx = \mathcal F_{\{\e_n\}}(\xi \cdot\,; Q(x_0;\delta)) = \lim_{n \to +\infty}
\int_{Q(x_0,\delta)} f\left(x,\frac{x}{\e_n},\frac{x}{\e_n^2};\xi
+\nabla u_n(x) \right) \, dx.$$ \noindent As before, the
Decomposition Lemma (see Fonseca, M\"uller and Pedregal \cite{FMP}
or Fonseca and Leoni \cite{FL}) let us to assume that $\{ |\nabla
u_n|^p\}$ is equi-integrable. We split $Q(x_0,\delta)$ into $h^N$
small disjoint cubes $Q_{i,h}$ such that
\begin{equation}\label{div}Q(x_0,\delta)=\bigcup_{i=1}^{h^N}Q_{i,h}\quad  \text{ and }
\quad \mathcal L^N(Q_{i,h}) = (\delta/h)^N.\end{equation} \noindent
Then
$$\int_{Q(x_0,\delta)} f_{\{\e_n\}}(x;\xi)\, dx =\lim_{h, n} \sum_{i=1}^{h^N} \int_{Q_{i,h}}
f\left(x,\frac{x}{\e_n},\frac{x}{\e_n^2};\xi +\nabla u_n(x) \right)
\, dx.$$ \noindent Let $\eta>0$.  By  Scorza-Dragoni's Theorem (see
Ekeland and Temam \cite{Ek&Te}), there exists a compact set $K_\eta
\subset \Omega$ such that
\begin{equation}\label{eta}
\mathcal L^N(\Omega \setminus K_\eta) < \eta,
\end{equation}
and the restriction of $f$ to $K_\eta \times \Rb^N \times \Rb^N
\times \Rb^{d \times N}$ is a continuous function. Given
$\lambda>0$, we introduce
$$R^\lambda_n:=\{x \in \Omega : |\xi +\nabla u_n(x)| \leq \lambda\},$$
\noindent for all $n\in \Nb$ and we note that due to Chebyshev's
inequality, we have
\begin{equation}\label{lambda}
\mathcal L^N(\Omega \setminus R^\lambda_n) \leq
\frac{C}{\lambda^p},
\end{equation}
for some constant $C>0$ independent of $n$ and $\lambda$.  Then
$$\int_{Q(x_0,\delta)} f_{\{\e_n\}}(x;\xi)\, dx  \geq  \limsup_{\lambda,\eta,h,n}\sum_{i=1}^{h^N} \int_{Q_{i,h} \cap K_\eta \cap R^\lambda_n}
f\left(x,\frac{x}{\e_n},\frac{x}{\e_n^2};\xi +\nabla u_n(x) \right)
\, dx.$$ In view of condition $(H_3)$, $f$ is uniformly continuous
on $K_\eta \times \Rb^N \times \Rb^N \times \overline B(0,\lambda)$.
Denoting by $\omega_{\eta,\lambda} :\Rb^+ \to \Rb^+$ the modulus of
continuity of $f$ on $K_\eta \times \Rb^N \times \Rb^N \times
\overline B(0,\lambda)$, for every $(x,x') \in [Q_{i,h} \cap K_\eta
\cap R^\lambda_n] \times [Q_{i,h} \cap K_\eta]$,
\begin{eqnarray}\label{023}
\left|f\left(x,\frac{x}{\e_n},\frac{x}{\e_n^2};\xi +\nabla u_n(x) \right)-f\left(x',\frac{x}
{\e_n},\frac{x}{\e_n^2};\xi +\nabla u_n(x) \right)\right|& \leq & \omega_{\eta,\lambda}(|x-x'|)\nonumber\\
& \leq & \omega_{\eta,\lambda}\left(\frac{\sqrt{N}\delta}{h} \right).
\end{eqnarray}
\noindent From (\ref{div}) and (\ref{023}), after integrating in
$(x,x')$ over $[Q_{i,h} \cap K_\eta \cap R^\lambda_n] \times
[Q_{i,h} \cap K_\eta]$, we get since $\omega_{\eta,\lambda}$ is
continuous and satisfies $\omega_{\eta,\lambda}(0)=0$
\begin{eqnarray*}
&& \sum_{i=1}^{h^N}\frac{h^N}{\delta^N}\int_{Q_{i,h}\cap K_\eta} \left\{ \int_{Q_{i,h} \cap K_\eta \cap R^\lambda_n} \left|f\left(x,\frac{x}{\e_n},\frac{x}{\e_n^2};\xi +\nabla u_n(x) \right)\right. \right.\\
&& \hspace{8.0cm} \left. \left. -f\left(x',\frac{x}{\e_n},\frac{x}{\e_n^2};\xi +\nabla u_n(x) \right)\right|\, dx\right\}dx'\\
&& \hspace{4.0cm} \leq \delta^N\omega_{\eta,\lambda}\left(\frac{\sqrt{N}\delta}{h} \right) \xrightarrow[h \to +\infty]{} 0,
\end{eqnarray*}
uniformly in $n \in \Nb$, for all $\eta>0$ and $\lambda>0$. Hence,
by  Fubini's Theorem
\begin{eqnarray}\label{1}
&&\int_{Q(x_0,\delta)} f_{\{\e_n\}}(x;\xi)\, dx  \nonumber\\
&&\hspace{0.0cm} \geq
\limsup_{\lambda,\eta,h,n}\frac{h^N}{\delta^N}\sum_{i=1}^{h^N}
\int_{Q_{i,h} \cap K_\eta \cap R^\lambda_n} \left\{ \int_{Q_{i,h}
\cap K_\eta} f\left(x',\frac{x}{\e_n},\frac{x}{\e_n^2};\xi +\nabla
u_n(x) \right) \, dx' \right\}\, dx.
\end{eqnarray}
\noindent However, as a consequence of  $(H_4)$ and (\ref{eta}) we
have that for all $\lambda>0$,
\begin{eqnarray}\label{2}
&&\hspace{-1cm}\frac{h^N}{\delta^N} \sum_{i=1}^{h^N} \int_{Q_{i,h}
\cap K_\eta \cap R^\lambda_n} \left\{ \int_{Q_{i,h} \setminus
K_\eta} f\left(x',\frac{x}{\e_n},\frac{x}{\e_n^2};\xi +\nabla u_n(x)
\right) \, dx' \right\}\, dx\nonumber\\
&&\leq \beta\frac{h^N}{\delta^N} \sum_{i=1}^{h^N}
\frac{\delta^N}{h^N} (1+\lambda^{p}) \mathcal L^N (Q_{i,h}\setminus K_{\eta}) \nonumber\\
&&\leq \beta   (1+\lambda^p) \mathcal L^N (\O\setminus K_{\eta})
\xrightarrow[\eta \to 0]{} 0
\end{eqnarray}
\noindent uniformly in $n \in \Nb$ and $h \in \Nb$, and similarly
\begin{equation}\label{3}
\frac{h^N}{\delta^N} \sum_{i=1}^{h^N}
\int_{Q_{i,h}  \cap R^\lambda_n \setminus K_\eta} \left\{ \int_{Q_{i,h}}
f\left(x',\frac{x}{\e_n},\frac{x}{\e_n^2};\xi +\nabla u_n(x) \right) \, dx' \right\}\, dx
\leq \beta(1+\lambda^p)\eta \xrightarrow[\eta \to 0]{} 0,
\end{equation}
uniformly in $n \in \Nb$ and $h \in \Nb$. Moreover,
(\ref{div}) and (\ref{lambda}), together with the
equi-integrability of $\{|\nabla u_n|^p\}$,  yield
\begin{eqnarray}\label{4}
&\ds \sup_{n,h,\eta}\frac{h^N}{\delta^N} \sum_{i=1}^{h^N}
\int_{Q_{i,h}  \setminus R^\lambda_n} \left\{ \int_{Q_{i,h}}
f\left(x',\frac{x}{\e_n},\frac{x}{\e_n^2};\xi +\nabla u_n(x) \right) \, dx' \right\}\, dx
\nonumber \\
&\ds \leq \sup_{n \in \Nb} \beta\int_{Q(x_0,\delta) \setminus R^\lambda_n}(1+|\nabla u_n(x)|^p)\, dx  \xrightarrow[\lambda \to +\infty]{} 0.
\end{eqnarray}
Finally, (\ref{1})-(\ref{4}) and Fubini's Theorem lead to
\begin{eqnarray}\label{fin}
&&\int_{Q(x_0,\delta)} f_{\{\e_n\}}(x;\xi)\, dx  \nonumber\\
&& \hspace{0.5cm}\geq \limsup_{h,
n}\frac{h^N}{\delta^N}\sum_{i=1}^{h^N} \int_{Q_{i,h}} \left\{
\int_{Q_{i,h}}
f\left(x',\frac{x}{\e_n},\frac{x}{\e_n^2};\xi +\nabla u_n(x) \right) \, dx\right\}\, dx' \nonumber\\
&& \hspace{0.5cm} \geq \limsup_{h \to
+\infty}\frac{h^N}{\delta^N}\sum_{i=1}^{h^N} \int_{Q_{i,h}}
\left\{ \liminf_{n \to +\infty}\int_{Q_{i,h}}
f\left(x',\frac{x}{\e_n},\frac{x}{\e_n^2};\xi +\nabla u_n(x)
\right) \, dx\right\}\, dx',
\end{eqnarray}
where we have used Fatou's Lemma. Fix $x' \in Q_{i,h}$ such that
$\overline f_{\rm hom}(x';\xi)$ is well defined and apply Theorem
\ref{jf-m2} to the continuous function $(y,z,\xi) \mapsto
f(x',y,z;\xi)$. Since $u_n \to 0$ in $L^p(Q(x_0,\delta);\Rb^d)$, we
can use the $\Gamma$-$\liminf$ inequality to get
$$\liminf_{n \to +\infty}\int_{Q_{i,h}}
f\left(x',\frac{x}{\e_n},\frac{x}{\e_n^2};\xi +\nabla u_n(x)
\right) \, dx \geq \frac{\delta^N}{h^N} \overline f_{\rm
hom}(x';\xi).$$ Then, in view of (\ref{fin}) we conclude
(\ref{009}). $\hfill\blacksquare$

\begin{lemma}\label{ineq>}
For all $\xi \in \Qb^{d \times N}$ and all $x_{0} \in L$,
$$\int_{Q(x_0,\delta)} f_{\{\e_n\}}(x;\xi)\, dx   \leq
\int_{Q(x_0,\delta)}\overline f_{\rm hom}(x;\xi)\, dx.$$
\end{lemma}

\noindent{\it Proof.} As in Lemma \ref{ineq<}, we decompose
$Q(x_0,\delta)$ into $h^N$ small disjoints cubes $Q_{i,h}$
satisfying (\ref{div}). Since $f$ and $\overline f_{\rm hom}$ are
Carath\'eodory functions, by Scorza-Dragoni's Theorem (see Ekeland
and Temam \cite{Ek&Te}) for each $\eta>0$, we can find a compact
set $K_\eta \subset Q(x_0,\delta)$ such that
\begin{equation}\label{SD}
\mathcal L^N(Q(x_0,\delta) \setminus K_\eta) < \eta,
\end{equation}
$f$ is continuous on $K_\eta \times \Rb^N \times \Rb^N \times \Rb^{d
\times N}$ and $\overline f_{\rm hom}$ is continuous on $K_\eta
\times \Rb^{d \times N}$. Let $$I_{h,\eta}:=\left\{ i\in
\{1,\cdots,h^N \}: \quad K_{\eta} \cap Q_{i,h} \neq \emptyset
\right\}.$$ For $i\in I_{h,\eta}$, choose $x_i^{h,\eta} \in K_\eta
\cap Q_{i,h}$. Theorem \ref{jf-m2}, together with e.g.\! Proposition
11.7 in Braides and Defranceschi \cite{BD},  implies the existence
of a sequence $\{u_i^{n,h,\eta}\} \subset W^{1,p}_0(Q_{i,h};\Rb^d)$
such that $u_i^{n,h,\eta} \to 0$ in $L^p(Q_{i,h};\Rb^d)$ as $n \to
+\infty$ and
$$\int_{Q_{i,h}} \overline f_{\rm hom}(x_i^{h,\eta};\xi)\, dx =
\Big(\frac{\delta}{h}\Big)^{N} \overline f_{\rm
hom}(x_i^{h,\eta};\xi) = \lim_{n \to +\infty} \int_{Q_{i,h}}f\left(
x_i^{h,\eta},\frac{x}{\e_n},\frac{x}{\e_n^2};\xi + \nabla
u_i^{n,h,\eta}\right)\, dx.$$ \noindent  Set
\begin{equation}\label{defi}u_n^\eta(x):= \left\{
\begin{array}{ll}
u_i^{n,h,\eta}(x) & \text{if }x \in Q_{i,h}\text{ and } i\in I_{h,\eta},\\
&\\
0 & \text{otherwise}.
\end{array}
\right.\end{equation} \noindent Then $\{u_n^\eta\} \subset
W^{1,p}_0(Q(x_0,\delta);\Rb^d)$, $u_n^\eta \to 0$ in
$L^p(Q(x_0,\delta);\Rb^d)$ as $n \to +\infty$ and
\begin{eqnarray}\label{037}
\liminf_{h \to +\infty}\sum_{i \in I_{h,\eta}} \int_{Q_{i,h}}
\overline f_{\rm hom}(x_i^{h,\eta};\xi)\, dx =\liminf_{h \to
+\infty}\lim_{n \to +\infty} \sum_{i \in I_{h,\eta}}
\int_{Q_{i,h}}f\left(
x_i^{h,\eta},\frac{x}{\e_n},\frac{x}{\e_n^2};\xi + \nabla
u_n^\eta\right)\, dx.
\end{eqnarray}
\noindent In view of (\ref{fhom5}) and (\ref{SD}) we have
\begin{equation}\label{039}\sup_{h \in \Nb}\sum_{i \in I_{h,\eta}}
\int_{Q_{i,h} \setminus K_\eta} \overline f_{\rm
hom}(x_i^{h,\eta};\xi)\, dx \leq \beta (1+|\xi|^p)\mathcal
L^N(Q(x_0,\delta) \setminus K_\eta) \xrightarrow[\eta \to 0]{}
0,\end{equation} \noindent thus from (\ref{037}) and  (\ref{039}) it
comes  that
\begin{eqnarray}\label{f1}
&&\hspace{-2.0cm}\liminf_{\eta, h}\sum_{i \in I_{h,\eta}} \int_{Q_{i,h} \cap K_\eta} \overline f_{\rm hom}(x_i^{h,\eta};\xi)\, dx\nonumber\\
&&\geq \liminf_{\eta, h, n} \sum_{i \in I_{h,\eta}} \int_{Q_{i,h}
\cap K_\eta}f\left(
x_i^{h,\eta},\frac{x}{\e_n},\frac{x}{\e_n^2};\xi + \nabla
u_n^\eta\right)\, dx.
\end{eqnarray}
\noindent Since $\overline f_{\rm hom}(\cdot\,;\xi)$ is continuous
on $K_\eta$, it is uniformly continuous. Thus, denoting by
$\omega_\eta$ its modulus of continuity on $K_\eta$, we have for all
$x \in Q_{i,h} \cap K_\eta$,
\begin{equation}\label{f3}
|\overline f_{\rm hom}(x;\xi)-\overline f_{\rm
hom}(x_i^{h,\eta};\xi)|\leq \omega_\eta(|x-x_i^{h,\eta}|) \leq
\omega_\eta \left(\frac{\sqrt N \delta}{h}\right) \xrightarrow[h \to
+\infty]{}0.
\end{equation}
\noindent In view of (\ref{SD}), (\ref{f3}) and (\ref{f1}), we get
since $Q_{i,h} \cap K_\eta =\emptyset$  for $i \not\in I_{h,\eta}$,
\begin{eqnarray}\label{f2}
\int_{Q(x_0,\delta)} \overline f_{\rm hom}(x;\xi)\, dx& = &\lim_{\eta \to 0} \int_{K_\eta} \overline f_{\rm hom}(x;\xi)\, dx\nonumber\\
& = &\liminf_{\eta, h}\sum_{i \in I_{h,\eta}} \int_{Q_{i,h} \cap K_\eta} \overline f_{\rm hom}(x;\xi)\, dx\nonumber\\
& = &\liminf_{\eta, h}\sum_{i \in I_{h,\eta}} \int_{Q_{i,h} \cap K_\eta} \overline f_{\rm hom}(x_i^{h,\eta};\xi)\, dx\nonumber\\
& \geq & \liminf_{\lambda, \eta, h, n} \sum_{i \in I_{h,\eta}}
\int_{Q_{i,h} \cap K_\eta \cap R_{n,\eta}^\lambda}f\left(
x_i^{h,\eta},\frac{x}{\e_n},\frac{x}{\e_n^2};\xi + \nabla
u_n^\eta\right)\, dx,
\end{eqnarray}
where $R_{n,\eta}^\lambda:= \left\{ x \in Q(x_0,\delta) : |\xi
+\nabla u_n^\eta(x)|\leq \lambda \right\}$. From (\ref{037}) and the
fact that $u_n^\eta \equiv 0$ on $Q(x_0,\delta) \setminus \bigcup_{i
\in I_{h,\eta}} Q_{i,h}$, we get
\begin{equation}\label{1853}
\sup_{n \in \Nb,\eta>0} \int_{Q(x_0,\delta)}|\nabla u^\eta_n|^p dx
<+ \infty.
\end{equation}
In particular, according to Chebyshev's
inequality, we have
\begin{equation}\label{cheb}
\mathcal L^N (Q(x_0,\delta) \setminus R_{n,\eta}^\lambda) \leq
\frac{C}{\lambda^p},
\end{equation}
for some constant $C>0$ independent of $n$, $\eta$ and $\lambda$.
Since $f$ is continuous on $K_\eta \times \Rb^N \times \Rb^N \times
\Rb^{d \times N}$ and separately $Q$-periodic in its second and
third variable (see assumption $(H_3)$), it is uniformly continuous
on $K_\eta \times \Rb^N \times \Rb^N \times \overline B(0,\lambda)$.
Thus, denoting by $\omega_{\eta,\lambda}$ its modulus of continuity
on $K_\eta \times \Rb^N \times \Rb^N \times \overline B(0,\lambda)$,
we have for all $x \in Q_{i,h} \cap K_\eta \cap R_{n,\eta}^\lambda$,
\begin{eqnarray*}
&&\left| f\left( x,\frac{x}{\e_n},\frac{x}{\e_n^2};\xi + \nabla u_n^\eta(x) \right)-f\left( x_i^{h,\eta},\frac{x}{\e_n},\frac{x}{\e_n^2};\xi + \nabla u_n^\eta(x) \right) \right|\\
&&\hspace{2.0cm} \leq \omega_{\eta,\lambda} (|x-x_i^{h,\eta}|)\\
&&\hspace{2.0cm} \leq \omega_{\eta,\lambda} \left(\frac{\sqrt N
\delta}{h}\right) \xrightarrow[h \to +\infty]{} 0.
\end{eqnarray*}
Then, according to (\ref{f2}) and the fact that $Q_{i,h} \cap K_\eta =\emptyset$ for $i \not\in I_{h,\eta}$,
\begin{eqnarray*}
\int_{Q(x_0,\delta)} \overline f_{\rm hom}(x;\xi)\, dx & \geq &
\liminf_{\lambda, \eta, h, n} \sum_{i \in I_{h,\eta}} \int_{Q_{i,h}
\cap K_\eta \cap R_{n,\eta}^\lambda}f\left(
x,\frac{x}{\e_n},\frac{x}{\e_n^2};\xi + \nabla
u_n^\eta\right)\, dx,\\
& = & \liminf_{\lambda, \eta, n}  \int_{K_\eta \cap
R_{n,\eta}^\lambda}f\left( x,\frac{x}{\e_n},\frac{x}{\e_n^2};\xi +
\nabla u_n^\eta\right)\, dx.
\end{eqnarray*}
In view of the $p$-growth condition $(H_4)$, (\ref{SD}) and
the definition of $R_{n,\eta}^\lambda$,
$$\sup_{n \in \Nb} \int_{R_{n,\eta}^\lambda \setminus K_\eta}f\left( x,\frac{x}{\e_n},\frac{x}{\e_n^2};\xi + \nabla u_n^\eta\right)\, dx \leq \beta(1+\lambda^p) \eta \xrightarrow[\eta \to 0]{} 0,$$
then
$$\int_{Q(x_0,\delta)} \overline f_{\rm hom}(x;\xi)\, dx \geq \liminf_{\lambda, \eta, n}
\int_{R_{n,\eta}^\lambda}f\left(
x,\frac{x}{\e_n},\frac{x}{\e_n^2};\xi + \nabla u_n^\eta\right)\,
dx.$$ Let $\lambda_k \nearrow +\infty$ and $\eta_k \searrow 0^+$,
by a diagonalization procedure, it is possible to find a
subsequence $\{n_k\}$ of $\{n\}$ such that, upon setting $v_k :=
u_{n_k}^{\eta_k}$ and $R_k := R^{\lambda_k}_{n_k,\eta_k}$, then
$v_k \in W^{1,p}_0(Q(x_0,\delta);\Rb^d)$, $v_k \to 0$ in
$L^p(Q(x_0,\delta);\Rb^d)$ and
$$\int_{Q(x_0,\delta)} \overline f_{\rm hom}(x;\xi)\, dx \geq \liminf_{k \to +\infty}
\int_{R_k}f\left(x,\frac{x}{\e_{n_k}},\frac{x}{\e_{n_k}^2};\xi +
\nabla v_k \right)\, dx.$$ By (\ref{1853}) and the Poincar\'e
Inequality, the sequence $\{v_k\}$ is bounded in
$W^{1,p}(Q(x_0,\delta);\Rb^d)$  uniformly with respect to $k \in
\Nb$ so that, according to the Decomposition Lemma (see Fonseca,
M\"uller and Pedregal \cite{FMP} or Fonseca and Leoni \cite{FL}),
there is no loss of generality to assume that $\{|\nabla v_k|^p\}$
is equi-integrable. It turns out, in view of the $p$-growth
condition $(H_4)$ and (\ref{cheb}) that
$$\int_{Q(x_0,\delta) \setminus R_k}
f\left( x,\frac{x}{\e_{n_k}},\frac{x}{\e_{n_k}^2};\xi + \nabla v_k
\right)\, dx \leq \beta \sup_{l \in \Nb} \int_{Q(x_0,\delta)
\setminus R_k}(1+|\nabla v_l|^p)\, dx \xrightarrow[k \to
+\infty]{}0.$$ Thus, using the $\G$-$\liminf$ inequality,
\begin{eqnarray*}
\int_{Q(x_0,\delta)}  \overline f_{\rm hom}(x;\xi)\, dx & \geq & \liminf_{k \to +\infty}
\int_{Q(x_0,\delta)}  f\left( x,\frac{x}{\e_{n_k}},\frac{x}{\e_{n_k}^2};\xi + \nabla
v_k \right)\, dx\\
& \geq & \int_{Q(x_0,\delta)}  f_{\{\e_n\}}(x;\xi)\, dx.
\end{eqnarray*}
$\hfill\blacksquare$\\

\noindent{\it Proof of Theorem \ref{jf-m1}. }As a consequence of
Lemma \ref{ineq<} and \ref{ineq>}, we have  $\overline f_{\rm
hom}(x;\xi)=f_{\{\e_n\}}(x;\xi)$ for all $x \in L$ and all $\xi
\in \Qb^{d \times N}$. By Lemma \ref{intrep2} and the fact that
$\overline f_{\rm hom}$ is (equivalent to) a Carath\'eodory
function,  it follows that the equality holds for all $\xi \in
\Rb^{d \times N}$ and a.e. $x \in \O$. Therefore, we have
$\mathcal F_{\{\e_n\}}(u;A)=\mathcal F_{\rm hom}(u;A)$ for all $A
\in \mathcal A(\O)$ and all $u \in W^{1,p}(A;\Rb^d)$. Since the
result does not depend upon the specific choice of the
subsequence, we conclude thanks to Proposition 8.3 in Dal Maso
\cite{DM} that the whole sequence $\mathcal F_\e(\cdot\,;A)$
$\G(L^p(A))$-converges to $\mathcal F_{\rm hom}(\cdot\,;A)$.
Taking $A=\O$ we conclude the proof of Theorem \ref{jf-m1}.
$\hfill\blacksquare$


\subsection{Some remarks in the convex case}\label{cx-subsection}

\noindent We start this section by noticing  that  under the
additional hypothesis that $f(x,y,z, \cdot)$ is convex for a.e.\!
$x$ and all $(y,z)$, in which case $(H_{1})$ is equivalent to
requiring that $f(x,\cdot,\cdot,\xi)$ is continuous for a.e.\! $x$
and all $\xi$, equality  (\ref{fhombarra}) and (\ref{fhom})
simplify to read
$$\overline f_{\rm hom}(x;\xi)=\inf_{\phi} \left\{\int_{Q} f_{\rm
hom}(x,y;\xi+\nabla \phi(y))\, dy, \quad \phi \in W_{0}^{1,p}(Q
;\Rb^{d})~~ \big(\text{or equivalently}~\phi \in W_{\rm
per}^{1,p}(Q ;\Rb^{d})\big)\right\}$$ \noindent for all $\xi\in
\Rb^{d\times N}$ and a.e. $x \in \Omega$, and
$$ f_{\rm hom}(x,y;\xi)= \inf_{\phi} \left\{
\int_{Q}  f(x,y,z;\xi+\nabla \phi(z))\, dz, \quad \phi \in
W_{0}^{1,p}(Q ;\Rb^{d})~~ \big(\text{or equivalently}~\phi \in
W_{\rm per}^{1,p}(Q ;\Rb^{d})\big)\right\}$$ \noindent  for a.e
$x\in \O$ and all $(y,\xi)\in \Rb^N \times \Rb^{d\times N}$ (see
M\"{u}ller \cite {Mu} and Braides and Defranceschi \cite{BD}).

Our objective here is to  present an alternative proof of Lemma
\ref{ineq>}  in the convex case.  Namely we would like to show that
$ f_{\{\e_n\}}(x_{0};\xi) \leq \overline f_{\rm hom}(x_{0};\xi)$
a.e.\! $x_{0}\in \O$ and all $\xi\in \Rb^{d\times N}$, without
appealing to Theorem \ref{jf-m2}. For this purpose let us denote by
${\cal S}$ (resp.\! ${\cal C}$) a countable set of  functions in
$\mathcal C_{c}^{\infty}(Q;\Rb^{d})$ (resp.\! $\mathcal
C_{c}^{\infty}(Q \times Q;\Rb^{d})$) dense in
$W_{0}^{1,p}(Q;\Rb^{d})$ (resp.\!
$L^{p}(Q;W_{0}^{1,p}(Q;\Rb^{d}))$). Define $L$ to be the set of
Lebesgue points $x_{0}$ for all functions
\begin{equation}\label{C1)}
f_{\{\e_n\}}(\cdot\,;\xi), \quad \overline f_{\rm
hom}(\cdot\,;\xi)
\end{equation}
\noindent and
\begin{equation}\label{C2)}
x\rightarrow \int_{Q}\int_{Q} f(x,y,z;\xi+\nabla_{y}
\phi(y)+\nabla_{z} \psi (y,z))\, dy\, dz,
\end{equation}
\noindent with  $\phi \in {\cal S}$, $\psi \in \mathcal C$ and $\xi
\in \Qb^{d\times N}$, and for which $\overline{f}_{\rm
hom}(x_{0};\,\cdot\,)$  is well defined. Note that  $\mathcal
L^N(\O\setminus L)=0$. Let $x_{0}\in L$ and $\xi\in \Qb^{d\times
N},$ then
\begin{equation}\label{lpg}f_{\{\e_n\}}(x_{0};\xi) = \lim_{\delta \to  0}
\frac{1}{\delta^{N}}\int_{Q(x_{0},\delta)}f_{\{\e_n\}}(x;\xi)\,dx =
\ds \lim_{\delta \to  0}\frac{{\cal F}_{\{\e_n\}}(\xi
\,\cdot\,;Q(x_{0},\delta))}{\delta^{N}}.
\end{equation}
\noindent Given $m \in \Nb$ consider  $\phi_m \in {\cal S}$  such
that
\begin{equation}\label{desdens}\overline{f}_{\rm
hom}(x_{0};\xi)+\frac{1}{m} \geq \int_{Q} f_{\rm
hom}(x_{0},y;\xi+\nabla {\phi}_m (y))\, dy.
\end{equation}
\noindent Then by Theorem III.30 in Castaing and Valadier \cite{CV}
and following a similar argument as in Lemma 4.6 in Fonseca and
Zappale \cite{FZ}, there exist $\Phi_{m}\in
L^{p}\big(Q;W_{0}^{1,p}(Q;\Rb^{d})\big)$ such that
$$f_{\rm hom}(x_{0},y;\xi+\nabla \phi_{m}(y))+\frac{1}{m}\geq \int_{Q}
f(x_{0},y,z;\xi+\nabla_{y} \phi_m(y)+ \nabla_{z} \Phi_{m} (y,z))\,
dz.$$ \noindent We now choose $\Phi_{m,k}\in \mathcal C$  such that
\begin{equation}\label{conv}
\|\Phi_{m,k}-\Phi_{m}\|_{L^{p}(Q;W^{1,p}_0(Q;\Rb^d))}\xrightarrow[k
\to +\infty]{} 0,
\end{equation}
\noindent and we extend $\phi_{m}$ and $\Phi_{m,k}$ periodically to
$\Rb^{N}$  and $\Rb^{N}\times \Rb^{N}$, respectively. For each $x\in
\Rb^{N}$  define
$$u^{n}_{m,k}(x):=\xi\cdot x+\e_{n}\phi_{m}\Big(\frac{x}{\e_{n} }\Big)+\e_{n}^{2}
\Phi_{m,k}\Big(\frac{x}{\e_{n}}, \frac{x}{\e_{n}^{2}}\Big)$$
\noindent  and consider  $\delta>0$  small enough so that
$Q(x_{0},\delta)\in {\mathcal A}(\O).$  For fixed $m$ and $k$ we
have $u^{n}_{m,k} \to  v$ in $L^{p}(Q(x_{0};\delta);\Rb^{d})$ as $n
\to +\infty$, where $v(x)=\xi\cdot x$.  Hence by (\ref{lpg}) and the
$p$-Lipschitz property of $f(x,y,\cdot)$ (see Marcellini \cite{Ma1})
\begin{eqnarray}\label{igua1}
f_{\{\e_n\}}(x_{0};\xi) &\leq & \debaixodoliminf {} {k, \delta, n
} \frac{1}{\delta^{N}}\int_{Q(x_{0};\delta)}
f\left(x,\frac{x}{\e_{n}},\frac{x}{\e_{n}^{2}};\xi+\nabla_{y}\phi_{m}\left(\frac{x}{\e_{n}}\right)
\right.\nonumber\\
&& \hspace{3.0cm}\left. + \e_n \nabla_y \Phi_{m,k}\Big(\frac{x}{\e_{n}},
\frac{x}{\e_{n}^{2}}\Big) + \nabla_{z}\Phi_{m,k}\Big(\frac{x}{\e_{n}},
\frac{x}{\e_{n}^{2}}\Big)\right)\, dx\nonumber\\
&\leq  & \debaixodoliminf {} {k, \delta, n }
\frac{1}{\delta^{N}}\int_{Q(x_{0};\delta)}
f\left(x,\frac{x}{\e_{n}},\frac{x}{\e_{n}^{2}};\xi+\nabla_{y}\phi_{m}\left(\frac{x}{\e_{n}}\right)
\right.\nonumber\\
&&\hspace{3.0cm}\left.+\nabla_{z}\Phi_{m,k}\Big(\frac{x}{\e_{n}},\frac{x}{\e_{n}^{2}}\Big)\right)\, dx.
\end{eqnarray}
\noindent Arguing as in Proposition 4.10 in  Ba\'{\i}a and Fonseca
\cite{B&Fo}, we define
$$h_{m,k}(x,y,z):= f\big(x,y,z; \xi+\nabla_{y}\phi_{m}(y)
+\nabla_{z}\Phi_{m,k}(y,z)\big).$$ \noindent Then (see e.g.\!
Allaire and Briane \cite{AB} or Donato \cite{Do})  since $h_{m,k}
\in L^p(Q(x_0,\delta); \mathcal C_{\rm per}(Q \times Q))$, we get
\begin{eqnarray}\label{igua2} &&\debaixodoliminf {}{n \to  +\infty}
\int_{Q(x_{0};\delta)}
f\left(x,\frac{x}{\e_{n}},\frac{x}{\e_{n}^{2}};\xi+\nabla_{y}\phi_{m}\left(\frac{x}{\e_{n}}\right)
+\nabla_{z}\Phi_{m,k}\Big(\frac{x}{\e_{n}},
\frac{x}{\e_{n}^{2}}\Big)\right)\, dx\nonumber\\
&&\ds= \debaixodolim {} {n \to  +\infty}\int_{Q(x_{0};\delta)}
h_{m,k}\left(x,\frac{x}{\e_{n}},\frac{x}{\e_{n}^{2}}\right)\, dx\nonumber\\
&&\ds =\int_{Q(x_{0};\delta)}\int_{Q} \int_{Q} h_{m,k}(x,y,z)\, dz\,dy\, dx\nonumber\\
& &\ds =\int_{Q(x_{0};\delta)}\int_{Q} \int_{Q} f\big(x,y,z;
\xi+\nabla_{y}\phi_{m}\left(y\right) +\nabla_{z}\Phi_{m,k}(y,
z)\big)dz\, dy\, dx .\end{eqnarray} \noindent Therefore by
(\ref{C2)})
$$\begin{array}{ll}
&\ds \debaixodoliminf {}{\delta \to  0}\debaixodoliminf {}{n \to  +
\infty} \frac{1}{\delta^{N}}\int_{Q(x_{0};\delta)}
f\left(x,\frac{x}{\e_{n}},\frac{x}{\e_{n}^{2}};\xi+\nabla_{y}\phi_{m}\left(\frac{x}{\e_{n}}\right)
+\nabla_{z}\Phi_{m,k}\Big(\frac{x}{\e_{n}},
\frac{x}{\e_{n}^{2}}\Big)\right)dz\, dy\, dx\\\\
&\ds  = \int_{Q} \int_{Q} f(x_{0}, y,  z; \xi+\nabla_{y}\phi_{m}(y)+
\nabla_{z}\Phi_{m,k}(y,z))\,dz\,dy,
\end{array}$$
\noindent and thus, by (\ref{conv}), $(H_{4})$, (\ref{igua1}),
(\ref{desdens}) and Fubini's Theorem, we obtain
\begin{eqnarray*}
f_{\{\e_n\}}(x_{0};\xi)&\leq &\int_{Q}
\int_{Q}f(x_{0},y,z;\xi+\nabla_{y} {\phi}_m
(y)+\nabla_{z}\Phi_{m}(y,z))\,dz\,dy\\
&=& \int_{Q} \left[\int_{Q} f(x_{0},y,z;\xi+\nabla_{y} {\phi}_m (y)+\nabla_{z}\Phi_{m}(y,z))\,dz\right] \,dy\\
&\leq& \int_{Q} f_{\rm
hom}(x_{0},y;\xi+\nabla \phi_{m}(y))\,dy +\frac{1}{m}\\
&\leq& \overline{f}_{\rm hom}(x_{0};\xi)+\frac{2}{m}.
\end{eqnarray*}
\noindent Letting $m \to  + \infty$ we deduce that
$f_{\{\e_n\}}(x_{0};\xi)\leq \overline{f}_{\rm
 hom}(x_{0};\xi)$.
$\hfill\blacksquare$


\section{Application to thin films}\label{part3}

\noindent This part is devoted to the study of a reiterated
homogenization  problem in the framework of 3D-2D dimensional
reduction. Our main result is stated in  Theorem \ref{jf-m3}.  We
organize this section as follows.  In Subsection
\ref{continuity3d2d} we discuss the main properties of $W_{\rm
hom}$ and $\overline{W}_{\rm hom}$.  Then, in Subsection
\ref{3homo} we address the case where $W$ is independent of the
macroscopic in-plane variable  $x_{ \alpha}$ (Theorem
\ref{jf-m4}). Finally,  Theorem \ref{jf-m3} is proved in
Subsection \ref{3hetero}.

\begin{rmk}{\rm
Without loss of generality, we  assume that $W$ is non negative
upon replacing $W$ by $W+\beta$ which is non negative in view of
$(A_4)$. }
\end{rmk}


\subsection{Properties of $\overline{W}_{\rm hom}$}\label{continuity3d2d}

\noindent  As in Section \ref{continuity} we can see that the
function $W_{\rm hom}$ given in (\ref{Whom}) is well defined and
it is (equivalent to) a Carath\'eodory function:
\begin{eqnarray}
&&W_{\rm hom}(\cdot,\cdot\,;\xi) \text{ is }\mathcal L^3 \otimes \mathcal L^2\text{-measurable for all  }\xi \in \Rb^{3 \times 3}, \label{Whom1}\\
&&W_{\rm hom}(x,y_\a;\cdot) \text{ is continuous for }\mathcal L^3
\otimes \mathcal L^2\text{-a.e. }(x,y_\a) \in \Omega \times \Rb^2.
\label{Whombis2}
\end{eqnarray}
By condition $(A_3)$ it follows that
\begin{equation}\label{Whom3}
W_{\rm hom}(x,\cdot\,;\xi) \text{ is }Q'\text{-periodic for a.e.
}x \in \O \text{ and all }\xi \in \Rb^{3 \times 3}.
\end{equation}
\noindent Moreover, $W_{\rm hom}$ is  quasiconvex in the $\xi$
variable and satisfies the same $p$-growth and $p$-coercivity
condition $(A_4)$ as $W$:
\begin{equation}\label{Whom4}
\frac{1}{\beta}|\xi|^p-\beta\leq W_{\rm hom}(x,y_\a;\xi) \leq
\beta(1+|\xi|^p)\quad \text{for a.e. }x \in \O\text{ and all
}(y_\a,\xi) \in \Rb^{2} \times \Rb^{3\times 3},
\end{equation}
\noindent where $\beta$ is the constant in $(A_{4})$. Arguing as in
Remark 2.2 in Babadjian and Ba\'{\i}a \cite{BB1}, (\ref{Whom1}),
(\ref{Whombis2}) and (\ref{Whom4}) imply that the function
$\overline W_{\rm hom}$ given in (\ref{Whombarra}) is also well
defined, and is (equivalent to) a Carath\'eodory function, which
implies that the definition of $\W_{\rm hom}$ makes sense on
$W^{1,p}(\Omega;\Rb^3)$. Finally, $\overline W_{\rm hom}$ is also
quasiconvex in the $\overline \xi$ variable and satisfies the same
$p$-growth and $p$-coercivity condition $(A_4)$ as $W$ and $W_{\rm
hom}$:
\begin{equation}\label{Whom5}
\frac{1}{\beta}|\overline \xi|^p-\beta \leq \overline W_{\rm
hom}(x_\a;\overline \xi) \leq \beta(1+|\overline \xi|^p)\quad
\text{for a.e. }x_\a \in \o\text{ and all }\overline \xi \in
\Rb^{3\times 2},
\end{equation}
\noindent where, as before, $\beta$ is the constant in $(A_{4})$.


\subsection{Independence of the in-plane macroscopic variable}\label{3homo}

\noindent In this section, we assume that $W$   does not depend
explicitly on $x_\a$, namely  $W: I \times \Rb^{3}\times
\Rb^{2}\times \Rb^{3\times 3}\rightarrow \Rb^{+}.$ For each
$\e>0$, consider the functional $\W_{\e}: L^{p}(\O;\Rb^{3})
\rightarrow [0,+\infty]$ defined by
\begin{equation}\label{Wepsilon2}
\W_{\e}(u):=\left\{\begin{array}{ll} \ds \int_{\O}
W\left(x_3,\frac{x}{\e},\frac{x_\a}{\e^{2}};\nabla_\a
 u(x)\Big|\frac{1}{\e} \nabla_3 u(x) \right)\, dx & \text{if } u  \in W^{1,p}(\O;\Rb^{3}),
\\
+\infty & \text{otherwise}.
\end{array}\right.
\end{equation}
\noindent Our objective is to prove the following result.
\begin{thm}\label{jf-m4}
Under assumptions  $(A_1)$-$(A_4)$ the  $\Gamma (L^p(\O))$-limit
of the family $\{\W_\e \}_{\e>0}$ is given by
$$\W_{\rm hom}(u) =
\left\{\begin{array}{ll} \ds 2\int_{\o} \overline W_{\rm
hom}(\nabla_\a u(x_\a))\, dx_\a & \text{if } u  \in
W^{1,p}(\o;\Rb^{3}),\\
+ \infty &\text{otherwise},
\end{array}\right.$$
\noindent where $\overline W_{\rm hom}$ is defined, for all
$\overline \xi\in \Rb^{3\times 2}$, by
\begin{eqnarray}\label{Whombarra2}
\overline W_{\rm hom}(\overline \xi) & := & \lim_{T \to  + \infty}
\inf_{\phi} \Big\{ \frac{1}{2T^{2}}\int_{(0,T)^{2} \times I} W_{\rm
hom}(y_3,y_\a;\overline \xi+\nabla_\a \phi(y)|\nabla_3 \phi(y))\,
dy:\nonumber\\
&&\hspace{2.0cm}\phi \in W^{1,p}((0,T)^2 \times I;\Rb^3) \text{ and
} \phi=0 \text{ on }\partial (0,T)^2 \times I \Big\},\end{eqnarray}
\noindent and
\begin{eqnarray}\label{Whombarra21}
W_{\rm hom}(y_3,y_\a;\xi) & := & \lim_{T \to  + \infty}
\inf_{\phi} \Big\{ \frac{1}{T^{3}}\int_{(0,T)^{3}}
W(y_3,y_\a,z_3,z_\a;\xi+\nabla \phi(z))\,dz:\nonumber\\
&&\hspace{4.0cm} \phi \in W_{0}^{1,p}((0,T)^{3} ;\Rb^{3})\Big\},
\end{eqnarray}
for all $(y,\xi)\in \Rb^3 \times \Rb^{3\times 3}$.
\end{thm}
Since the proofs are very similar to that of Section \ref{2homo}, we
just sketch them highlighting the main differences.  For the
detailed proofs we refer to Babadjian \cite{Bab} Chapter 2.


\subsubsection{Existence and integral representation of the $\G$-limit}

\noindent For the same reason than in the proof of Theorem
\ref{jf-m2} in Section \ref{2homo}, we localize the functionals
given in (\ref{Wepsilon2}) on  the class of bounded open subsets
of $\Rb^2$, denoted by $\mathcal A_0$. For each $\e>0$, consider
$\W_{\e}: L^{p}(\Rb^2 \times I;\Rb^{3}) \times {\cal A}_0
\rightarrow [0,+\infty]$ defined by
\begin{equation}\label{We2}
\W_{\e}(u;A):=\left\{\begin{array}{ll} \ds \int_{A \times I}
W\left(x_3,\frac{x}{\e},\frac{x_\a}{\e^{2}};\nabla_\a
 u(x)\Big|\frac{1}{\e} \nabla_3 u(x) \right)\, dx & \text{if } u \in W^{1,p}(A \times I;\Rb^{3}),
\\
+\infty & \text{otherwise}.
\end{array}\right.
\end{equation}
Given $\{\e_j\} \searrow 0^+$  and $A \in \mathcal A_0$, consider
the $\Gamma$-lower limit  of $\{  \W_{\e_j} (\cdot\,;A) \}_{j \in
\Nb}$ for the $L^{p}(A \times I;\Rb^{3})$-topology, defined  for $u
\in L^{p}(\Rb^2 \times I;\Rb^3)$, by
$$\W_{\{\e_j\}}(u;A) := \inf_{\{u_j\}}\left\{\liminf_{j \to
+\infty} \W_{\e_j}(u_j;A) : \quad u_j \to u \text{ in }L^p(A
\times I;\Rb^3)\right\}.$$ \noindent In view of the $p$-coercivity
condition $(A_{4})$, for each $A \in \mathcal A_0$ it follows that
$\W_{\{\e_j\}}(u;A)$ is infinite whenever $u \in L^p(\Rb^2 \times
I;\Rb^3) \setminus W^{1,p}(A;\Rb^3)$, so it suffices  to study the
case where $u \in W^{1,p}(A;\Rb^3)$.  Arguing  as in the proof of
Lemma 2.6 in Braides, Fonseca and Francfort \cite{BFF}, we can
prove the existence of a subsequence $\{\e_{j_{n}}\}
\equiv\{\e_{n}\}$ such that $\W_{\{\e_n\}}(\cdot\,;A)$ is the
$\G(L^p(A \times I)$-limit of $\{\W_{\e_n} (\cdot\,;A) \}_{n \in
\Nb}$ for each $A \in \mathcal A_0$.  In addition,  following the
lines of Lemma 2.5 in Braides, Fonseca and Francfort \cite{BFF},
it is possible to show the following result.

\begin{lemma}
For each $A \in \mathcal A_0$ and all $u \in W^{1,p}(A;\Rb^3)$, the
restriction of $\W_{\{\e_n\}}(u;\cdot)$ to $\mathcal A(A)$ is a
Radon measure, absolutely continuous with respect to the
two-dimensional Lebesgue measure.
\end{lemma}

But as in Section \ref{e-gl}, one has to ensure that the integral
representation given by Theorem 1.1 in Buttazzo-Dal Maso \cite{DMB}
is independent of the open set $A \in \mathcal A_0$. The following
result, prevents this dependence from holding since it leads to an
homogeneous integrand as it will be seen in Lemma \ref{intrep3d2d}
below.

\begin{lemma}\label{independency3d2d}
For all $\overline \xi \in \Rb^{3\times 2}$, $y_\a^0$ and $z_\a^0\in \Rb^2$, and
$\delta>0$
$$\W_{\{\e_{n}\}}(\overline \xi\,\cdot\,;Q'(y_\a^0,\delta))=\W_{\{\e_{n}\}}(\overline \xi\,\cdot\,;
Q'(z_\a^0,\delta)).$$
\end{lemma}

\noindent {\it Proof. }It is obviously enough to show that
$$\W_{\{\e_{n}\}}(\overline \xi\,\cdot\,;Q'(y_\a^0,\delta))\geq \W_{\{\e_{n}\}}(\overline \xi\,\cdot\,;
Q'(z_\a^0,\delta)).$$ \noindent According to Theorem 1.1 in Bocea
and Fonseca \cite{Bo&Fo} together with Lemma 2.6 in Braides, Fonseca
and Francfort \cite{BFF}, there exists a sequence $\{u_{n}\}\subset
W^{1,p}(Q'(y_\a^0,\delta) \times I;\Rb^{3})$ such that $\big\{\big|
\big(\nabla_\a u_n |\frac{1}{\e_n} \nabla_3 u_n \big)\big|^p \big\}$
is equi-integrable, $u_n=0$ on $\partial Q'(y_\a^0,\delta) \times
I$, $u_{n}\rightarrow 0$ in $L^{p}(Q'(y_\a^0,\delta) \times
I;\Rb^{3})$ and
$$\W_{\{\e_{n}\}}(\overline \xi\cdot\,;Q'(y_\a^0,\delta))=\lim_{n \to
+\infty}\int_{Q'(y_\a^0,\delta) \times
I}W\left(x_3,\frac{x}{\e_{n}}, \frac{x_\a}{\e_{n}^{2}};\overline
\xi+\nabla_\a u_{n}(x)\Big|\frac{1}{\e_n} \nabla_3 u_n(x)
\right)dx.$$ \noindent We argue exactly as in the proof of Lemma
\ref{independency} with $y_\a^0$ and $z_\a^0$ in place of $y_0$ and
$z_0$. For all $n\in \Nb$, extend $u_n$ by zero to the whole $\Rb^2
\times I$ and set $v_n(x_\a,x_3)=u_n(x_\a+x_\a^{\e_n},x_3)$ for
$(x_\a,x_3) \in Q'(z_\a^0,\delta) \times I$, where
$x_\a^{\e_{n}}:=m_{\e_{n}}\e_{n}-\e_{n}^{2}l_{\e_{n}}$. Then
$\{v_n\} \subset W^{1,p}(Q'(z_\a^0,\delta) \times I;\Rb^3)$, $v_n
\to 0$ in $L^p(Q'(z_\a^0,\delta) \times I;\Rb^3)$, the sequence
$\big\{\big| \big(\nabla_\a v_n \big|\frac{1}{\e_n} \nabla_3 v_n
\big)\big|^p \big\}$ is equi-integrable and
\begin{eqnarray}\label{1239}
&&\W_{\{\e_{n}\}}(\overline \xi\cdot\,;Q'(y_\a^0,\delta))\nonumber\\
&&\hspace{1.0cm} = \limsup_{n \to +\infty}\int_{Q'(z_\a^0,\delta)
\times I}W\left(x_3,\frac{x_\a}{\e_{n}} - \e_n l_{\e_n},
\frac{x_3}{\e_n},\frac{x_\a}{\e_{n}^{2}}; \overline \xi+\nabla_\a
v_{n}(x) \Big|\frac{1}{\e_n}\nabla_3 v_{n}(x)\right)dx,
\end{eqnarray}
\noindent where we have used the $p$-growth condition $(A_4)$ and
the fact that $\mathcal L^2(Q'(z_\a^0,\delta) \setminus
Q'(y_\a^0-x_\a^{\e_n},\delta)) \to 0$. To eliminate the term $\e_n
\, l_{\e_n}$ in (\ref{1239}), we would like to apply a uniform
continuity argument. Since for a.e. $x_3 \in I$ the function
$W(x_3,\cdot,\cdot\,;\cdot)$ is continuous on $\Rb^3 \times \Rb^2
\times \Rb^{3 \times 3}$, then $(A_3)$ implies that it is uniformly
continuous on $\Rb^3 \times \Rb^2 \times \overline{B}(0,\lambda)$
for any $\lambda>0$. We define
$$R^\lambda_n:= \left\{ x \in Q'(z_\a^0,\delta) \times I : \left|\left(\overline \xi  + \nabla_\a
v_n(x)\Big|\frac{1}{\e_n} \nabla_3 v_n(x) \right) \right| \leq
\lambda\right\},$$ \noindent and we note that by Chebyshev's
inequality
\begin{equation}\label{che-13d2d}\mathcal L^3([Q'(z_\a^0,\delta) \times I] \setminus R^\lambda_n) \leq
C/\lambda^p,\end{equation} \noindent for some constant $C>0$
independent of $\lambda$ or $n$. Thus, in view  of (\ref{1239}) and
the fact that $W$ is nonnegative,
\begin{eqnarray*}
\W_{\{\e_{n}\}}(\overline \xi\cdot\,;Q'(y_\a^0,\delta)) \geq
\limsup_{\lambda,
n}\int_{R^\lambda_n}W\left(x_3,\frac{x_\a}{\e_{n}}-\e_n
l_{\e_n},\frac{x_3}{\e_n},\frac{x_\a}{\e_n^2};\overline
\xi+\nabla_\a v_{n}(x)\Big|\frac{1}{\e_n} \nabla_3 v_n(x)
\right)\,dx.
\end{eqnarray*}
\noindent Denoting by $\omega_\lambda(x_3,\cdot) : \Rb^+ \to \Rb^+$
the modulus of continuity of $W(x_3,\cdot,\cdot\,;\cdot)$ on $\Rb^3
\times \Rb^2 \times \overline{B}(0,\lambda)$, we can check that for
a.e. $x_3 \in I$, the function $t \mapsto \o_\lambda(x_3,t)$ is
continuous, increasing and satisfies $\o_\lambda(x_3,0)=0$ while,
for all $t \in \Rb^+$, the function $x_3 \mapsto \o_\lambda(x_3,t)$
is measurable (as the supremum of measurable functions). We get, for
any $x \in R^\lambda_n$
\begin{eqnarray*}
&&\left| W\left(x_3,\frac{x}{\e_{n}},
\frac{x_\a}{\e_{n}^{2}};\overline\xi+\nabla_\a
v_{n}(x)\Big|\frac{1}{\e_n} \nabla_3 v_n(x) \right)\right. \\
&&\hspace{5.0cm} \left. -W\left(x_3,\frac{x_\a}{\e_{n}}-\e_n
l_{\e_n},\frac{x_3}{\e_n},
\frac{x_\a}{\e_{n}^{2}};\overline\xi+\nabla_\a
v_{n}(x)\Big|\frac{1}{\e_n} \nabla_3 v_n(x) \right)\right|\\
&&\hspace{2.5cm} \leq \omega_\lambda(x_3,\e_n l_{\e_n}).
\end{eqnarray*}
The properties of $\o_\lambda$, Beppo-Levi's Monotone Convergence
Theorem and (\ref{1239}) yield
\begin{eqnarray*}
\W_{\{\e_{n}\}}(\overline \xi\cdot\,;Q'(y_\a^0,\delta))&\geq&
\limsup_{\lambda, n}\left\{
\int_{R^\lambda_n}W\left(x_3,\frac{x}{\e_{n}},
\frac{x_\a}{\e_{n}^{2}};\overline\xi+\nabla_\a
v_{n}(x)\Big|\frac{1}{\e_n} \nabla_3 v_n(x) \right)dx \right.\\
&&\hspace{4.0cm} \left. - \delta^2 \int_{-1}^1 \omega_\lambda(x_3,\e_n l_{\e_n})\, dx_3\right\}\\
&  = &  \liminf_{n \to +\infty}\int_{Q'(z_\a^0,\delta) \times
I}W\left(x_3,\frac{x}{\e_{n}},
\frac{x_\a}{\e_{n}^{2}};\overline\xi+\nabla_\a
v_{n}(x)\Big|\frac{1}{\e_n} \nabla_3 v_n(x) \right)dx\\
& \geq & \W_{\{\e_{n}\}}(\overline \xi\cdot\,;Q'(z_\a^0,\delta)),
\end{eqnarray*}
where we have used the equi-integrability of $\big\{\big|
\big(\nabla_\a v_n \big|\frac{1}{\e_n} \nabla_3 v_n \big)\big|^p
\big\}$, the $p$-growth condition $(A_4)$, (\ref{che-13d2d}) and the
fact that $v_n \to 0$ in $L^p(Q'(z_\a^0,\delta) \times I;\Rb^3)$.
$\hfill\blacksquare$\\

As a consequence of this lemma and adapting the argument used in the
proof of Lemma \ref{intrep}, we deduce that

\begin{lemma}\label{intrep3d2d}
There exists a continuous function $W_{\{\e_n\}}: \Rb^{3 \times 2} \rightarrow \Rb^+$
such that for all $A \in \mathcal A_0$ and all $u \in W^{1,p}(A;\Rb^3)$,
$$\W_{\{\e_{n}\}} (u;A)=2\int_{A} W_{\{\e_n\}}(\nabla_\a u(x_\a))\, dx_\a.$$
\end{lemma}


\subsubsection{Characterization of the $\Gamma$-limit}\label{c-gl3d2d}

\noindent  In view of Lemma \ref{intrep3d2d}, we only need to prove
that ${\overline W}_{\rm hom}(\overline \xi)=
W_{\{\e_{n}\}}(\overline \xi)$ for all $\overline \xi \in \Rb^{3
\times 2}$, and thus it suffices to work with affine functions
instead of with general Sobolev functions.

We state, without proof, an equivalent result to Proposition
\ref{prop} for the  dimension reduction  case.

\begin{proposition}\label{prop1}
Given $M>0$, $\eta>0$,  and $\varphi: [0,+\infty) \rightarrow
[0,+\infty]$  a continuous and increasing function satisfying
$\varphi (t)/t \to +\infty$ as $t \to +\infty$, there exists
$\e_{0}\equiv \e_{0}(M,\eta)>0$ such that for every  $0 <\e <
\e_{0}$, every $a \in \Rb^2$ and every $u\in W^{1,p}((a+Q') \times
I;\Rb^3)$ with
\begin{equation}\label{prop-of-u}
\displaystyle\int_{(a+Q') \times I}\varphi (|\nabla u|^{p})\,
dx\leq M,
\end{equation}
\noindent there exists $v \in W^{1,p}_{0}((a+Q') \times I;\Rb^3)$
with $\|v\|_{L^{p}((a+Q') \times I;\Rb^3)}\leq \eta$ satisfying
$$\displaystyle\int_{(a+Q') \times I}
W\Big(x_3,x_\a,\frac{x_3}{\e},\frac{x_\a}{\e};\nabla u \Big)\,dx
\geq \displaystyle\int_{(a+Q') \times I} W_{\rm
hom}(x_3,x_\a;\nabla u + \nabla v)\, dx -\eta.$$
\end{proposition}

\begin{lemma}\label{lem3d2d<}
For all $\overline \xi\in \Rb^{3\times 2}$, ${\overline W}_{\rm hom}(\overline \xi)
\leq  W_{\{\e_{n}\}}(\overline \xi)$.
\end{lemma}

\noindent\textit{Proof. }From Lemma \ref{intrep3d2d}, Theorem 1.1 in Bocea and Fonseca and Lemma 2.6 in Braides, Fonseca and Francfort, we may find a sequence  $\{w_n\} \subset W^{1,p}(Q' \times I;\Rb^3)$ such that $\big\{\big|\big(\nabla_\a w_n \big|\frac{1}{\e_n} \nabla_3 w_n \big) \big|^p\big\}$ is equi-integrable, $w_n=0$ on $\partial Q' \times I$, $w_n \to 0$ in $L^p(Q' \times I;\Rb^3)$ and
$$2 W_{\{\e_{n}\}}(\overline \xi)=\lim_{n\to +\infty}\int_{Q' \times I}
W\left(x_3,\frac{x}{\e_{n}}, \frac{x_\a}{{\e_{n}}^{2}};  \overline \xi+ \nabla_\a
w_{n}(x)\Big|\frac{1}{\e_n} \nabla_3 w_n(x) \right)\, dx.$$
Thus, from De La Vall\'ee Poussin criterion (see e.g. Proposition 1.27 in Ambrosio, Fusco and
Pallara \cite{AFP}) there exists an increasing continuous function
$\varphi: [0,+\infty)\rightarrow [0,+\infty]$ satisfying $\varphi(t) /
t \to  +\infty$ as $t \to  +\infty$ and such that
$$\sup_{n \in \Nb}\int_{Q' \times I} \varphi\left (\left|\left(\nabla_\a w_n \Big|\frac{1}{\e_n} \nabla_3 w_n \right) \right|^p \right)\,dx\leq 1.$$
\noindent Changing variables yields
$$W_{\{\e_n\}}(\overline \xi)= \lim_{n \to +\infty}
\frac{1}{2 T_n^2}\int_{(0,T_n)^2 \times I} W\left(
x_3,x_\a,\frac{x_3}{\e_n},\frac{x_\a}{\e_n^2};\overline \xi+
\nabla_\a z_n(x)|\nabla_3 z_n(x) \right)\, dx$$ and
$$\sup_{n \in \Nb} \frac{1}{T_{n}^{2}}\int_{(0,T_{n})^{2}\times I} \varphi
(|\nabla z_{n}|^{p})\,dx \leq 1,$$ \noindent where we set
$T_n:=1/\e_n$ and $z_n(x):=T_n w_n (x_\a/T_n,x_3)$. Note that $z_n
\in W^{1,p}((0,T_n)^2 \times I;\Rb^3)$ and $z_n = 0$ on $\partial
(0,T_n)^2 \times I$. For all $n \in \Nb$, define
$I_n:=\left\{1,\cdots,\lb T_n \rb ^2\right\}$ and for any $i \in
I_n$, take $a_i^n \in \Zb^2$ such that
$$\bigcup_{i \in I_n} (a_i^n+Q') \subset (0,T_n)^2.$$
Moreover, for all $M>0$, let
$$I_n^M:=\left\{ i \in I_n : \int_{(a_i^n+Q') \times I} \varphi
(|\nabla z_{n}|^{p})\,dx \leq M \right\}.$$ Applying Proposition
\ref{prop1}, we get for any $\eta>0$ and any $i \in I_{n}^{M}$ the
existence of $v_i^{n,M,\eta} \in W_{0}^{1,p}((a_{i}^n+Q') \times
I;\Rb^{3})$ with $\|{v}_i^{n,M,\eta}\|_{L^{p}((a_{i}^n+Q') \times
I;\Rb^{3})}\leq \eta$ and
\begin{eqnarray*}
&&\int_{(a_{i}^n+Q' \times I}W\left(x_3,x_\a, \frac{x_3}{\e_n},\frac{x_\a}{\e_n}; \overline \xi+ \nabla_\a {z}_n| \nabla_3 z_n\right)\, dx\\
&&\hspace{1.0cm}\geq  \frac{1}{T_n^2} \int_{(a_i^n +Q') \times I}W_{\rm hom}\left( x_3,x_\a; \overline \xi +\nabla_\a{z}_n+ \nabla_\a v_i^{n,M,\eta}|\nabla_3 {z}_n+ \nabla_3 v_i^{n,M,\eta} \right) dx - \eta.
\end{eqnarray*}
\noindent Hence,
\begin{eqnarray}\label{3d2d5}
 W_{\{\e_n\}}(\overline \xi) \geq \limsup_{M, \eta, n} \frac{1}{2 T_n^2}
\sum_{i \in I_n^M}\int_{(a_i^n +Q') \times I}W_{\rm hom}\left( x_3,x_\a; (\overline \xi|0) + \nabla \phi^{n,M,\eta} \right) dx
\end{eqnarray}
\noindent where $\phi^{n,M,\eta} \in W^{1,p}((0;T_n)^2 \times
I;\Rb^3)$ is defined by
$$\phi^{n,M,\eta}(x):=\left\{\begin{array}{ll}
z_{n}(x)+v_i^{n,M,\eta}(x) & \text{if } x \in (a_{i}^n+Q') \times I \text{ and } i \in I_{n}^{M},\\
&\\
z_n(x) & \text{otherwise}
\end{array}\right.$$
and satisfies $\phi^{n,M,\eta}=0$ on $\partial (0,T_n)^2 \times I$.
\noindent In view of the definition of $\phi^{n,M,\eta}$, the
$p$-growth condition (\ref{Whom4}) and the equi-integrability of
$\big\{\big|\big(\nabla_\a w_n \big|\frac{1}{\e_n} \nabla_3 w_n
\big) \big|^p\big\}$, we get arguing exactly as in Lemma \ref{lem<},
$$W_{\{\e_n\}}(\overline \xi) \geq \limsup_{M, \eta, n}
\frac{1}{2 T_n^2}\int_{(0,T_n)^2 \times I} W_{\rm hom}( x_3,x_\a;\overline \xi +\nabla_\a {\phi}^{n,M,\eta}| \nabla_3 {\phi}^{n,M,\eta} )\, dx\\
\geq \overline W_{\rm hom}(\overline \xi).$$
$\hfill\blacksquare$\\

Let us now prove the converse inequality.
\begin{lemma}\label{lem3d2d>}
For all $\xi\in \Rb^{3\times 2}$, ${\overline W}_{\rm hom}(\overline \xi)
\geq  W_{\{\e_{n}\}}(\overline \xi)$.
\end{lemma}
\noindent\textit{Proof. } In view of (\ref{Whombarra2}), for
$\delta>0$ fixed take  $T\equiv T_{\delta} \in \Nb$, with
$T_{\delta} \to  +\infty$ as $\delta  \to  0$, and $\phi\equiv
{\phi_{\delta}}\ \in W^{1,p}((0,T)^{2} \times I;\Rb^{3})$ be such that
$\phi=0$ on $\partial (0,T)^2 \times I$ and
\begin{equation}\label{delta3d2d}
{\overline W}_{\rm hom}(\overline \xi)+\delta \geq \frac{1}{{2T}^{2}}
\int_{(0,T)^{2} \times I}W_{\rm hom}(x_3,x_\a; \overline \xi+\nabla_\a {\phi}(x)|\nabla_\a \phi(x))\, dx.
\end{equation}
\noindent From Theorem 1.1 in Ba\'{\i}a and Fonseca \cite{B&Fo}
(with $f(y,z;\xi)=W(y_3,y_\a,z_3,z_\a;\xi)$), Proposition 11.7 in
Braides and Defranceschi \cite{BD} and Decomposition Lemma (see
Fonseca, M\"uller and Pedregal \cite{FMP} or Fonseca and Leoni
\cite{FL}) there exists $\{\phi_{n}\}\subset W_{0}^{1,p}((0,T)^{2}
\times I;\Rb^{3})$ such that $\{|\nabla \phi_n|^p\}$ is
equi-integrable, $\phi_{n}\rightarrow \phi$ in $L^{p}((0,T)^{2}
\times I;\Rb^{3})$ and
\begin{eqnarray}\label{b-n3d2d}
&&\int_{(0,T)^{2} \times I}W_{\rm hom}(x_3,x_\a; \overline \xi+\nabla_\a {\phi}(x)|\nabla_3\phi(x) )\,
dx\nonumber\\
&&\hspace{2.0cm}=\lim_{n \to +\infty}\int_{(0,T)^{2}\times I}
W\left(x_3,x_\a,\frac{x_3}{\e_{n}},\frac{x_\a}{\e_{n}};\overline
\xi+\nabla_\a \phi_{n}(x)|\nabla_3 \phi_n(x) \right)
dx.\end{eqnarray} \noindent Fix $n\in \Nb$ such that  $\e_{n} \ll
1$. For all $i \in \Zb^2$ let $a_{i}^n \in \e_n \Zb^2 \cap (i(T+1) +
[0,\e_n)^2)$ (uniquely defined). Set
$$\tilde \phi_n(x):=\left\{
\begin{array}{cl}
\phi_n(x_\a-a_{i}^n,x_3) & \text{if }  x \in (a_i^n + (0,T)^2) \times I \text{ and } i\in \Zb^{2},\\
&\\
0 & \text{otherwise},\hspace{0.0cm}
\end{array}
\right.$$ then $\tilde \phi_n \in W^{1,p}(\Rb^2 \times I;\Rb^3)$.
Let $I_n:=\{ i \in \Zb^2 : (0,T/\e_n)^2 \cap (a_i^n + (0,T)^2) \ne
\emptyset \}$. If $\psi_n(x):=\e_n \tilde \phi_n(x_\a/\e_n,x_3)$
then $\psi_n \to 0$ in $L^p((0,T)^2 \times I;\Rb^3)$, as $n  \to
+\infty$. Consequently, the $p$-growth condition $(A_4)$ implies
that
\begin{eqnarray}\label{1503}
&&\hspace{-1.0cm}\W_{\{\e_n\}}(\overline \xi \cdot\,;(0,T)^2)\nonumber\\
&& \leq \liminf_{n \to
+\infty} \int_{(0,T)^2 \times I}
W \left(x_3, \frac{x}{\e_n},\frac{x_\a}{\e_n^2};\overline\xi +\nabla_\a \psi_n(x)\Big|\frac{1}{\e_n} \nabla_3 \psi_n(x) \right)\, dx\nonumber\\
&& \leq  \liminf_{n \to +\infty} \e_n^2 \left\{ \sum_{i \in
I_n}\int_{(a_i^n + (0,T)^2) \times I} W \left(
x_3,x_\a,\frac{x_3}{\e_n},\frac{x_\a}{\e_n};
\overline \xi +\nabla_\a \tilde \phi_n(x)| \nabla_3 \tilde \phi_n(x) \right)\, dx \right.\nonumber\\
&& \hspace{2.0cm} \left. + c \mathcal L^2 \left(\Big(0,\frac{T}{\e_n}\Big)^2 \setminus \bigcup_{i \in I_n}(a_i^n + (0,T)^2) \right) \right\}\nonumber\\
&& = \liminf_{n \to +\infty} \e_n^2 \sum_{i \in I_n} \int_{(0,T)^2
\times I} W \left(x_3,x_\a+a_i^n-i(T+1)
,\frac{x_3}{\e_n},\frac{x_\a}{\e_n};\overline \xi
+\nabla_\a \phi_n(x)|\nabla_3 \phi_n(x) \right)\, dx \nonumber\\
&& \hspace{2.0cm} + c T^{2} \left(1- \left(\frac{T}{T+1} \right)^{2}
\right)\,
\end{eqnarray}
\noindent where we have used $(A_3)$, the fact that $T \in \Nb$ and
$a_i^n/\e_n \in \Zb^2$. We now use the same uniform continuity
argument than in the proof of Lemmas \ref{independency3d2d} and
\ref{lem>}, we get
$$ W_{\{\e_n\}}(\overline \xi) \leq \overline W_{\rm hom} (\overline \xi) +\delta +  c \left(1-
\left(\frac{T}{T+1}\right)^{2}\right).$$ The result follows by
letting $\delta$ tend to zero.
$\hfill\blacksquare$\\

\noindent{\it Proof of Theorem \ref{jf-m4}.} From Lemma
\ref{lem3d2d<} and Lemma \ref{lem3d2d>}, we conclude that
${\overline W}_{\rm hom}(\overline \xi) = W_{\{\e_{n}\}}(\overline
\xi)$ for all $\overline \xi\in \Rb^{3\times 2}$. As a consequence,
$\W_{\{\e_n\}}(u;A)= \W_{\rm hom}(u;A)$ for all $A \in \mathcal A_0$
and all $u \in W^{1,p}(A;\Rb^3)$. Since the $\G$-limit does not
depend upon the extracted subsequence, Proposition 8.3 in Dal Maso
\cite{DM} implies that the whole sequence $\W_\e(\cdot\,;A)$
$\G(L^p(A \times I))$-converges to $\W_{\rm hom}(\cdot\,;A)$.
\hfill$\blacksquare$


\subsection{The general case}\label{3hetero}

\noindent Our aim here is to study the case where the function $W$
depends also on the in-plane variable.


\subsubsection{Existence and integral representation of the $\Gamma$-limit}

\noindent As in Section \ref{2hetero}, to prove Theorem
\ref{jf-m3} it is convenient  to localize the functionals
$\W_{\e}$  in (\ref{Wepsilon}) on the class of all bounded open
subsets of $\o$, denoted by $\mathcal A(\o)$. For each $\e>0$ we
consider the family of functionals $\W_{\e}: L^{p}(\O;\Rb^{3})
\times {\cal A}(\o) \to [0,+\infty]$ defined by
\begin{equation}\label{We}
\W_{\e}(u;A):=\left\{\begin{array}{ll} \ds \int_{A \times I}
W\left(x,\frac{x}{\e},\frac{x_\a}{\e^{2}};\nabla_\a
u(x)\Big|\frac{1}{\e} \nabla_3 u(x) \right)\, dx & \text{if } u \in
W^{1,p}(A \times I;\Rb^{3}),\\
+\infty & \text{otherwise}.
\end{array}\right.
\end{equation}
\noindent Given $\{\e_j\} \searrow 0^+$  and $A \in \mathcal
A(\o)$ we define the $\Gamma$-lower limit  of $\{  \W_{\e_j}
(\cdot\,;A) \}_{j \in \Nb}$ with respect to the  $L^{p}(A \times
I;\Rb^{3})$-topology  by
$$\W_{\{\e_j\}}(u;A) := \inf_{\{u_j\}}\left\{\liminf_{j \to
+\infty} \W_{\e_j}(u_j;A) : \quad u_j \to u \text{ in }L^p(A \times
I;\Rb^3)\right\}$$ \noindent  for all $u \in L^{p}(\O;\Rb^3)$. Our
main objective is to show that
\begin{equation}\label{1856}
\W_{\{\e_j\}}=\W_{\rm hom}
\end{equation}
\noindent where $\W_{\rm hom} : L^p(\O;\mathbb R^3) \times \mathcal
A(\o) \to [0,+\infty]$ is given by
$$\W_{\rm hom}(u;A) =
\left\{\begin{array}{ll} \ds 2\int_{A} \overline W_{\rm
hom}(x_\a;\nabla_\a u(x_\a))\, dx_\a & \text{if } u  \in
W^{1,p}(A;\Rb^{3}),\\ + \infty &\text{otherwise}.
\end{array}\right.$$
\noindent The conclusion of Theorem \ref{jf-m3} would follow taking
$A=\o$. By  hypotheses $(A_{4})$ it follows that
$\W_{\{\e_j\}}(u;A)=+ \infty$ for each $A \in \mathcal A(\o)$
whenever $u \in L^p(\O;\Rb^3) \setminus W^{1,p}(A;\Rb^3)$.  As a
consequence of Theorem 2.5 in Braides, Fonseca and Francfort
\cite{BFF}, given $\{\e_{j}\}\searrow 0^{+}$ there exists  a
subsequence $\{\e_{j_n}\} \equiv \{\e_{n}\}$ of $\{\e_{j}\}$ for
which the functional $\W_{\{\e_{n}\}}(\cdot \,;A)$ is the $\G(L^p(A
\times I))$-limit of $\{\W_{\e_n} (\cdot\,;A) \}_{n \in \Nb}$ for
each $A \in \mathcal A(\o)$.  Moreover given  $u \in
W^{1,p}(A;\Rb^3)$
$$\W_{\{\e_n\}} (u;A)=2\int_A W_{\{\e_n\}}(x_\a;\nabla_\a(x_\a))\, dx_\a,$$
\noindent for some  Carath\'eodory function $W_{\{\e_n\}} : \o
\times \Rb^{3 \times 2} \to \Rb$. Accordingly, to prove equality
(\ref{1856}) it suffices to show that $W_{\{\e_n\}}(x_\a;\overline
\xi)=\overline W_{\rm hom}(x_\a;\overline \xi)$ for a.e. $x_\a \in
\o$ and all $\overline \xi \in \Rb^{3 \times 2}$, which allow us to
work with affine functions instead of with general  Sobolev
functions.

The following proposition, that is of use in the sequel, allow us to
extend  continuously Carath\'eodory integrands. It relies on
Scorza-Dragoni's  Theorem (see Ekeland and Temam \cite{Ek&Te}) and
on Tietze's  Extension Theorem (see Theorem 3.1 in DiBenedetto
\cite{DiB}).

\begin{proposition}\label{tietze}
Let $W: \O \times \Rb^3 \times \Rb^2 \times \Rb^{3 \times 3} \to
\Rb$ satisfying $(A_1)$-$(A_4)$. Then for any $m\in \Nb$, there
exists a compact set $C_{m} \subset \O$ and a continuous function
$W^m: \O \times \Rb^3 \times \Rb^2 \times \Rb^{3 \times 3} \to
\Rb$ such that
$W^m(x,\cdot,\cdot\,;\cdot)=W(x,\cdot,\cdot\,;\cdot)$ for all $x
\in C_{m}$ and
\begin{equation}\label{TEE}
\mathcal L^3(\O \setminus C_{m}) <\frac{1}{m}.
\end{equation}
Moreover,
\begin{itemize}
\item[-]$y_\a \mapsto W^m(x,y_\a,y_3,z_\a;\xi)$ is $Q'$-periodic
for all $(z_\a,y_3,\xi)\in \Rb^{3} \times \Rb^{3\times 3}$ and
a.e. $x \in \Omega$, \item[-] $(z_\a,y_3) \mapsto
W^m(x,y_\a,y_3,z_\a ;\xi)$ is $Q$-periodic for all $(y_\a,\xi)\in
\Rb^{2}\times\Rb^{3\times 3}$ and a.e. $x \in \Omega$;\
\end{itemize}
and for some $\beta>0$, we have
\begin{equation}\label{PG}
 -\beta\leq W^m(x,y,z_\a;\xi)\leq
\beta(1+|\xi|^p) \quad \text{for  all }(y,z_\a,\xi) \in \Rb^{3}
\times \Rb^2 \times \Rb^{3\times 3} \text{ and a.e }x \in \Omega.
\end{equation}
\end{proposition}

\noindent {\it Proof. }By Scorza-Dragoni's Theorem (see Ekeland and
Temam \cite{Ek&Te}) for any $m\in \Nb$  there exists a compact set
$C_{m} \subset \O$ with $\mathcal L^3(\O \setminus C_{m}) <1/m$ such
that $W$ is continuous on $C_{m} \times \Rb^3 \times \Rb^2 \times
\Rb^{3 \times 3}$. Since $C_{m} \times \Rb^3 \times \Rb^2 \times
\Rb^{3 \times 3}$ is a closed set, according to Tietze's Extension
Theorem (see DiBenedetto \cite{DiB}) one can extend $W$ into a
continuous function $W^m$ outside $C_{m} \times \Rb^3 \times \Rb^2
\times \Rb^{3 \times 3}$. By the construction of $W^m$ it can be
seen that it  satisfies the same periodicity and growth condition
than $W$ and that it is bounded from below by $-\beta$.
$\hfill\blacksquare$\\

We  remark  that the above  result improve Lemma 4.1 in Babadjian
and Ba\'{\i}a \cite{BB1} in which we only obtained a separately
continuous function.


\subsubsection{Characterization of the $\Gamma$-limit}

\noindent For each  $T>0$ consider ${\cal S}_{T}$ a countable set
of functions in $\mathcal C^\infty([0,T]^2 \times [-1,1];\Rb^{3})$
that is dense in

$$\{ \varphi \in W^{1,p}((0,T)^{2}\times I;\Rb^{3}):\,\, \varphi=0 \,\text{
on } \,  \partial{(0,T)^{2}} \times I\}.$$

\noindent Let $L$  be the set of Lebesgue points $x_{\alpha}^{0}$
for all functions $W_{\{\e_n\}}(\cdot;\overline{\xi})$, $\overline
W_{\rm  hom}(\cdot;\overline{\xi})$ and
$$x_{\alpha}\mapsto \int_{(0,T)^{2} \times I}W_{\rm
hom}(x_{\alpha},y_3,
y_{\alpha};\overline{\xi}+\nabla_{\alpha}\varphi(y)
|\nabla_{3}\varphi(y))\, dy,$$ \noindent with $T \in \Nb$,
$\varphi\in {\cal S}_{T}$ and $\overline{\xi}\in \Qb^{3\times 2}$,
and for which $\overline W_{\rm hom}(x_{\alpha}^{0};\,\cdot\,)$ is
well defined. Note that $\mathcal L^2(\o \setminus L)=0$.

We start by proving the following inequality.
\begin{lemma}\label{3d2d-general-ineq<}
For all $x_\a^0 \in L$ and all $\overline \xi \in \Qb^{3 \times
2}$, we have $W_{\{\e_n\}}(x^{0}_{\alpha};\overline{\xi})\geq
\overline W_{\rm hom}(x^{0}_{\alpha};\overline{\xi}).$
\end{lemma}

\noindent {\it Proof. }  Let $\delta>0$ small enough so that
$Q'(x_\a^0,\delta)\in \mathcal A (\o)$. By Theorem 1.1 in Bocea and
Fonseca \cite{Bo&Fo} we can find a sequence $\{u_n\} \subset
W^{1,p}(Q'(x_\a^0,\delta) \times I;\Rb^3)$ with $u_n \to 0$ in
$L^p(Q'(x_\a^0,\delta) \times I;\Rb^3)$, such that the sequence of
scaled gradients  $\big\{  \big( \nabla_\a u_n | \frac{1}{\e_n}
\nabla_3 u_n \big)  \big\}$ is $p$-equi-integrable and
\begin{eqnarray*}\W_{\{\e_n\}}(\overline \xi \cdot\,;Q'(x_\a^0,\delta)) &
=& 2\int_{Q'(x_{\alpha}^0,\delta)}
W_{\{\e_n\}}(x_{\alpha};\overline \xi)\, dx_{\alpha}\\
& = &\lim_{n \to +\infty} \int_{Q'(x_\a^0,\delta) \times I}W
\left(x,\frac{x}{\e_n},\frac{x_\a}{\e_n^2};\overline \xi + \nabla_\a
u_n(x) \Big| \frac{1}{\e_n} \nabla_3 u_n(x)
\right)dx.\end{eqnarray*} \noindent Given $m\in \Nb$ let $C_{m}$ and
$W^m$ be given by Proposition \ref{tietze}.  Then since $W\geq 0$
and $W=W^m$ on $C_{m}\times \Rb^{3}\times \Rb^{2}\times \Rb^{3\times
3}$ we get
$$\W_{\{\e_n\}}(\overline \xi \cdot\,;Q'(x_\a^0,\delta)) \geq \limsup_{m, n}
\int_{[Q'(x_\a^0,\delta) \times I]\cap C_{m}}W^m
\left(x,\frac{x}{\e_n},\frac{x_\a}{\e_n^2};\overline \xi +
\nabla_\a u_n(x) \Big|
 \frac{1}{\e_n} \nabla_3 u_n(x) \right)dx.$$
By the  $p$-growth condition (\ref{PG}), the equi-integrability of
$\big\{ \big| \big( \nabla_\a u_n | \frac{1}{\e_n} \nabla_3 u_n
\big) \big|^p \big\}$ and relation (\ref{TEE}), we obtain
\begin{eqnarray*}
&&\int_{[Q'(x_\a^0,\delta) \times I] \setminus C_{m}} W^m \left(x,\frac{x}{\e_n},\frac{x_\a}{\e_n^2};\overline \xi + \nabla_\a u_n(x) \Big| \frac{1}{\e_n} \nabla_3 u_n(x) \right)dx\\
&&\hspace{2.0cm} \leq \beta  \int_{[Q'(x_\a^0,\delta) \times I]
\setminus C_{m}} \left( 1 + \left| \left( \overline \xi +
\nabla_\a u_n(x) \Big| \frac{1}{\e_n} \nabla_3 u_n(x)
\right)\right|^p \right)dx \xrightarrow[m \to +\infty]{} 0,
\end{eqnarray*}
uniformly with respect to $n \in \Nb$. Then, we get that
$$\W_{\{\e_n\}}(\overline \xi \cdot\,;Q'(x_\a^0,\delta)) \geq \limsup_{m, n}
\int_{Q'(x_\a^0,\delta) \times I}W^m
\left(x,\frac{x}{\e_n},\frac{x_\a}{\e_n^2};\overline \xi + \nabla_\a
u_n(x) \Big| \frac{1}{\e_n} \nabla_3 u_n(x) \right)dx.$$ \noindent
For any $h\in \Nb$, we split $Q'(x_\a^0,\delta)$ into $h^{2}$
disjoint cubes $Q'_{i,h}$ of side length $\delta/h$ so that
$$Q'(x_\a^0,\delta) =   \bigcup_{i=1}^{h^{2}}Q'_{i,h} $$
\noindent and
\begin{eqnarray*}
\W_{\{\e_n\}}(\overline \xi \cdot\,;Q'(x_\a^0,\delta)) &\geq &
\limsup_{m,h,n} \sum_{i=1}^{h^{2}} \int_{ Q'_{i,h} \times I}W^m
\left(x,\frac{x}{\e_n},\frac{x_\a}{\e_n^2};\overline \xi + \nabla_\a
u_n(x) \Big| \frac{1}{\e_n} \nabla_3 u_n(x) \right)dx\\
&\geq & \limsup_{m,\lambda,h,n} \sum_{i=1}^{h^{2}} \int_{ [Q'_{i,h}
\times I]\cap R^{\lambda}_{n}}W^m
\left(x,\frac{x}{\e_n},\frac{x_\a}{\e_n^2};\overline \xi + \nabla_\a
u_n(x) \Big| \frac{1}{\e_n} \nabla_3 u_n(x) \right)dx
\end{eqnarray*}
\noindent where given $\lambda >0$ we define
$$R^{\lambda}_{n}:=\left\{x\in Q'(x_\a^0,\delta)\times I:\quad  \left | \left (\overline \xi +
\nabla_\a u_n(x) \Big| \frac{1}{\e_n} \nabla_3 u_n(x) \right )
\right | \leq \lambda \right\}.$$ \noindent Since $W^m$ is
continuous and separately periodic it is in particular uniformly
continuous on $\overline \O \times \Rb^3 \times \Rb^2 \times
\overline{B}(0,\lambda)$. With similar arguments to that used in the
proof of Lemma 3.5 in Babadjian and Ba\'{\i}a \cite{BB1} (with $W^m$
in place  of  $W^{m,\lambda}$), we obtain
\begin{eqnarray}\label{507}
&&\W_{\{\e_n\}}(\overline \xi \cdot\,;Q'(x_\a^0,\delta))\nonumber\\
&& \geq \limsup_{h, n}  \frac{h^2}{\delta^2} \sum_{i=1}^{h^2} \int_{Q'_{i,h}} \int_{Q'_{i,h} \times I}W \left(x'_\a,x_3,\frac{x}{\e_n},\frac{x_\a}{\e_n^2};\overline \xi + \nabla_\a u_n(x) \Big| \frac{1}{\e_n} \nabla_3 u_n(x) \right)dx\, dx'_\a \nonumber \\
&& \geq \hspace{-0.1cm}\limsup_{h \to +\infty}
\frac{h^2}{\delta^2}  \sum_{i=1}^{h^2} \hspace{-0.1cm}
\int_{Q'_{i,h}} \hspace{-0.1cm} \liminf_{n \to +\infty}
\hspace{-0.1cm} \int_{Q'_{i,h} \times I} \hspace{-0.1cm} W
\left(x'_\a,x_3,\frac{x}{\e_n},\frac{x_\a}{\e_n^2};\overline \xi +
\nabla_\a u_n(x) \Big| \frac{1}{\e_n} \nabla_3 u_n(x) \right)dx\,
dx'_\a,
\end{eqnarray}
\noindent where we have used Fatou's Lemma.  We now fix $x'_\a \in
Q'_{i,h}$ such that $\overline W_{\rm hom}(x'_\a;\overline \xi)$ is
well defined,  then by Theorem \ref{jf-m4} we get that
\begin{equation}\label{508}
\liminf_{n \to +\infty} \int_{Q'_{i,h} \times I}W \left(x'_\a,x_3,\frac{x}{\e_n},
\frac{x_\a}{\e_n^2};\overline \xi + \nabla_\a u_n(x) \Big| \frac{1}{\e_n} \nabla_3 u_n(x) \right)dx
 \geq 2 \frac{\delta^2}{h^2} \overline W_{\rm hom}(x'_\a;\overline \xi).
\end{equation}
\noindent Gathering (\ref{507}) and (\ref{508}), it turns out that
$$\int_{Q'(x_\a^0,\delta)} W_{\{\e_n\}}(x_\a;\overline \xi)\, dx_\a \geq \int_{Q'(x_\a^0,\delta)}
\overline W_{\rm hom}(x'_\a;\overline \xi)\, dx'_\a.$$ \noindent As
a consequence the claim follows  by the choice of $x^{0}_{\alpha}$,
after  dividing the previous inequality by $\delta^{2}$ and letting
$\delta\go 0$.
$\hfill\blacksquare$\\

We  now  prove the converse inequality.

\begin{lemma}\label{3d2d-general-ineq>}
For all $\xi \in \Qb^{3\times 2}$ and all $x^{0}_{\alpha} \in L$,
$ W_{\{\e_n\}}(x^{0}_\a;\overline \xi)\leq \overline  W_{\rm
hom}(x_\a;\overline \xi).$
\end{lemma}

\noindent{\it Proof.}  For every  $m\in \Nb$,  consider  the set
$C_{m}$ and the function $W^{m}$   given by Proposition
\ref{tietze}, and define $\overline{(W^{m})}_{\rm hom}$ and
${(W^{m})}_{\rm hom}$ as (\ref{Whombarra}) and (\ref{Whom}), with
$W^{m}$ in place of $W$. For fixed $\eta>0$ and any $m\in \Nb$ let
$K_{\eta}^{m}$ be a compact subset of $\o$ given by Scorza-Dragoni's
Theorem (see Ekeland and Temam \cite{Ek&Te}) with ${\mathcal
L}^{2}(\o\setminus K_{\eta}^{m})\leq \eta$ and such that
$\overline{(W^{m})}_{\rm hom}: K_{\eta}^{m}\times
\Rb^{3\times 2}\rightarrow \Rb$ is continuous .\\

{\it Step 1.} We claim that
\begin{eqnarray}\label{642}
2\liminf_{m \go +\infty}\int_{Q^{\prime}(x_{\alpha}^{0},\delta)}
\overline{(W^m)}_{\rm hom} (x_{\alpha};\overline \xi)\, dx_{\alpha}
&\geq & 2\int_{Q^{\prime}(x_{\alpha}^{0},\delta)}
W_{\{\e_{n}\}}\left( x_{\alpha},\overline{\xi} \right)\,
dx_{\alpha}.
\end{eqnarray}
\noindent To show this inequality we follow a similar argument  to
that of Lemma \ref{ineq>}. As before, we first  decompose
$Q^{\prime}(x_{\alpha}^{0},\delta)$ into $h^2$ small disjoint cubes
$Q^{\prime}_{i,h}$ and we set
$$I^{m}_{h,\eta}:=\left\{ i\in \{1,\cdots,h^2 \}:
\quad K_{\eta}^{m} \cap Q^{\prime}_{i,h} \neq \emptyset \right\}.$$
\noindent For $i\in I^{m}_{h,\eta}$ choose $x_i^{h,\eta,m} \in
K^{m}_\eta \cap Q^{\prime}_{i,h}$. By Theorem \ref{jf-m4} together
with Lemma 2.6 in Braides, Fonseca and Francfort \cite{BFF} there
exists   a sequence $\{u_i^{n,h,\eta,m}\} \subset
W^{1,p}(Q^{\prime}_{i,h}\times I;\Rb^3)$ with $u_i^{n,h,\eta,m}=0$
on $\partial  Q^{\prime}_{i,h} \times I $,  $u_i^{n,h,\eta,m}
\xrightarrow [n\go +\infty] {} 0$ in $L^p(Q^{\prime}_{i,h}\times
I;\Rb^3)$, and such that
\begin{eqnarray*}
&&  2\int_{Q^{\prime}_{i,h}} \overline{(W^{m})}_{\rm
hom}(x_i^{h,\eta,m};\overline \xi)\, dx_{\alpha}\\
&& \hspace{2cm} =  \lim_{n \to +\infty} \int_{Q^{\prime}_{i,h}\times
I}W^{m}\left(x_i^{h,\eta,m},x_{3},\frac{x}{\e_n},\frac{x_{\alpha}}{\e_n^2};\overline\xi
+ \nabla_{\alpha} u_i^{n,h,\eta,m}\Big | \frac{1}{\e_n}\nabla_{3}
u_i^{n,h,\eta,m} \right)\, dx.
\end{eqnarray*}
\noindent  Setting
$$u_n^{\eta,m}(x):= \left\{
\begin{array}{ll}
u_i^{n,h,\eta,m}(x) & \text{if }x_{\alpha} \in Q^{\prime}_{i,h}\text{ and } i\in I^{m}_{h,\eta},\\
&\\
0 & \text{otherwise},
\end{array}
\right.$$ \noindent it follows that  $\{u_n^{\eta,m}\} \subset
W^{1,p}(Q^{\prime}(x_{\alpha}^{0},\delta)\times I;\Rb^3)$ and
$u_n^{\eta,m} \xrightarrow [n\go +\infty] {} 0$ in
$L^p(Q^{\prime}(x_{\alpha}^{0},\delta)\times I;\Rb^3)$.  Thus
\begin{eqnarray*}
&&\hspace{-1cm} 2\liminf_{\eta, h}\sum_{i \in I^{m}_{h,\eta}}
\int_{Q^{\prime}_{i,h}}
\overline{ (W^{m})}_{\rm hom}(x_i^{h,\eta,m};\overline \xi)\, dx_{\alpha}\nonumber\\
&& \hspace{2cm} \geq \liminf_{ \eta, h, n} \sum_{i \in
I^{m}_{h,\eta}} \int_{Q^{\prime}_{i,h}\times I}W^{m}\left(
x_i^{h,\eta,m},x_{3,}\frac{x}{\e_n},\frac{x_{\alpha}}{\e_n^2};\overline\xi
+ \nabla_{\alpha} u_n^{\eta,m} \Big | \frac{1}{\e_{n}}\nabla_{3}
u_n^{\eta,m} \right)\, dx.
\end{eqnarray*}
\noindent As in Lemma \ref{ineq>} we obtain
\begin{eqnarray}
&&\hspace{-0.5cm}2\int_{Q^{\prime}(x_{\alpha}^{0},\delta)}
\overline{( W^{m})}_{\rm hom}(x_{\alpha};\overline\xi)\, dx_{\alpha}\nonumber\\
&& \geq \liminf_{\lambda, \eta, h, n} \sum_{i \in I^{m}_{h,\eta}}
\int_{(Q^{\prime}_{i,h}\times I) \cap
R_{n,\eta,m}^\lambda}W^{m}\left(
x_i^{h,\eta,m},x_{3,}\frac{x}{\e_n},\frac{x_{\alpha}}{\e_n^2};\overline\xi
+ \nabla_{\alpha} u_n^{\eta,m}\Big | \frac{1}{\e_{n}}\nabla_{3}
u_n^{\eta,m} \right)\, dx,\nonumber
\end{eqnarray}
\noindent where $$R_{n,\eta,m}^\lambda:= \left\{ x \in
Q^{\prime}(x_{\alpha}^{0},\delta)\times I  :
\left|\left(\overline\xi + \nabla_{\alpha} u_n^{\eta,m}(x) \Big |
\frac{1}{\e_{n}}\nabla_{3} u_n^{\eta,m}(x)\right)\right|\leq \lambda
\right\}$$ \noindent with
\begin{equation}\label{1242}
{\mathcal
L}^{3}\left(\left[Q^{\prime}(x_{\alpha}^{0},\delta)\times I\right
] \setminus R_{n,\eta,m}^\lambda \right)\leq
\frac{C}{\lambda^{p}},
\end{equation}
\noindent for some constant $C>0$ independent of $n,\, \eta,\, m$
and $\lambda$. Taking into account that $W^{m}$ is continuous we get
that
$$2\int_{Q^{\prime}(x_{\alpha}^{0},\delta)} \overline{ (W^{m})}_{\rm hom}(x_{\alpha};\overline\xi)\, dx_{\alpha} \geq
\liminf_{\lambda, \eta, n} \int_{R_{n,\eta,m}^\lambda} W^{m}\left(
x,\frac{x}{\e_n},\frac{x_{\alpha}}{\e_n^2};\overline\xi +
 \nabla_{\alpha}
u_n^{\eta,m}\Big | \frac{1}{\e_{n}}\nabla_{3} u_n^{\eta,m}
\right)\, dx.$$ \noindent  By a diagonalization argument, given
$\lambda_m \nearrow +\infty$, and $\eta_m \searrow 0^{+}$ there
exists $n_m\nearrow +\infty$ such that
$$2\liminf_{m \go +\infty}\int_{Q^{\prime}(x_{\alpha}^{0},\delta)} \overline{ (W^m)}_{\rm hom}(x_{\alpha};
\overline\xi)\, dx_{\alpha} \geq \liminf_{m \go +\infty} \int_{R_m}
W^m\left(
x,\frac{x}{\e_{n_m}},\frac{x_{\alpha}}{\e_{n_m}^2};\overline\xi +
\nabla_{\alpha} v_m\Big | \frac{1}{\e_{n_m}}\nabla_{3} v_m \right)\,
dx,$$ \noindent where $v_m:= u_{n_m}^{\eta_m,m}\in
W^{1,p}(Q^{\prime}(x_{\alpha}^0,\delta)\times I;\Rb^3)$ with $v_m
\to 0$ in $L^p(Q^{\prime}(x_{\alpha}^0,\delta)\times I;\Rb^3)$, and
where $R_m:=R_{n_m,\eta_m,m}^{\lambda_m}.$ Using the Bocea and
Fonseca decomposition lemma for scaled gradients (Theorem 1.1 in
\cite{Bo&Fo}) we can assume, without loss of generality, that the
sequence $\big\{\big| \big( \nabla_{\alpha} v_m  |
\frac{1}{\e_{n_m}}\nabla_{3} v_m \big)\big|^{p}\big\}$ is
equi-integrable. Then, since $W^m=W$ on $C_m \times \Rb^{3}\times
\Rb^{2}\times \Rb^{3\times 3}$ it comes that
\begin{eqnarray*}
&& 2\liminf_{m\go +\infty}\int_{Q^{\prime}(x_{\alpha}^{0},\delta)}
\overline{ (W^m)}_{\rm hom}(x_{\alpha}; \overline\xi)\,
dx_{\alpha}\\
&& \hspace{2cm} \geq    \liminf_{m \to +\infty}
\int_{Q^{\prime}(x_{\alpha}^{0},\delta)\times I} W^m \left(
x,\frac{x}{\e_{n_m}},\frac{x_{\alpha}}{\e_{n_m}^2};\overline \xi +
\nabla_{\alpha} v_m \Big| \frac{1}{\e_{n_m}}\nabla_{3}
v_m \right)\, dx\\
&& \hspace{2cm} \geq    \liminf_{m \to +\infty}
\int_{[Q^{\prime}(x_{\alpha}^{0},\delta)\times I]\cap C_m} W\left(
x,\frac{x}{\e_{n_m}},\frac{x_{\alpha}}{\e_{n_m}^2};\overline \xi +
\nabla_{\alpha} v_m \Big| \frac{1}{\e_{n_m}}\nabla_{3}
v_m \right)\, dx\\
&& \hspace{2cm} = \liminf_{m \to +\infty}
\int_{Q^{\prime}(x_{\alpha}^{0},\delta)\times I}  W\left(
x,\frac{x}{\e_{n_m}},\frac{x_{\alpha}}{\e_{n_m}^2};\overline \xi +
\nabla_{\alpha} v_m\Big| \frac{1}{\e_{n_m}}\nabla_{3} v_m \right)\,
dx
\end{eqnarray*}
\noindent by the growth conditions on $W$, the
$p$-equi-integrability  of the above sequence of scaled gradients,
(\ref{TEE}) and  (\ref{1242}). As a result we get inequality
(\ref{642}).\\

{\it Step 2.}  Fixed $\rho>0$, let  $T\in \Nb$ and $\varphi \in
\mathcal S_{T}$ be such that
\begin{equation}\label{1747}
\overline W_{\rm hom}(x^{0}_\a;\overline \xi) + \rho \geq
\frac{1}{2 T^{2}}\int_{(0,T)^{2} \times I} W_{\rm
hom}(x^{0}_\a,y_3,y_\a;\overline \xi+\nabla_\a \varphi(y)|
\nabla_3 \varphi(y))\, dy.\end{equation} \noindent   Taking
$(T,\varphi)$ in the definition of $ \overline{(W^{m})}_{\rm hom}$
and  recalling Remark \ref{infwhom} (with $W^m$ in place of $W$)
it follows that
\begin{eqnarray}\label{1731}
&&\int_{Q^{\prime}(x_{\alpha}^{0},\delta)} \overline{(W^m)}_{\rm
hom}(x_\a;\overline \xi) \,
d{x_{\alpha}}\nonumber\\
&&\hspace{1.5cm} \leq \frac{1}{2
T^{2}}\int_{Q^{\prime}(x_{\alpha}^{0},\delta)} \int_{(0,T)^{2}
\times I} (W^m)_{\rm hom}(x_\a,y_3,y_\a;\overline \xi+\nabla_\a
\varphi(y)| \nabla_3 \varphi(y))\, dy\, dx_{\alpha}.\end{eqnarray}
\noindent Define $E_m:= \{(x_{\alpha},y_{\alpha},y_{3})\in
Q^{\prime}(x_{\alpha}^{0},\delta) \times (0,T)^{2}\times I: \quad
(x_{\alpha},y_{3})\in C_m  \}$. From (\ref{TEE}) it follows that
\begin{equation}\label{1732}
{\mathcal L}^{2}\otimes {\mathcal
L}^{3}\big([Q^{\prime}(x_{\alpha}^{0},\delta)\times (0,T)^{2}\times
I]\setminus E_m\big) \leq T^2/m.
\end{equation}
\noindent Since  $(W^m)_{\rm hom} = W_{\rm hom}$ on $C_m \times
\Rb^{2}\times \Rb^{3\times 3}$ it comes that
\begin{eqnarray}\label{1729}&&\int_{Q^{\prime}(x_{\alpha}^{0},\delta)}
\int_{(0,T)^{2} \times I} (W^m)_{\rm hom}(x_\a,y_3,y_\a;\overline
\xi+\nabla_\a \varphi(y)| \nabla_3
\varphi(y))\, dy\, dx_{\alpha}\nonumber \\
&& \hspace{2cm} = \int_{E_m} W_{\rm hom}(x_\a,y_3,y_\a;\overline
\xi+\nabla_\a \varphi(y)| \nabla_3
\varphi(y))\, dy\, dx_{\alpha}\nonumber\\
&  & \hspace{3cm} + \int_{[Q^{\prime}(x_{\alpha}^{0},\delta)\times
(0,T)^{2} \times I]\setminus E_m} (W^m)_{\rm
hom}(x_\a,y_3,y_\a;\overline \xi+\nabla_\a \varphi(y)| \nabla_3
\varphi(y))\, dy\, dx_{\alpha}\nonumber\\
&&\hspace{2cm}\leq  \int_{Q^{\prime}(x_{\alpha}^{0},\delta)\times
(0,T)^{2} \times I}W_{\rm hom}(x_\a,y_3,y_\a;\overline
\xi+\nabla_\a \varphi(y)| \nabla_3
\varphi(y))\, dy\, dx_{\alpha}\nonumber\\
& & \hspace{3cm} +  C \int_{[Q^{\prime}(x_{\alpha}^{0},\delta)\times
(0,T)^{2} \times I]\setminus E_m} (1+|\nabla \varphi (y)|^{p}) \,
dy\, dx_{\alpha}
\end{eqnarray}
\noindent  by property (\ref{Whom4}) with $W^m$ in place of $W$.
Passing to the limit as $m \go +\infty$, relations (\ref{1731}),
(\ref{1732}) and (\ref{1729}) yield to
\begin{eqnarray*}&&\limsup_{m \go
+\infty}\int_{Q^{\prime}(x_{\alpha}^{0},\delta)}
\overline{(W^m)}_{\rm hom}(x_\a;\overline
\xi) \, d{x_{\alpha}} \\
&& \hspace{2cm} \leq  \frac{1}{2
T^{2}}\int_{Q^{\prime}(x_{\alpha}^{0},\delta)} \int_{(0,T)^{2}
\times I}   W_{\rm hom}(x_\a,y_3,y_\a;\overline \xi+\nabla_\a
\varphi(y)| \nabla_3 \varphi(y))\, dy\, dx_{\alpha}.\end{eqnarray*}
\noindent Hence by (\ref{642}) we  obtain
\begin{eqnarray*}
\int_{Q^{\prime}(x_{\alpha}^{0},\delta)}  W_{\{\e_{n}\}}\left(
x_{\alpha},\overline{\xi} \right)\, dx_{\alpha} &\leq  &  \frac{1}{2
T^{2}}\int_{Q^{\prime}(x_{\alpha}^{0},\delta)} \int_{(0,T)^{2}
\times I} \hspace{-0.1cm} W_{\rm hom}(x_\a,y_3,y_\a;\overline
\xi+\nabla_\a \varphi(y)| \nabla_3 \varphi(y))\, dy\, dx_{\alpha}.
\end{eqnarray*}
\noindent As a consequence, by the choice of $x^{0}_{\alpha}$
together with (\ref{1747}) we finally get, after  dividing the
previous inequality by $\delta^{2}$ and letting $\delta\go 0$, that
\begin{eqnarray*}
W_{\{\e_{n}\}}\left( x^{0}_{\alpha},\overline{\xi} \right) &\leq  &
\frac{1}{2 T^{2}} \int_{(0,T)^{2} \times I} \hspace{-0.1cm} W_{\rm
hom}(x^{0}_\a,y_3,y_\a;\overline
\xi+\nabla_\a \varphi(y)| \nabla_3 \varphi(y))\, dy\\
&\leq  & \overline W_{\rm hom}(x^{0}_\a;\overline \xi) + \rho
\end{eqnarray*}
\noindent and the result follows by letting $\rho\go 0$.
$\hfill\blacksquare$\\

\noindent{\it Proof of Theorem \ref{jf-m3}.} As a consequence of
Lemmas \ref{3d2d-general-ineq>} and \ref{3d2d-general-ineq<}, we
have $\overline W_{\rm hom}(x_\a;\overline
\xi)=W_{\{\e_n\}}(x_\a;\overline \xi)$ for all $x_\a \in L$ and
all $\overline \xi \in \Qb^{3 \times 2}$. Since $\overline W_{\rm
hom}$ and $W_{\{\e_n\}}$ are Carath\'eodory functions,  it follows
that the equality holds for all $\overline \xi \in \Rb^{3 \times
2}$ and a.e. $x_\a \in \o$. Therefore, we have
$\W_{\{\e_n\}}(u;A)=\W_{\rm hom}(u;A)$ for all $A \in \mathcal
A(\o)$ and all $u \in W^{1,p}(A;\Rb^3)$. Since the result does not
depend upon the specific choice of the subsequence, we conclude
thanks to Proposition 8.3 in Dal Maso \cite{DM} that the whole
sequence $\W_\e(\cdot\,;A)$ $\G(L^p(A \times I))$-converges to
$\W_{\rm hom}(\cdot\,;A)$. Taking $A=\o$ we conclude the proof of
Theorem \ref{jf-m3}.
$\hfill\blacksquare$\\

To conclude, let us state an interesting consequence of Theorem
\ref{jf-m3}.

\begin{coro}Let $W: \O\times \Rb^3 \times \Rb^{3 \times 3} \to \Rb$ be a  function satisfying $(A_1)$, $(A_2)
$ and $(A_4)$, and such that $W(x,\cdot\,;\xi)$ is $Q$-periodic
for all $\xi \in \Rb^{3 \times 3}$ and a.e.\! $x \in \O$. Define
the functional $\W_\e : L^p(\O;\Rb^3) \to [0,+\infty]$ by
$$\W_{\e}(u):=\left\{\begin{array}{ll} \ds \int_{\O}
W\left(x,\frac{x}{\e};\nabla_\a  u(x) \Big|\frac{1}{\e} \nabla_3 u(x)\right) dx & \text{if } u\in W^{1,p}(\O;\Rb^{3}),\\
&\\
+\infty & \text{otherwise}.
\end{array}\right.$$
Then   the   $\G(L^p(\O))$-limit of the family  $\{\W_\e\}_{\e>0}$
is given by
 $\W_{\rm hom} : L^p(\O;\Rb^3) \to [0,+\infty]$ with
$$\mathcal \W_{\rm hom}(u): =
 \left\{\begin{array}{ll} \ds 2 \int_{\omega} \overline W_{\rm hom}(x_{\alpha},\nabla_\a u(x_\a))\, dx_\a &
  \text{if } u\in W^{1,p}(\omega;\Rb^{3}), \\&\\
+ \infty &\text{otherwise},
\end{array}\right.$$

\noindent where, for all $\overline \xi\in \Rb^{3\times 2}$ and
a.e.\! $x_{\alpha} \in \o$

\begin{eqnarray*}
\overline W_{\rm hom}(x_{\alpha},\overline \xi) & := & \lim_{T\go
+ \infty} \inf_{\phi} \Big\{ \frac{1}{2 T^{2}}\int_{(0,T)^{2}
\times I} W_{\rm hom}(x_{\alpha},y_{3},y_\a;\overline
\xi+\nabla_\a \phi(y)| \nabla_3 \phi(y)
)\, dy: \nonumber\\
&&\hspace{2.0cm}\phi \in W^{1,p}((0,T)^{2} \times I ;\Rb^{3}),
\quad \phi=0 \text{ on }\partial (0,T)^2 \times
I\Big\}\end{eqnarray*}

\noindent and, for all $y_\a \in \Rb^2$ and a.e.\! $x \in \O$

$$W_{\rm hom}(x, y_\a;\xi):=\lim_{T\go + \infty} \inf_{\phi} \left\{
\frac{1}{T^{3}}\int_{(0,T)^{3}}  W(x,y_\a,z_3;\xi+\nabla
\phi(z))\, dz: \; \phi \in W_{0}^{1,p}((0,T)^{3}
;\Rb^{3})\right\}.$$
\end{coro}

\noindent {\bf Acknowledgements:}  We would like to thank Irene
Fonseca, Gilles Francfort and Giovanni Leoni for their fruitful
comments and improvements, and Dag Lukkassen for useful
bibliographical suggestions. We  acknowledge the LPMTM at Paris 13
University and the Center for Nonlinear Analysis at Carnegie
Mellon University for their hospitality and support. The research
of M. Ba\'{\i}a was partially supported by Funda\c{c}\~{a}o para a
Ci\^{e}ncia e Tecnologia (Grant PRAXIS XXI SFRH
$\hspace{-0.1cm}\backslash\hspace{-0.1cm}$ BD
$\hspace{-0.1cm}\backslash\hspace{-0.1cm}$ 1174
$\hspace{-0.1cm}\backslash\hspace{-0.1cm}$ 2000), Fundo Social
Europeu, the Department of Mathematical Sciences of Carnegie
Mellon University and its Center for Nonlinear Analysis, and by
the MULTIMAT Network.


\vspace{0.5cm}

\begin{center}
\begin{small}

Jean-Fran\c{c}ois Babadjian\\
\textsc{L.P.M.T.M., Universit\'e Paris Nord, 93430, Villetaneuse, France}\\
\textit{E-mail address}: {\bf jfb@galilee.univ-paris13.fr}

\vspace{0.5cm}

Margarida Ba\'{\i}a\\
\textsc{Department of Mathematical Science, Carnegie Mellon University,}\\
\textsc{Pittsburgh, PA 15213, USA}\\
\textit{E-mail address}: {\bf mbaia@math.ist.utl.pt}

\end{small}

\end{center}

\end{document}